\documentclass[11pt, reqno, oneside]{amsart}
\usepackage[a4paper, total={6.7in, 9.5in}]{geometry}
\usepackage{amsmath,amsfonts,amsthm,amssymb,nicefrac}
\usepackage{graphicx,color,url}
\usepackage{xcolor}
\usepackage{accents}
\usepackage{enumerate}
\usepackage{enumitem}
\usepackage{comment,hyperref}
\includecomment{versionb}
\newtheorem{theorem}{Theorem}
\newtheorem{lemma}[theorem]{Lemma}
\newtheorem{corollary}[theorem]{Corollary}
\newtheorem{proposition}[theorem]{Proposition}
\theoremstyle{definition}
\newtheorem{remark}[theorem]{Remark}

\newtheorem{definition}[theorem]{Definition}
\newtheorem{claim}[theorem]{Claim}

\numberwithin{equation}{section}\numberwithin{theorem}{section}
\newcounter{stepctr}
{\end{list}}

\def\XXint#1#2#3{{\setbox0=\hbox{$#1{#2#3}{\int}$}
 \vcentre{\hbox{$#2#3$}}\kern-.5\wd0}}

\providecommand{\abs}[1]{\lvert#1\rvert}

\newcommand{\mbb}[1]{\mathbb{#1}}
\newcommand{\circo}{\accentset{\circ}}

\newcommand{\R}{\mbb{R}}

\newcommand{\e}{\varepsilon}

\newcommand\simtimes{\mathbin{%
    \stackrel{\sim}{\smash{\times}\rule{0pt}{0.8ex}}%
    }}

\DeclareMathOperator{\tr}{tr}

\DeclareMathOperator{\vol}{vol}

\DeclareMathOperator{\Rm}{Rm}
\DeclareMathOperator{\dist}{dist}

\DeclareMathOperator{\pl}{pl}

\begin{document}

\title[High codimension MCF with surgery]{High codimension mean curvature flow with surgery}

\author{Stephen Lynch}
\address{Department of Mathematics, King's College London, Strand, London, WC2R 2LS, UK}
\email{stephen.lynch@kcl.ac.uk}

\author{Huy The Nguyen}
\address{School of Mathematical Sciences\\
	Queen Mary University of London\\
	Mile End Road\\
	London E1 4NS}
\email{h.nguyen@qmul.ac.uk}


\maketitle

\begin{abstract}
We construct a mean curvature flow with surgery for submanifolds of arbitrary codimension. The theory applies to closed submanifolds satisfying a natural quadratic pinching condition, which serves as the high-codimension analogue of $2$-convexity and is preserved under the flow in dimensions $n \geq 8$. Our results therefore are in line with the current state-of-the-art in codimension one (where at present $2$-convexity is required for surgery). Central to our analysis is a collection of new a priori estimates for the second fundamental form, uniform across surgeries, which yield a precise description of high-curvature regions and permit controlled surgeries. This provides the first notion of mean curvature flow through singularities with topological control in higher codimensions. As a consequence we obtain a sharp classification: Every closed quadratically $2$-convexity submanifold is diffeomorphic either to $\mathbb{S}^n$ or to a finite connected sum of $\mathbb{S}^{n-1}$-bundles over $\mathbb{S}^1$.

\end{abstract}

\tableofcontents

\section{Introduction}

Over the last 50 years, curvature flows have become central objects of study in geometry. Among these, the mean curvature flow is one of the most important: It is the natural analogue of the heat equation for submanifolds and, in principle, offers a powerful mechanism for evolving extrinsic geometries into canonical forms. However solutions are subject to finite-time singularity formation. A profound problem, with many potential applications, is to determine the structure of these singularities and reveal the geometric information they encode. 

For hypersurfaces, the deep work of Huisken–Sinestrari \cite{HuSi09} developed a mean curvature flow with surgery for $2$-convex hypersurfaces, demonstrating that singularities can be controlled, excised, and replaced in a manner that allows the flow to continue. Their construction produced a decomposition of arbitrary $2$-convex hypersurfaces into a collection of model geometries. This had striking topological consequences, including a proof of an important special case of the Schoenflies conjecture: Every $2$-convex embedding of $\mathbb{S}^3$ into $\mathbb{R}^4$ bounds a standard $4$-ball. These ideas parallel the role of the Ricci flow with surgery in the Hamilton--Perelman program, which resolved the geometrization conjecture for $3$-manifolds.


In contrast, the high-codimension mean curvature flow is far more intricate, and until now no analogue of the Huisken–Sinestrari theory has been available. Various new phenomena arise for flows of codimension greater than one: the normal bundle typically has nonzero curvature, the second fundamental form has additional components orthogonal to the mean curvature direction, embeddedness fails to be preserved, and the evolution equations for the second fundamental form contain an array of additional reaction terms which are not easily understood. For these reasons even the qualitative structure of singularities in high codimension has remained largely elusive.\footnote{It is worth mentioning that there is a huge array of possible singularity models in higher codimensions. For example, every minimal submanifold of a sphere of any dimension generates a self-similarly shrinking mean curvature flow in Euclidean space.}

In this work we overcome these long-standing difficulties and construct, for the first time, a mean curvature flow with surgery for submanifolds of arbitrary codimension. This is the first notion of mean curvature flow through singularities with topological control in codimensions greater than one.\footnote{There exist several weak notions of mean curvature flow \cite{Brakke1978, Evans1991, Chen1991, Ilmanen1994, Ambrosio1996}. These do not allow for topological control through singularities.} Our results apply to closed immersions $F : \mathcal M^n \to \mathbb{R}^{n+m}$ satisfying the \emph{quadratic $2$-convexity} condition 
    \[|A|^2 < \frac{1}{n-2} |H|^2,\]
where $A$ is the second fundamental form and $H$ is the mean curvature vector. This condition is the natural high-codimension analogue of the $2$-convexity assumption made in \cite{HuSi09}. It is known, by work of Andrews–Baker \cite{Andrews2010}, to be preserved under the flow in dimensions $n \geq 8$. 


\begin{theorem}\label{thm_main}
Let $n \geq 8$. If $F : \mathcal M^n \to \mathbb{R}^{n+m}$ satisfies $|A|^2 < \frac{1}{n-2}|H|^2$ then there exists a mean curvature flow with surgery starting from $F$ which terminates after finitely many steps. 
\end{theorem}

Under the flow with surgery, the submanifold evolves by the smooth mean curvature flow until its maximum curvature becomes large. We then stop the flow, carefully remove regions of large curvature, restart the flow and repeat. Theorem~\ref{thm_main} asserts that we can iterate this process until it terminates---that is, after finitely many surgeries every remaining component of the submanifold belongs to a finite list of model geometries. As a consequence we obtain the following complete classification of quadratically $2$-convex submanifolds in every codimension. 

\begin{corollary}\label{cor_topological conclusion}
If $n \geq 8$ and $\mathcal M^n$ admits a smooth immersion into some $\mathbb{R}^{n+m}$ satisfying $|A|^2 < \frac{1}{n-2} |H|^2$, then it is diffeomorphic either to $\mathbb{S}^n$ or to a finite connected sum of $\mathbb{S}^{n-1}$-bundles over $\mathbb{S}^1$. 
\end{corollary}

\begin{remark}
Corollary~\ref{cor_topological conclusion} is sharp in the sense that, if instead of quadratic $2$-convexity we assume $|A|^2 < c|H|^2$ with $c>\frac{1}{n-2}$, further topologies appear. Indeed, the codimension-2 submanifolds $\mathbb{S}^{n-2}(r)\times\mathbb{S}^2(1) \subset \mathbb R^{n-1}\times\mathbb R^3$ are such that $|A|^2/|H|^2 \searrow \frac{1}{n-2}$ as $r\to 0$. 
\end{remark}

To contextualise Corollary~\ref{cor_topological conclusion}, we recall that by Nash's embedding theorem \cite{Nash1956} any Riemannian manifold can be isometrically embedded in a Euclidean space of sufficiently large dimension. This means the topology and geometry of high-codimension submanifolds are unrestricted, in stark contrast with the codimension-one case. Corollary~\ref{cor_topological conclusion} is the first local-to-global classification result in higher codimensions which only assumes a one-sided curvature condition and is sufficiently general to include nontrivial topologies. 

Our classification of quadratically $2$-convex submanifolds exhibits a richer variety of topologies than its counterpart in codimension one (see \cite{HuSi09}). These arise from nontrivial sphere bundles. Up to diffeomorphism there are precisely two $\mathbb{S}^{n-1}$-bundles over $\mathbb{S}^1$, the product $\mathbb{S}^{n-1} \times \mathbb{S}^1$ and a single nontrivial example which we denote by $\mathbb{S}^{n-1} \simtimes \mathbb{S}^1$. The conclusion of Corollary~\ref{cor_topological conclusion} can be restated as follows: every manifold satisfying the hypotheses there is diffeomorphic either to $\mathbb{S}^n$ or to a finite connected sum of copies of $\mathbb{S}^{n-1} \times \mathbb{S}^1$ and $\mathbb{S}^{n-1} \simtimes \mathbb{S}^1$. The space $\mathbb{S}^{n-1} \simtimes \mathbb{S}^1$ does not appear in \cite{HuSi09} because it cannot be immersed in $\mathbb{R}^{n+1}$ with nowhere vanishing mean curvature (it is not orientable).  

To conclude this first part of the introduction we mention some important results related to our own. Following \cite{HuSi09}, flows with surgery have been developed for mean-convex surfaces of codimension one in \cite{BrHu16, Brendle2018, Haslhofer2017a}. These works rely heavily on noncollapsing estimates which are not available in higher codimensions. A Ricci flow with surgery was constructed by Perelman  \cite{Perelman2002, Perelman2003} in dimension~3, and under natural curvature conditions in higher dimensions by Hamilton~\cite{Hamilton1997} and Brendle~\cite{Brendle2016, Brendle2017}. 

\subsection{Key ideas and new analytic tools} At the heart of our construction is a detailed analysis of high-curvature regions under quadratic $2$-convexity, based on new a priori estimates for the second fundamental form that are uniform across surgeries. The main analytic advances on which we build our theory are:

\begin{enumerate}
    \item A scale-invariant, pointwise \emph{gradient estimate} for the second fundamental form (Theorem~\ref{thm_grad est}).\footnote{This estimate first appeared in the preprint \cite{Nguyen2018a}, which has not been published and is superseded by the present paper. It was used in \cite{Naff_singularity_models, Naff_canonical_neighbourhood} to prove a canonical neighbourhood theorem, applicable to smooth flows, but not to the flow with surgery beyond the first surgery time.} This permits controlled comparison of curvature at nearby points, even near surgery regions. 
    \item A \emph{planarity estimate} for flows with surgery (Theorem~\ref{thm_planarity surgeries}), extending Naff’s result \cite{Naff_planarity} for smooth flows. It asserts that, at points of large curvature, the flow becomes asymptotically codimension one: the components of the second fundamental form orthogonal to the mean curvature direction are of lower order.
    \item A \emph{planarity improvement theorem} (Theorem~\ref{thm_planarity improvement}), ensuring that just before each surgery, the submanifold becomes close to a hypersurface with much stronger quantitative control than is provided by the planarity estimate alone. This is crucial for ensuring that the planarity estimate is uniform globally in time, since it deteriorates slightly under each surgery. 
    \item \emph{Cylindrical estimates} (Theorem~\ref{thm_cylindrical via Stampacchia surgery}) which quantify how close high-curvature regions are to a cylinder $\mathbb{S}^{n-1} \times \mathbb{R}$ after rescaling.
\end{enumerate}

These estimates allow us to prove that, close to a potential singular time, either the entire submanifold is a positively curved immersion of a sphere or else it contains large neck regions where surgery can be performed. This is the content of our \emph{neck detection lemma} (Lemma~\ref{lem_NDL}). In order to ensure that surgery removes all regions of large curvature we prove a \emph{neck continuation theorem} (Theorem~\ref{thm_neckcontinuation}) which characterises the geometry of the submanifold as we travel out from a neck: either the curvature drops by a fixed factor or else the neck eventually closes up in a positively curved cap. A similar analysis was conducted for $2$-convex hypersurfaces in \cite{HuSi09}, but in higher codimensions the situation is more subtle, because the submanifold may not lie close to a hypersurface after we leave the neck. 

\subsection{Dimension considerations} The restriction to $n \geq 8$ is natural. Andrews–Baker established preservation of quadratic $2$-convexity precisely in this range \cite{Andrews2010}, and our analysis relies critically on this monotonicity. In dimensions below eight, the reaction terms in the evolution equation of $|A|^2 - \frac{1}{n-2}|H|^2$ do not appear to have a sign, and we conjecture that quadratic $2$-convexity is not preserved under the flow in these dimensions.\footnote{The aforementioned result of Andrews--Baker from \cite{Andrews2010} asserts that $|A|^2 < c |H|^2$ is preserved under the flow if $c \leq \frac{4}{3n}$. This has not yet been improved upon, see however \cite{Baker2017}.}

Furthermore, explicit constructions in low dimensions demonstrate that surgeries of the type we perform above cannot exist in general. Indeed, for $n=3$ quadratic $2$-convexity is equivalent to having positive scalar curvature (since $|H|^2 - |A|^2$ is equal to the scalar curvature by the Gauss equations), but the set of $3$-manifolds admitting a metric with positive scalar curvature includes many examples not listed in Corollary~\ref{cor_topological conclusion}, such as $\mathbb{RP}^2 \times \mathbb{S}^1$.

All of the analysis we have discussed so far does apply in dimensions $n \in \{5,6,7\}$, but only if we assume a pinching condition which is somewhat stronger than quadratic $2$-convexity and does not seem to be optimal. 




\subsection{Other quadratic pinching conditions} In \cite{Andrews2010} Andrews--Baker showed that, for $n \geq 4$ and in every codimension, a mean curvature flow which satisfies $|A|^2 < \frac{1}{n-1}|H|^2$ at the initial time will contract to a point, becoming asymptotically round in the process.\footnote{For $n \in \{2,3\}$ they draw the same conclusion under the more restrictive hypothesis  $|A|^2 < \frac{4}{3n} |H|^2$.} In particular, this implies that the initial submanifold must have been diffeomorphic to $\mathbb{S}^n$. The pinching condition $|A|^2 < \frac{1}{n-1}|H|^2$ is a quadratic analogoue of convexity, so this theorem generalises to higher codimensions Huisken's seminal work \cite{Hu84} on convex hypersurface flows. 

If one instead assumes $|A|^2 < c |H|^2$ for some $c > \frac{1}{n-1}$ then the mean curvature flow will in general form local singularities rather than contracting to a point. The condition $|A|^2 < \frac{1}{n-2}|H|^2$ which we study in this paper is the most general which allows for neckpinches (singularities where some blow-up is $\mathbb{S}^{n-1}\times\mathbb{R}^1$) but rules out bubblesheet singularities (where some blow-up is $\mathbb{S}^{n-2} \times \mathbb{R}^2$). Future work should aim to understand flows satisfying the quadratic $k$-convexity condition $|A|^2 < \frac{1}{n-k}|H|^2$, which is preserved in dimensions $n \geq 4k$ by \cite{Andrews2010}.

\subsection{Outline} Section~\ref{sec_prelim} contains preliminaries, including our definition of the mean curvature flow with surgery. In Section~\ref{sec_standard surgery} we construct a standard surgery procedure for removing the central portion of a neck and replacing it with two caps. In Section~\ref{sec_estimates} we prove the pointwise gradient estimates for the second fundamental form. In Section~\ref{sec_planarity} we prove the planarity estimate for flows with surgery and the planarity improvement theorem. The cylindrical estimates are established in Section~\ref{sec_cylindrical estimates}. With all of these a priori estimates in place, in Section~\ref{sec_ND} we prove the neck detection lemma and detail its consequences. Finally in Section~\ref{sec_surgeries} we prove the neck continuation theorem and fix our surgery algorithm depending on initial data, before concluding with the proofs of Theorem~\ref{thm_main} and Corollary~\ref{cor_topological conclusion}.

\subsection*{Acknowledgements}
This research was funded by the EPSRC grant EP/S012907/1. The authors are grateful to Gerhard~Huisken, Mat~Langford, Tang-Kai~Lee and Keaton~Naff for many long conversations about the high-codimension mean curvature flow from which this work has benefited. 


\section{Preliminaries}\label{sec_prelim}
Let $F: \mathcal{M}\times[0,T)\rightarrow \R^{n+m}$ be a solution of the mean curvature flow with closed, smoothly immersed timeslices $\mathcal{M}_t = F(\mathcal{M},t)$. We denote by $g = g(t)$ the induced metric, by $\mu = \mu(t)$ the induced volume measure and by $A(X,Y) := (D_X Y)^\perp$ the (vector-valued) second fundamental form. The mean curvature vector is $H = \tr_g A$. Given locally defined frames $\{X_i\}$ and $\{\nu_\alpha\}$ for the tangent and normal bundles we write $A_{ij} = A(X_i, X_j)$ and $A_{ij}^\alpha = \langle A_{ij}, \nu_\alpha\rangle$ for the components of $A$. In every codimension, whenever $|H| > 0$, we may single out a principal normal direction $\nu :=H/|H|$. (In codimension one, $\nu$ is the inward unit normal.) Projecting the second fundamental form onto this principal direction yields the bilinear form $h := \langle A, \nu\rangle$. The remaining components of the second fundamental form we denote by $A^-= A- h \nu$. We recall from \cite{Naff_planarity} that, as consequence of the Codazzi equations, with respect to any orthonormal frame for the tangent space we have
    \begin{equation}\label{eq_principal torsion}
    (|H| g_{ik} - h_{ik})\nabla_k \nu = \nabla_k A^-_{ki} - \langle \nabla_k A^-_{ki}, \nu\rangle \nu.
    \end{equation}
This simple identity will play a key role at multiple points in our analysis. 

If $\mathcal{M}_t$ evolves by the mean curvature flow then we have the evolution equations\footnote{Both $\langle U ,V \rangle$ and $U \cdot V$ denote the Euclidean inner product of $U,V \in \mathbb{R}^{n+m}$.}
	\begin{align*}
	\partial_t g &= -2 H \cdot A,\qquad
\partial_t d \mu= - |H|^2 d \mu.
	\end{align*}
We have induced connections on the tangent and normal bundles of each $\mathcal{M}_t$, given by taking projections of the ambient Euclidean connection:
    \[\nabla_X Y := (D_X Y)^\top, \qquad \nabla_X N := (D_X N)^\perp\]
for $X, Y$ tangent and $N$ normal to $\mathcal{M}_t$. To differentiate vector fields in time we use the connection $\nabla_{t}$ constructed in \cite{Andrews2010} which acts by
    \[g(\nabla_t X, Y) = \partial_t X \cdot Y - H \cdot A(X,Y), \qquad \nabla_{t} N = (\partial_t N)^\perp\]
for $X, Y$ tangent and $N$ normal to $\mathcal{M}_t$. Note that we have $\nabla_{t} g = 0$. These connections extend to duals and tensor products in the usual manner.

We have Simons' identity 
    \begin{align*}
    \Delta A_{ij}&=\nabla_i\nabla_jH+H\cdot A_{ip}A_{pj}-A_{ij}\cdot A_{pq}A_{pq}\\
    &\qquad +2A_{jq}\cdot A_{ip}A_{pq}-A_{iq}\cdot A_{qp}A_{pj}-A_{jq}\cdot A_{qp}A_{pi},
    \end{align*}
which is used to derive the following evolution equation for $A$ (see for example \cite{Andrews2010}): 
    \begin{align*}
    \nabla_{t} A_{ij} &= \Delta A_{ij}+A_{ij}\cdot A_{pq}A_{pq}+A_{iq}\cdot A_{qp}A_{pj}\\
    &\qquad +A_{jq}\cdot A_{qp}A_{pi}-2A_{ip}\cdot A_{jq}A_{pq}.
    \end{align*}
The evolution equation for the mean curvature vector is found by taking the trace:
    \begin{align*}
    \nabla_{t} H=\Delta H+H\cdot A_{pq}A_{pq}.
    \end{align*}
Fixing in addition a local orthonormal frame for the normal bundle, the evolution equations for $|A|^2$ and $|H|^2$ may be expressed as
    \begin{align}
    \partial_t \abs{A}^2= \Delta\abs{A}^2-2\abs{\nabla A}^2+2 \sum_{\alpha, \beta}\bigg(\sum_{i,j} A_{ij\alpha}A_{ij\beta}\bigg)^2& \label{eqn:A2} \\
    \qquad +2 \sum_{i,j,\alpha,\beta}\bigg( \sum_p A_{ip\alpha}A_{jp\beta}-A_{jp\alpha}A_{ip\beta}\bigg)^2&, \notag \\
    \partial_t\abs{H}^2= \Delta\abs{H}^2-2\abs{\nabla H}^2+2\sum_{i,j}\bigg( \sum_{\alpha}H_{\alpha}A_{ij\alpha}\bigg)^2&. \label{eqn:H2}
    \end{align}
The last term in \eqref{eqn:A2} is the squared length of the normal curvature, which we denote by $\abs{\Rm^{\perp}}^2$. For convenience we label the reaction terms of the above evolution equations by
    \[R_1 = \sum_{\alpha, \beta}\bigg(\!\sum_{i,j}A_{ij\alpha}A_{ij\beta}\!\bigg)^2+\abs{\Rm^{\perp}}^2, \qquad R_2 = \sum_{i,j}\!\bigg(\!\sum_{\alpha}H_{\alpha}A_{ij\alpha}\bigg)^2. \]

\subsection{Preservation of pinching}
We consider $Q := |A|^2+a-c |H|^2$ where $a$ and $c$ are positive constants. Combining the evolution equations for $|A|^2$ and $|H|^2$ yields
    \begin{align}
    \label{eqn_pinchpres}
    \partial_t Q&= \Delta Q-2 ( |\nabla A |^2-c |\nabla H|^2)+2 (R_1 - cR_2).
    \end{align}
We have the following Kato-type inequality, which is a consequence of the Codazzi equations. This is proven in \cite{Andrews2010} (as in \cite{Hu84}) and shows that the gradient terms in \eqref{eqn_pinchpres} are nonpositive if $c \leq \frac{3}{n+2}$.

\begin{lemma}\label{lem_Kato}
For any $n$-submanifold of Euclidean space we have $|\nabla A|^2\geq \frac{3}{n+2}|\nabla H|^2$. 
\end{lemma}
 
If $c$ satisfies the stricter inequality $c \leq \frac{4}{3n}$ then we also have $R_1-cR_2 \leq 0$ at any point where $Q \leq 0$ (see \cite[Theorem~2]{Andrews2010}), which means the pinching condition $Q \leq 0$ is preserved, by the parabolic maximum principle.

\begin{lemma}\label{lem_pinching}
Let $F:\mathcal{M}^n \times[0,T)\rightarrow \R^{n+m}$ be a solution to the mean curvature flow such that $\mathcal{M}_0$ satisfies $|A|^2+a \leq c |H|^2$ for some $a > 0$ and $c \leq \frac{4}{3n}$. This condition is then preserved for all $t \in [0,T)$. 
\end{lemma}

We will require a slight refinement of Lemma~\ref{lem_pinching}. Following the computations in \cite[Theorem~2]{Andrews2010}, but keeping track of $Q$ and $a$, one obtains:

\begin{lemma}\label{lem_pinching refinement}
Let $F:\mathcal{M}^n \times[0,T)\rightarrow \R^{n+m}$ be a solution to the mean curvature flow with $|H| > 0$. Then for the quadratic quantity $Q := |A|^2 + a - c|H|^2$ we have 
    \begin{align}\label{eqn_pinchpres2}
    \begin{split}
    \partial_t Q&\leq \Delta Q-2 ( |\nabla A |^2-c |\nabla H|^2)+2 |h|^2 Q +\frac{2}{n}\frac{1}{c-\nicefrac{1}{n}}|A^-|^2 Q \\
    &\qquad -2 a |h|^2-\frac{2a}{n}\frac{1}{c-\nicefrac{1}{n}}| A^-|^2 \\
    &\qquad + \left(6-\frac{2}{n (c-\nicefrac{1}{n})} \right) |\circo h|^2 | A^-|^2+\left(3-\frac{2}{n (c-\nicefrac{1}{n})} \right)|A^-|^4.
    \end{split}
    \end{align}
\end{lemma}

\subsection{Surgery class} Fix $n \geq 5$. Given positive parameters $R$ and $\alpha = (\alpha_1,\alpha_2,\alpha_3)$, we define $\mathcal{C}_{n,m}(R,\alpha)$ to be the class of closed immersions $F: \mathcal M^n \rightarrow \R^{n+m}$ satisfying:
    \begin{enumerate}[label=(\roman*)]
    \item $|A|^2 + \alpha_2 R^{-2}  \leq (c_n- \alpha_1) |H|^2$ where $c_n := \frac{1}{n-2}$ in dimensions $n \geq 8$ and $c_n := \frac{3(n+2)}{2n(n+2)}$ in dimensions $5 \leq n \leq 7$, 
    \item $\mu(\mathcal M^n)\leq \alpha_3 R^n$. 
    \end{enumerate}
The parameters $\alpha$ are scale-invariant: If $F\in \mathcal{C}_{n,m}(R,\alpha)$ then $rF \in \mathcal{C}_{n,m}(r R,\alpha)$. In a slight abuse of terminology we also say that an immersed submanifold is in the class $\mathcal C_{n,m}(R,\alpha)$ if it admits a parameterization $F$ with the above properties.

For each $R$ and $\alpha$ the class $\mathcal{C}_{n,m}(R,\alpha)$ is invariant under the mean curvature flow. This follows immediately from Lemma~\ref{lem_pinching} and $\frac{d}{dt}d\mu = -|H|^2d\mu$. Our flow with surgery will be constructed so that it too preserves this class. 

\subsection{Mean curvature flow with surgery}
We now define the mean curvature flow with surgery in arbitrary codimensions. Our definition is similar to that in \cite{HuSi09}. Notice however that since submanifolds in $\mathcal C_{n,m}(R,\alpha)$ need not be orientable the twisted bundle $\mathbb{S}^{n-1} \simtimes \mathbb{S}^1$ necessarily appears in our definition. 

\begin{definition} The \emph{mean curvature flow with surgery} is determined by an algorithm that assigns to each initial submanifold $F_0: \mathcal{M}_1 \rightarrow \R^{n+m}$ in the class $\mathcal{C}_{n,m}(R,\alpha)$:
\begin{itemize}
\item a sequence of intervals $[0,T_1]$, $[T_1, T_2]$, $[T_2, T_3],\cdots, [T_{N-1}, T_N]$,
\item a sequence of manifolds $\mathcal{M}_i, 1 \leq i \leq N$, and
\item a sequence of smooth mean curvature flows $F_t^i: \mathcal{M}_i\rightarrow \R^{n+m}, t\in [ T_{i-1}, T_i]$.
\end{itemize}
Moreover, these smooth mean curvature flows are such that:
\begin{itemize}
\item the initial submanifold for the flow $F^1_t$ is given by $F_0: \mathcal{M}_1\rightarrow \R^{n+m}$, and
\item for $2\leq i \leq N$, the initial submanifold for $F^i_t:\mathcal{M}_i \rightarrow \R^{n+m}$ on $[T_{i-1},T_i]$ is obtained from $F_{T_{i-1}}^{i-1}$ by the following two-step procedure. First, using the standard surgery procedure defined in Section \ref{sec_standard surgery}, a submanifold $ \hat{F}_{T_{i-1}}^{i-1}: \mathcal{M}_i \rightarrow \R^{n+m}$ is obtained from $F_{T_{i-1}}^{i-1}: \mathcal{M}_{i-1}\rightarrow \R^{n+m}$ by replacing finitely many disjoint necks each with two caps. Then, finitely many connected components are recognised as being diffeomorphic to $ \mathbb{S}^n$, $ \mathbb{S}^{n-1} \times \mbb S^1$ or $\mathbb{S}^{n-1} \simtimes \mathbb{S}^1$ and are discarded. After these two steps are completed, the resulting submanifold is taken as the initial submanifold for the smooth mean curvature flow $F^i_t$ on $[T_{i-1},T_i]$.
\end{itemize}
The mean curvature flow with surgery \emph{terminates after finitely many steps} if either:
\begin{itemize}
\item at time $T_N$ all connected components of $F_{t}^N$ are recognised as being diffeomorphic to $ \mathbb{S}^n$, $ \mathbb{S}^{n-1} \times \mbb S^1$ or $\mathbb{S}^{n-1} \simtimes \mathbb{S}^1$, or
\item in the second step above all components of the surgically modified submanifold $\hat{F}_{T_N}^N$ are recognised as being diffeomorphic to $ \mathbb{S}^n$, $\mathbb{S}^{n-1} \times \mbb S^1$ or $\mathbb{S}^{n-1} \simtimes \mathbb{S}^1$.
\end{itemize}
 \end{definition}
A priori, there is no reason for the surgery algorithm to terminate after finitely many steps. For initial submanifolds in the class $\mathcal C_{n,m}(R,\alpha)$ we will ensure this by choosing the surgery times and scales so that each surgery removes a portion of the submanifold of fixed area. In particular, we will fix constants $\omega_1, \omega_2, \omega_3 > 1$ depending only on $n, m, \alpha$, choose curvature thresholds $H_1, H_2, H_3$ such that
	\begin{align*}
	H_1\geq \omega_1 R^{-1}, \quad H_2 = \omega_2 H_1, \quad H_3 = \omega_3 H_2,
	\end{align*}
and construct the mean curvature flow with surgery so that:
\begin{itemize}
\item Each surgery time $T_i, 1 \leq i \leq N$, is the first time in $[T_{i-1}, T_i]$ where $|H|$ has a maximum value of $H_3$.
\item For each $2 \leq i \leq N$ the maximum mean curvature of $F^i_t$ at its initial time $t = T_{i-1}$ is at most $H_2$. In particular, $|H|$ is bounded from above by $H_3$ at every time.
\item All surgeries are performed in regions where $\frac{10}{11} H_1\leq |H|\leq \frac{11}{10}H_1$.
\end{itemize}

\subsection{Necks}

We introduce the notion of an $(\e,k,L)$-neck region in a submanifold of higher codimension. A related intrinsic notion was introduced in \cite{Hamilton1997}, and a related notion for hypersurfaces was introduced in \cite{HuSi09}.

\begin{definition}\label{def_neck}
Let $\mathcal M$ be an $n$-submanifold of $\mathbb{R}^{n+m}$. We say that $\mathcal{M}$ is an \emph{$(\e,k)$-neck at scale $r >0$} if there exists $\Omega\subset\mathbb{S}^{n-1}\times\mathbb{R}$ such that $r^{-1}\mathcal{M}$ is a cylindrical graph over $\Omega$ in the following sense: Up to a choice of ambient Euclidean coordinates, $r^{-1} \mathcal{M}$ agrees with 
    \[\{(\omega+u(\omega,z)\omega,y(\omega,z),z):(\omega,z)\in\Omega\subset\mathbb{S}^{n-1}\times\mathbb{R}\},\]
where the maps $u:\Omega\to\mathbb{R}$ and $y:\Omega\to\mathbb{R}^{m-1}$ satisfy $\|u\|_{C^{k+2}}+\|y\|_{C^{k+2}} \leq \e$. Most often we take $\Omega$ to be a region of the form $\mathbb{S}^{n-1} \times [a,b]$ where $b - a$ is large. When there is no need to specify the scale $r$ we simply refer to $\mathcal{M}$ as an \emph{$(\e,k)$-neck}. In addition, $\mathcal{M}$ is an \emph{$(\e,k,L)$-neck} about $p\in\mathcal{M}$ if the intrinsic ball $\mathcal B(p,Lr) \subset \mathcal{M}$ is an $(\e,k)$-neck at scale $r$.
\end{definition}

For an $(\e, k)$-neck, if $\e$ is small then the second fundamental form and its derivatives (up to order $k$) are close to those of the standard cylinder after rescaling (cf. the notion of an `extrinsic curvature neck' used in \cite{HuSi09}). In particular, as $\e \to 0$ we have 
    \[|H| = (n-1 + O(\e)) r^{-1}, \qquad |h|^2 - \frac{1}{n-1} |H|^2 = O(\e)r^{-2}, \qquad |A^-|^2 = O(\e)r^{-2}\]
and $|\nabla^\ell A|^2 = O(\e)r^{-2(\ell+1)}$, $1 \leq \ell \leq k$,
uniformly for any fixed $n$ and $m$.

As in \cite[Section~3]{HuSi09}, our notion of an $(\e,k)$-neck is stronger than the related intrinsic notion developed in \cite{Hamilton1997} (see in particular \cite[Proposition~3.3]{HuSi09}). Therefore, we may utilise Hamilton's normal-form parameterization for necks based on CMC-spheres and harmonic mappings. Our discussion closely follows Section~3 of \cite{HuSi09}. For an $n$-dimensional submanifold $\mathcal{M} \subset \mathbb{R}^{n+m}$ and a local diffeomorphism $N:\mathbb S^{n-1}\times[a,b] \to \mathcal{M}$ we write $\Sigma_z = N(\mbb S^{n-1}\times\{z\})$. The average radius $r: [a,b] \to \mathbb{R}$ is then defined by requiring that $|\Sigma_z|_g = \sigma_{n-1}r(z)^{n-1}$, where $\sigma_{n-1}$ is the area of $\mathbb{S}^{n-1}$. The following notion of an $(\e, k)$-cylindrical submanifold neck is our replacement for the $(\varepsilon,k)$-cylindrical hypersurface necks in \cite{HuSi09} (see Definition~3.9 therein). 

\begin{definition}[cf.\ Definition~3.9 in \cite{HuSi09}] A local diffeomorphism $N: \mbb{S}^{n-1}\times[a,b]\to \mathcal{M}$ is an \emph{$(\e, k)$-cylindrical submanifold neck} if it satisfies the conditions of \cite[Definition~3.7]{HuSi09} and in addition
    \[\left||h|^2 - \frac{1}{n-1}|H|^2\right| \leq \e r(z)^{-2}, \qquad |A^-|^2 \leq \e r(z)^{-2}, \qquad |\nabla^\ell A|^2 \leq \e r(z)^{-2(\ell+1)}\]
for every $1 \leq \ell \leq k$ and $z \in [a,b]$.  We refer to $b-a$ as the length of the neck. Note that the length of a neck is invariant under scaling.
\end{definition}

\begin{definition}
A local diffeomorphism $N: [a,b] \times S^{n-1} \to \mathcal{M}$ is \emph{normal} if it satisfies the conditions of \cite[Definition~3.8]{HuSi09}. An $(\e,k)$-cylindrical submanifold neck $N$ is a \emph{maximal normal} $(\e,k)$-cylindrical submanifold neck if $N$ is normal and if whenever $N^*$ is another such neck with $N=N^*\circ G$ for some diffeomorphism $G$ then $G$ is surjective.
\end{definition}

\begin{remark}
Given a normal $(\varepsilon, k)$-cylindrical submanifold neck of length $b-a \geq 2\delta$, for any $(\varepsilon', k')$ we can choose $\varepsilon$ small and $k$ large so that every point with $z \in [b+\delta, a-\delta]$ lies at the centre of an $(\varepsilon',k',\delta/2)$-neck (in the graphical sense of Definition~\ref{def_neck}).
\end{remark}

Results from \cite{Hamilton1997} provide existence, uniqueness and overlapping properties for normal necks. In fact, \cite[Theorem~3.12]{HuSi09} carries over to our setting, provided we replace each instance of `$(\e, k)$-cylindrical hypersurface neck' in the statement of that theorem with `$(\e,k)$-cylindrical submanifold neck'. The only important difference is that in codimension $m \geq 2$, if $\mathcal{M}$ is compact and covered by a neck then it could be diffeomorphic to $\mathbb{S}^{n-1} \times \mathbb{S}^1$ or to $\mbb{S}^{n-1} \simtimes \mbb{S}^1$. The latter does not occur as an immersed $2$-convex hypersurface in $\mathbb{R}^{n+1}$ since it is not orientable. 

\begin{theorem}\label{neck overlapping properties}
For any $\delta > 0$ we can choose $\varepsilon$ and $k$ so that if $N : \mathbb{S}^{n-1}\times[a,b] \to \mathcal{M}$ is an $(\varepsilon,k)$-cylindrical submanifold neck with $b-a \geq 3\delta$ then we can find a normal neck $N^\ast$ and a diffeomorphism $G$ of the domain cylinder of $N^\ast$ onto a region in the domain cylinder of $N$ containing all points at least $\delta$ from the ends such that $N^\ast = N \circ G$. The normal neck $N^\ast$ is contained in a maximal normal $(\varepsilon,k)$-cylindrical submanifold neck unless the target $\mathcal{M}$ is diffeomorphic either to $\mathbb{S}^{n-1} \times \mathbb{S}^1$ or $\mathbb{S}^{n-1} \simtimes \mathbb{S}^1$.
\end{theorem}
\begin{proof}
This follows from Theorem~2.5 in Section~3 of \cite{Hamilton1997} (cf.\ \cite[Theorem~3.12]{HuSi09}). By that theorem $N^\ast$ is contained in a maximal normal neck unless $\mathcal M$ is diffeomorphic to a quotient of $\mathbb{S}^{n-1}\times\mathbb{R}$. Our necks have approximately umbilic and hence spherical cross-sections so the only allowed quotients are $\mathbb{S}^{n-1} \times \mathbb{S}^1$ and $\mathbb{S}^{n-1} \simtimes \mathbb{S}^1$.
\end{proof}

\subsection{Shrinking necks}
For a smooth mean curvature flow $F:\mathcal{M}\times[0,T] \to \mathbb{R}^{n+m}$ let $\mathcal{B}_{g(t)}(p,r)\subset\mathcal{M}$ be the intrinsic closed ball of radius $r>0$ about $p\in\mathcal{M}$ with respect to the metric $g(t)$. If $t$ and $\theta$ are given such that $0\leq t-\theta < t \leq T$, we define
    \[\mathcal{P}(p,t,r,\theta) = \left\{(q,s) : q \in \mathcal B_{g(t)}(p,r), s \in [t-\theta,t] \right\}.\]
Now consider instead a mean curvature flow with surgery. If $t$ is a surgery time then we write $g(t-)$ and $g(t+)$ for the metric before and after surgery and use the convention $g(t) = g(t-)$. As in \cite{HuSi09}, if there exist no points of $\mathcal B_{g(t)}(p,r)$ which belong to a region changed by a surgery at some time in the interval $(t-\theta,t]$, then we define $\mathcal{P}(p,t,r,\theta)$ as in the smooth case and say $\mathcal{P}(p,t,r,\theta)$ does \emph{not contain surgeries}. We also use the notation
    \[\hat{\mathcal P}(p,t,L,\theta) := \mathcal P(p,t,\hat r(p,t)L,\hat r(p,t)^2\theta), \qquad \hat{r}(p,t) := \frac{n-1}{|H(p,t)|}.\]
Observe that if $(p,t)$ lies on a neck then $\hat r(p,t)$ is approximately the radius of the neck.

Consider for $s \leq 0$ the function $\rho(r,s) = \sqrt{r^2 - 2(n-1)s}$, so that $\rho(r,s)$ is the radius at time $s$ of a standard $n$-dimensional cylinder which evolves by mean curvature flow and has radius $r$ at time $s=0$.

\begin{definition} A point $(p_0,t_0)$ lies at the centre of an \emph{$(\e,k,L,\theta)$-shrinking neck} if after setting $r_0 = \hat{r}_0(p_0,t_0)$ and $\mathcal B_0 = \mathcal B_{g(t_0)}(p_0,r_0 L)$ the following hold:
\begin{itemize}
\item the parabolic neighbourhood $\hat{\mathcal P}(p_0,t_0,L,\theta)$ does not contain surgeries;
\item for every $t\in [t_0-r_0^2\theta,t_0]$, the region $\mathcal B_0$ with respect to the immersion $F(\cdot,t)$ is an $(\e, k)$-neck at scale $\rho(r_0,t- t_0)$ (where if $t_0 - r_0^2\theta$ is a surgery time we consider the manifold after surgery).
\end{itemize}
When these properties hold we also say that $\hat{\mathcal P}(p_0,t_0,L,\theta)$ is an $(\e,k)$-shrinking neck.
\end{definition}

\subsection{Almost-planar submanifolds}

Let $\mathcal M$ be an $n$-dimensional submanifold of $\mathbb{R}^{n+m}$. We assume $|H| > 0$ so that the second fundamental form of $\mathcal M$ decomposes as $A = h \nu + A^-$. We assume in addition that 
    \begin{equation}\label{eq_scalar curvature pinching}
    h \leq (1-\delta)|H|g
    \end{equation}
for a positive constant $\delta$. Combined with \eqref{eq_principal torsion} this implies the useful inequality
    \begin{equation}\label{eq_principal torsion estimate}
    |H| |\nabla \nu| \leq C(n,\delta)|\nabla A^-|.
    \end{equation}

In \cite{Naff_planarity}, Naff observed that if $\mathcal M$ satisfies \eqref{eq_scalar curvature pinching} and $A^-$ vanishes identically then $\mathcal M$ is a hypersurface. That is, there is an $(n+1)$-dimensional affine subspace of $\mathbb{R}^{n+m}$ which contains $\mathcal M$. This explains why, in situations where the planarity estimate holds, blow-ups at a singularity of the mean curvature flow are hypersurfaces. In this final preliminary section we establish estimates asserting that if $A^-$ is small then the position vector of the submanifold lies close to an $(n+1)$-dimensional affine subspace. The subspace is chosen to be parallel to $\operatorname{span}\{\nu(p), T_p\mathcal M\}$ for some $p \in \mathcal M$. First, we compare the normal spaces at neighbouring points with the normal space at $p$. 

\begin{lemma}\label{lem_almost parallel frame}
Suppose $\mathcal M$ satisfies \eqref{eq_scalar curvature pinching} and fix a point $p \in \mathcal M$. Let $\{\omega_a\}_{a = n+1}^{n+m}$ be constant vectors in $\mathbb{R}^{n+m}$ which form an orthonormal basis for $N_p \mathcal M$. In addition, suppose $\omega_{n+1} = \nu(p)$. For every $q$ in the intrinsic ball $\mathcal B(p, \Lambda r) \subset \mathcal M$ and $a \in \{n+2, \dots, n+m\}$ we have 
    \[|\omega_a \cdot \nu(q)| + |\omega_a^\top(q)| \leq \exp(C\Lambda) \cdot  \sup_{\mathcal B(p,\Lambda r)} \Big(r|A^-| + r|H|^{-1}|\nabla A^-|\Big)\]
where $C$ depends only on $n$, $\delta$ and $\sup_{\mathcal B(p,\Lambda r)} r |A|$.
\end{lemma}
\begin{proof}
Let $\gamma = \gamma(s)$ be a minimizing unit-speed geodesic from $p$ to $q$. Let $\{e_i\}_{i=1}^n$ be a parallel orthonormal frame for $T_{\gamma(s)} M$. First observe that
    \[D_{\gamma'(s)} \nu = (D_{\gamma'(s)} \nu)^\top + (D_{\gamma'(s)} \nu)^\perp = -h(\gamma'(s), e_i)e_i + \nabla_{\gamma'(s)} \nu,\]
and
    \[D_{\gamma'(s)} e_i = (D_{\gamma'(s)} e_i)^\perp = h(\gamma'(s), e_i) \nu + A^-(\gamma'(s), e_i).\]
Using \eqref{eq_principal torsion estimate} to estimate $\nabla \nu$, we see that  
    \[f(s) := |\omega_a \cdot \nu(\gamma(s))| + \sum_{i=1}^n |\omega_a \cdot e_i(\gamma(s))|\]
satisfies 
    \[\frac{d}{ds} f \leq Cr^{-1}f + C \sup_{\mathcal B(p,\Lambda r)} (|A^-| + |H|^{-1}|\nabla A^-|).\]
The claim follows by ODE comparison, using the fact that $f(0) = 0$ for $a \in \{n+2, \dots, n+m\}$.  
\end{proof}

Next we consider the coordinate functions $F^a(x) = x \cdot \omega_a$. We expect the derivatives of $F^a$ to be small for indices $a \in \{n+2, \dots, n+m\}$ provided that $A^-$ and its derivatives are small. So we differentiate $F^a$ and isolate terms involving $\omega_a^\top$ and $\nu \cdot \omega_a$, since these can be estimated in terms of $A^-$ using Lemma~\ref{lem_almost parallel frame} for $a \in \{n+2, \dots, n+m\}$.

\begin{lemma}\label{lem_derivative estimates F alpha}
We consider the same setup as in Lemma~\ref{lem_almost parallel frame} and set $F^a(x) = x \cdot \omega_a$ for each $a \in \{n+1, \dots, n+m\}$. We then have
    \begin{enumerate}[label=(\roman*)]
    \item\label{first est F alpha} $|\nabla F^a| = |\omega_a^\top|$,\\
    \item\label{second est F alpha} $|\nabla^2 F^a| \leq  |h||\nu \cdot \omega_a| + |A^-|$,\\
    \item\label{third est F alpha} $|\nabla^3 F^a| \leq |A|^2|\omega_a^\top| + |\nu \cdot \omega_a||\nabla h| + C|\nabla A^-|$, and\\
    \item for each $\ell \geq 1$ 
    \begin{align}\label{higher ests F alpha}
    |\nabla^{\ell+3} F^a| & \leq |\nabla^{\ell+1} A^-| +|\nabla^{\ell+1} h||\nu \cdot \omega_a| +C\sum_{\ell_1 + \ell_2 + \ell_3 = \ell} |\nabla^{\ell_1} A||\nabla^{\ell_2} A||\nabla^{\ell_3+1}F^a|\notag\\
    &+ C\sum_{\ell_1 + \ell_2 = \ell} \,\sum_{k_1 + k_2 + k_3 = \ell_2} |\nabla^{\ell_1} h||\nabla^{k_1} P| |\nabla^{k_2} |H|^{-1}||\nabla^{k_3+1} A^-|
    \end{align}
    \end{enumerate}
where $P$ is the $(1,1)$-tensor such that $P^k_j(|H|\delta_k^i - h_k^i)=\delta^i_j$. The constants appearing in these estimates depend only on $n$ and $\delta$. 
\end{lemma}
\begin{proof}
We compute at a point $q \in \mathcal M$, using a local orthonormal frame $\{e_i\}$ for $T\mathcal{M}$ such that $\nabla_{e_i} e_j (q) = 0$. One easily checks that $\nabla F^a = \omega_a^\top$ and $\nabla^2 F^a = \langle A, \omega_a\rangle$, and hence \ref{first est F alpha} and \ref{second est F alpha} hold. To estimate the third derivatives of $F^a$ we begin by writing
    \[\nabla_i \nabla_j \nabla_k F^a = e_i \langle A_{jk},\omega_a\rangle = \langle e_i(A_{jk}),\omega_a\rangle,\]
and then split the derivative on the right-hand side into tangential and normal components,
    \begin{align*}
    e_i(A_{jk}) = -\langle A_{il}, A_{jk}\rangle e_l + \nabla_i A_{jk},
    \end{align*}
in order to obtain
    \begin{equation}\label{third derivative F alpha}
    \nabla^3 F^a = -\langle A(\cdot, \omega_a^\top), A\rangle + \langle \nabla A, \omega_a\rangle.
    \end{equation}
By expanding $\nabla A = \nabla h \nu + h \nabla \nu + \nabla A^-$ we see that 
    \[|\nabla^3 F^a| \leq |A|^2|\omega_a^\top| + |\nu \cdot  \omega_a||\nabla h| + |\nabla \nu||h| + |\nabla A^-|.\]
From this \ref{third est F alpha} is obtained using \eqref{eq_principal torsion estimate}.

Continuing to differentiate \eqref{third derivative F alpha}, one arrives at an expression of the form
    \[\nabla^{\ell+3} F^a = \sum_{\ell_1 + \ell_2 + \ell_3 = \ell} \nabla^{\ell_1} A \ast \nabla^{\ell_2} A \ast \nabla^{\ell_3}\omega_a^\top + \langle \nabla^{\ell+1} A, \omega_a\rangle.\]
Since $\omega_a^\top = \nabla F^a$ and 
    \[\nabla^{\ell+1} A = \nabla^{\ell+1} h \nu + \sum_{\ell_1 + \ell_2 = \ell} \nabla^{\ell_1} h \ast \nabla^{\ell_2+1} \nu + \nabla^{\ell+1} A^-,\]
we may equally write
    \begin{align*}
    \nabla^{\ell+3} F^a &= \sum_{\ell_1 + \ell_2 + \ell_3 = \ell} \nabla^{\ell_1} A \ast \nabla^{\ell_2} A \ast \nabla^{\ell_3+1}F^a + \langle \nu, \omega_a\rangle \nabla^{\ell+1} h\notag\\
    &+ \sum_{\ell_1 + \ell_2 = \ell} \nabla^{\ell_1} h \ast \nabla^{\ell_2+1} \nu \ast \omega_a + \langle \nabla^{\ell+1} A^-, \omega_a\rangle.
    \end{align*}
Together with \eqref{eq_principal torsion} this implies \eqref{higher ests F alpha}.
\end{proof}

Lemma~\ref{lem_almost parallel frame} and Lemma~\ref{lem_derivative estimates F alpha} together imply the following estimates, which will be made use of when we construct our standard surgery procedure in Section~\ref{sec_standard surgery}.

\begin{proposition}\label{prop_derivative estimates F alpha simple}
Suppose we are in the setting of Lemma~\ref{lem_almost parallel frame}. Let us assume in addition that $|H| \geq \kappa r^{-1}$ holds in $\mathcal B(p,\Lambda r)$ for a constant $\kappa > 0$. We then have, for each index $a \in \{n+2, \dots, n+m\}$,
    \[|\nabla F^a| + r|\nabla^2 F^a| + r^2|\nabla^3 F^a|\leq C_1 \sup_{\mathcal B(p,\Lambda r)} (r|A^-| + r^2|\nabla A^-|) \]
where $C_1$ depends only on $n, \Lambda, \delta$ and $\sup_{\mathcal B(p,\Lambda r)} r|A|$, and for $\ell \geq 4$,
    \begin{align*}
    r^{\ell-1}|\nabla^{\ell} F^a| &\leq C_2\sum_{k=2}^{\ell-2} r^{k+1}|\nabla^k A^-| + C_2 \sup_{\mathcal B(p,\Lambda r)} (r|A^-| + r^2|\nabla A^-|)
    \end{align*}
where $C_2$ depends only on $n, \Lambda, \delta, \kappa$ and $\sum_{k=0}^{\ell-2}\sup_{\mathcal B(p,\Lambda r)} r^{k+1}|\nabla^k A|$. These estimates hold everywhere in $\mathcal B(p, \Lambda r)$.  
\end{proposition}
\begin{proof}
The estimates for $\nabla^\ell F^a$, $1 \leq \ell \leq 3$, follow immediately from Lemma~\ref{lem_almost parallel frame} and Lemma~\ref{lem_derivative estimates F alpha}. Given that we are assuming $|H| \geq \kappa r^{-1}$, for $\ell \geq 4$ the inequality
    \begin{align*}
    r^{\ell-1}|\nabla^{\ell} F^a| &\leq C\sum_{k=0}^{\ell-2} r^{k+1}|\nabla^k A^-| +C|\nu \cdot \omega_a|+C|\omega_a^\top|
    \end{align*}
follows from Lemma~\ref{lem_derivative estimates F alpha} and a simple proof by induction, with $C$ depending only on $n, \delta, \kappa$ and $\sum_{k=0}^{\ell-2}\sup_{\mathcal B(p,\Lambda r)} r^{k+1}|\nabla^k A|$. The claim now follows when combine this with Lemma~\ref{lem_almost parallel frame}. 
\end{proof}

\section{The standard surgery procedure}\label{sec_standard surgery}

In this section we construct our standard surgery procedure, which removes the middle portion of a sufficiently long neck and replaces it with two smoothly attached caps. This procedure must be designed extremely carefully to maintain control on the first and second fundamental forms---otherwise, our flow with surgery will not terminate after finitely many steps. The most involved part of the analysis is an estimate for the component $A^-$. We cannot show that every portion of the neck becomes closer to a hypersurface under surgery, so we devote our analysis to estimating the amount by which $A^-$ increases; a precise estimate is needed here in order for our planarity estimate to be preserved under surgery (see Section~\ref{sec_planarity} for further discussion). This step is involved because $A^- = A - \langle A, H\rangle H/|H|^2$ is rather nonlinear in the derivatives of the immersion. 

Consider a smooth closed $n$-dimensional submanifold $\mathcal{M} \subset \R^{n+m}$, and let 
    \[N: \mathbb{S}^{n-1}\times [a,b] \to \mathcal{M}\]
be a normal $(\e,k)$-cylindrical submanifold neck, where $\varepsilon > 0$ is sufficiently small and $k \geq 4$. We choose $z_0 \in [a,b]$ such that $z\in [z_0 - 7\Lambda, z_0 + 7\Lambda]\subset [a,b]$ for some $\Lambda \geq 10$. Let us set $r_0 = r(z_0)$. The following definition specifies what it means to perform standard surgery on the submanifold neck $N$. Compared with the corresponding procedure for hypersurfaces used in \cite{HuSi09}, we introduce an extra step (see Part~\ref{surgery blend} immediately below) in which we interpolate between $N$ and a suitable hypersurface. This is natural in higher codimensions and simplifies our analysis substantially.

\begin{definition}[Standard surgery] \label{defn_standard}
Given a normal $(\e,k)$-cylindrical submanifold neck $N$ and $z_0$ as above, as well as parameters $0<\tau<1$ and $B \geq 12 \Lambda$, we define the \emph{standard surgery} procedure with parameters $(\tau, B)$ at the cross-section $ \Sigma_{z_0} = N(\mathbb{S}^{n-1} \times \{z_0\})$ in the following manner:
\begin{enumerate}[label=(\alph*)]
\item The two collars $  \mathbb{S}^{n-1} \times [a, z_0 - 4\Lambda]$ and $ \mathbb{S}^{n-1} \times [z_0+4\Lambda, b]$ are unchanged by surgery.

\item The two cylinders $ N(\mathbb{S}^{n-1} \times [z_0 - 4 \Lambda, z_0])$ and $ N(\mathbb{S}^{n-1} \times [z_0, z_0+ 4 \Lambda])$ are each replaced by an $n$-ball attached smoothly to $ \Sigma_{z_0 - 4\Lambda}$ and $ \Sigma_{z_0 + 4\Lambda}$ respectively. We only describe the procedure for the left portion $[z_0- 5 \Lambda, z_0]$, the right portion $[z_0, z_0+5\Lambda]$ being analogous. For convenience, let $ z_0 = 5 \Lambda$ and consider a normal parameterization $N:\mathbb{S}^{n-1}\times [0, 5 \Lambda]\rightarrow \mathcal{M}$. In addition, suppose without loss of generality that $N(\hat \omega, z_0) = 0$, where $\hat \omega$ is some point in $ \mathbb{S}^{n-1}$ which we fix arbitrarily, but consistently---choosing $\hat \omega$ to be the north pole will suffice. 

\item\label{surgery blend} To blend our neck into a hypersurface we first let $\{\omega_a\}_{a = 1}^{n+m}$ be an orthonormal frame for $\mathbb R^{n+m}$ such that $\omega_{n+1}, \dots, \omega_{n+m}$ are normal to $N$ at $(\hat \omega, z_0)$ and $\omega_{n+1} = \nu(\hat \omega, z_0)$, where $\nu$ is the principal normal to $N$. Next, we choose a smooth transition function $\varphi: [0,5\Lambda] \rightarrow \R$ with $ \varphi= 1$ on $[0, 2\Lambda]$, $\varphi = 0$ on $[3\Lambda, 5\Lambda]$ and $\varphi' \leq 0$. We then define 
    \[\bar N(\omega,z) := (N^1(\omega,z), \dots, N^{n+1}(\omega,z), \varphi(z) N^{n+2}(\omega,z), \dots, \varphi(z)N^{n+m}(\omega,z)),\]
where $N^a := \langle N, \omega_a\rangle$. This ensures that $\bar N$ is a hypersurface in the subspace $\omega_1 \wedge \dots \wedge \omega_{n+1}$ for $z \in [3\Lambda, 5\Lambda]$. Our choice of $\varphi$ only depends on $\Lambda$ and is defined in such a way that all of its derivatives are smaller if $\Lambda$ is larger. If $\Lambda \geq 10$, each derivative of $\varphi$ is bounded by some fixed constant.

\item Next we bend our submanifold in the direction of its principal normal over the region where $z \in [\Lambda, 5\Lambda]$. Let $u(z) := r_0 \exp \left( - \frac{B}{z - \Lambda} \right) $ for $z \in [\Lambda, 5\Lambda]$ where $B \geq 12\Lambda$. For each $\tau \in [0,1]$ we define
    \begin{align}\label{eqn_N}
    \tilde N(\omega,z,\tau) : = \bar N(\omega,z) + \tau u (z) \bar \nu(\omega, z),
    \end{align}
where $\bar \nu$ is the principal normal to $\bar N$. 

\item\label{surgery blend 2} Next we interpolate between $\tilde N$ and an axially symmetric hypersurface. We fix an approximating cylinder $C_{z_0} :  \mathbb{S}^{n-1} \times \R \rightarrow \R^{n+m}$ such that: the image of $C_{z_0}$ lies in the subspace $\omega_1 \wedge \dots \wedge \omega_{n+1}$, the radius of $C_{z_0}$ is $r_0$, the axis of $C_{z_0}$ contains the centre of mass of $\Sigma_{z_0}$ and is parallel to the normal space of $\Sigma_{z_0}$ inside $\mathcal M$ at $(\hat \omega,z_0)$. In addition, we fix a smooth transition function $\psi: [0,5\Lambda] \rightarrow \R$ such that $ \psi= 1$ on $[0, 3\Lambda]$, $\psi =0$ on $[4\Lambda, 5\Lambda]$ and $\psi'\leq 0$. Denote by $\tilde C_{z_0}(\cdot, \tau):  \mathbb{S}^{n-1} \times [0, 5\Lambda] \rightarrow \R^{n+m}$ the bending of the approximating cylinder defined above along its principal normal, i.e.\ $\tilde C_{z_0} := C_{z_0}(\omega,z) + \tau u(z) \nu_{C_{z_0}}(\omega, z)$. We then define
    \begin{align*}
    \hat N(\omega, z, \tau) := \psi(z) \tilde N(\omega, z, \tau) + (1- \psi(z))\tilde C_{z_0} (\omega,z).
    \end{align*}
Our choice of $\psi$ only depends on $\Lambda$ and is defined in such a way that all of its derivatives are smaller if $\Lambda$ is larger. If $\Lambda \geq 10$, each derivative of $\psi$ is bounded by some fixed constant.

\item In this last step we suitably change $u$ on $[4\Lambda, 5\Lambda]$ to a function $\hat u$ to ensure $\tau\hat u(z) \rightarrow r_0$ as $z$ approaches some $z_1 \in (4\Lambda, 5\Lambda]$, such that $\tilde C_{z_0} ([4\Lambda, 5\Lambda])$ is a smoothly attached, axially symmetric, strictly convex cap. Since this last deformation on $[4\Lambda, 5\Lambda]$ only concerns the axisymmetric case, it can be performed for each pair of parameters $(\tau, B)$ in such a way that on the resulting strictly convex cap there exists some fixed upper bound for the curvature and each of its derivatives, independent of $\Lambda \geq 10$ and the surgery parameters $(\tau, B)$. Moreover, we can assume that under this final deformation the principal curvatures increase pointwise at each $z$, in such a way that the quantities $|H|^2$ and $\frac{1}{n-1}|H|^2 - |A|^2$ also increase. 
\end{enumerate}
\end{definition}

We prove that if all of the parameters are chosen to lie in suitable ranges, then throughout the standard surgery procedure just described we maintain fine control on the geometry of the submanifold.  To achieve this we choose our parameters in the following order: For any $\Lambda \geq 10$ and $k \geq 4$ we first choose $B$ to be large depending on $\Lambda$, then set $\tau$ equal to some $\tau_0 > 0$ which depends on $n$ and $\Lambda$, and then finally restrict $\varepsilon \ll \tau_0$ depending on $n$ and $\Lambda$. We assume throughout that the mean radius of the neck is almost constant: $|r(z)/r_0 - 1| \leq \frac{1}{100}$ for every $z \in [-7\Lambda, 7\Lambda]$. This is always true if $\varepsilon$ is small enough. 

\subsection{Blending with a hypersurface} We examine $\bar N$ defined as in Definition~\ref{defn_standard}. That is we study the step in the standard surgery procedure which blends $N$ with a hypersurface over the interval $[2\Lambda, 3\Lambda]$. We will bound the change in the first and second fundamental forms under this blending purely in terms of $A^-$ and its derivatives. For each $\ell \in \{0, \dots, k\}$ let us define
    \[\Theta_\ell^- := \sup_{\mathcal B_g(p_0,6\Lambda r_0)} \sum_{i=0}^\ell r_0^{i+1} |\nabla^i A^-|,\]
where $p_0 := (\hat \omega, z_0)$ and $\mathcal B_g(p_0,6\Lambda r_0)$ is the intrinsic ball of radius $6 \Lambda r_0$ with respect to the induced metric $g$. We may assume $\Lambda$ is large enough and $\varepsilon$ is small enough so that 
    \[\mathbb{S}^{n-1} \times [-5\Lambda, 5\Lambda] \subset \mathcal B_g (p_0, 6\Lambda r_0) \subset  \mathbb{S}^{n-1} \times [-7\Lambda, 7\Lambda].\] 

\begin{remark}
On an $(\varepsilon,k)$-neck we have $\Theta_\ell^- = O(\varepsilon)$. We will construct our surgery algorithm such that whenever surgery is performed we have additional estimates for $\Theta_\ell^-$, coming from the planarity estimate and planarity improvement theorem (see Section~\ref{sec_planarity}). 
\end{remark}

\begin{lemma}\label{lem_coord derivs blended neck}
For each $\ell \in \{0, \dots, k\}$ the $C^{\ell+2}$-norm of $\bar N - N$ in $\mathbb{S}^{n-1} \times [0, 5\Lambda]$ is bounded by $C \Theta_\ell^- r_0$ where $C=C(n,\ell, \Lambda)$.
\end{lemma}
\begin{proof}
Given that $N$ is an $(\e, k)$-cylindrical submanifold neck we have $r_0^{i+1}|\nabla^i A| \leq C$ for each $0 \leq i \leq k$ and $r_0|H| \geq 1/C$. Applying Proposition~\ref{prop_derivative estimates F alpha simple}, we deduce that for each index $a \in \{n+2, \dots, n+m\}$ we have the pointwise estimate
    \begin{align*}
    \sum_{i = 1}^{\ell+2} r_0^{i-1}|\nabla^i N^a| &\leq C\sum_{i=2}^{\ell} r_0^{i+1}|\nabla^i A^-| + C \sup_{\mathcal B_g(p,6\Lambda r_0)} (r_0|A^-| + r_0^2|\nabla A^-|)
    \end{align*}
in the region where $z \in [0, 5\Lambda]$, and hence for $\ell \in \{0, \dots, k\}$
    \[\sup_{\mathbb{S}^{n-1}\times[0,5\Lambda]} \bigg(\sum_{i = 1}^{\ell+2} r_0^{i-1}|\nabla^{i} N^a|\bigg) \leq C \Theta_\ell^-.\]
The norm on the left is with respect to the induced metric $g$, which is comparable to the standard product metric on $ \mathbb{S}^{n-1} \times [0, 5\Lambda]$ scaled by $r_0^2$. Therefore, the estimate above implies that the $C^{\ell+2}$-norm of $N^a$ in $\mathbb{S}^{n-1} \times [0,5\Lambda]$ is bounded by $C\Theta_\ell^- r_0$ for each index $a \in \{n+2, \dots, n+m\}$. Here we have also used the normalisation $N(p_0) = 0$ to bound $N^a$ in terms of its gradient. The claim now follows, because
    \[\bar N(\omega,z) - N(\omega, z) = (0, \dots, 0, (\varphi(z)-1) N^{n+2}(\omega,z), \dots, (\varphi(z)-1)N^{n+m}(\omega,z))\]
and the derivatives of $\varphi$ are bounded by constants which depend only on $\Lambda$. 
\end{proof}

It follows that geometric quantities which can be computed from the first two derivatives of $\bar N$ are close to their counterparts for $N$ up to errors which are controlled by $\Theta_1^-$. We state the result in normal coordinates for the standard metric on $\mathbb{S}^{n-1} \times \mathbb R$ and use the notation $\alpha = \beta + O(\gamma)$ to mean that $|\alpha - \beta| \leq C\gamma$ for some $C = C(n,\Lambda)$. 

\begin{proposition}\label{geom quants blended neck}
If $\varepsilon$ is sufficiently small then at every point in $ \mathbb{S}^{n-1} \times [0, 5\Lambda]$ we have:
\begin{align*}
&\bar g_{ij} = g_{ij} + O(\Theta_1^- r_0^2), \quad \bar g^{ij} = g^{ij} + O(\Theta_1^- r_0^{-2}), \quad \bar A_{ij} = A_{ij} + O(\Theta_1^- r_0),\\
&\bar H = H + O(\Theta_1^- r_0^{-1}), \quad \bar \nu = \nu + O(\Theta_1^-), \quad \bar A_{ij}^- = A_{ij}^- + O(\Theta_1^- r_0).
\end{align*}
\end{proposition}
\begin{proof}
The first four statements are easily deduced from Lemma~\ref{lem_coord derivs blended neck} using $\bar g_{ij} = \langle \partial_i \bar N, \partial_j \bar N\rangle$, $\bar A_{ij} = \partial_i \partial_j \bar N - \bar \Gamma_{ij}^k \partial_k \bar N$, and the expression for $\bar \Gamma_{ij}^k$ in terms of second derivatives of $\bar N$. The remaining two statements then follow by writing $\bar \nu = \bar H/ |\bar H|$ and $\bar A_{ij}^- = \bar A_{ij} - \langle \bar A_{ij}, \bar \nu\rangle \bar \nu$.
\end{proof}

We note the following consequence of Proposition~\ref{geom quants blended neck}.  

\begin{proposition}\label{geom quants blended neck 2}
If $\varepsilon$ is sufficiently small then at every point in $ \mathbb{S}^{n-1} \times [0, 5\Lambda]$ we have:
\[|\bar A|^2_{\bar g} = |A|^2_g + O(\Theta_1^- r_0^{-2}), \quad |\bar H|^2 = |H|^2 + O(\Theta_1^- r_0^{-2}), \quad |\bar A^-|^2_{\bar g} = O(|\Theta_1^-|^2r_0^{-2}),\]
and for $1 \leq \ell \leq k$ 
\[|\bar \nabla^\ell \bar A|^2_{\bar g} \leq |\nabla^\ell A|^2 + C(n,\ell,\Lambda)\Theta_\ell^- r_0^{-2(1+\ell)}, \quad |\bar \nabla^\ell \bar A^-|^2_{\bar g} \leq  C(n,\ell,\Lambda) |\Theta_\ell^-|^2r_0^{-2(1+\ell)}.\]
\end{proposition}
\begin{proof}
The first three statements are immediate from Proposition~\ref{geom quants blended neck}. To obtain the final two statements, we differentiate the coordinate expressions for $\bar A_{ij}$ and $\bar A_{ij}^-$ to obtain 
    \[\bar \nabla^\ell \bar A_{ij} = \nabla^\ell A_{ij} + O(\Theta_\ell^- r_0), \qquad \bar \nabla^\ell \bar A_{ij}^- = \nabla^\ell A_{ij}^- + O(\Theta_\ell^- r_0)\]
from Lemma~\ref{lem_coord derivs blended neck} and then use the estimate for $\bar g^{ij}$ in Proposition~\ref{geom quants blended neck}. 
\end{proof}
    
\subsection{Bending along the principal normal} We consider the bent neck $\tilde N(\omega,z,\tau)$ given by $\bar N(\omega,z) + \tau u (z) \bar \nu(\omega, z)$ where $u = u(z)$ is the function introduced in Definition~\ref{defn_standard}. We assume throughout that $\tau \in [0,1]$ is small enough to ensure
    \begin{equation}\label{u is small}
    \tau(|u| + |u'| + |u''|) \leq 10^{-3}r_0
    \end{equation}
for $z \in [0, 5\Lambda]$. We write $C$ for any constant which depends only on $n$ and $\Lambda$. Around each point in $\mathbb{S}^{n-1}\times [0,5\Lambda]$ we can choose normal coordinates $\{x^i\}$ for the standard metric on $\mathbb{S}^{n-1} \times \mathbb{R}$ such that $\partial_1$ coincides with $\partial_z$.

The following statement is obtained by differentiating the definition of $\tilde N$ and inserting the identities $\partial_i \bar \nu = - \bar h_i^k \partial_k \bar N + \bar \nabla_i \bar \nu$ and $\partial_i \partial_k \bar N = \bar h_{ik}\bar \nu + \bar A_{ik}^- + \bar \Gamma_{ik}^l\partial_l \bar N$.
    
\begin{lemma}\label{lem_derivatives of tilde N}
We have
    \begin{align*}
    \partial_i \tilde N &= \partial_i \bar N + \tau \delta^1_i u' \bar \nu - \tau u \bar h_i^k \partial_k \bar N + \tau u \bar \nabla_i \bar \nu,\\
    \partial_i \partial_j \tilde N &=\partial_i\partial_j \bar N + \tau (\delta^1_i\delta^1_j u'' -u\bar h^l_j\bar h_{il})\bar \nu\\
    & -\tau (\delta^1_ju' \bar h_i^k + \delta^1_i u' \bar h_j^k + u \partial_i \bar h_j^k + u \bar h_j^l\bar \Gamma_{il}^k) \partial_k \bar N\\
    &- \tau u \bar h_j^l \bar A_{il}^- + \tau \delta^1_j u' \bar \nabla_i \nu + \tau \delta^1_i u' \bar \nabla_j \bar \nu + \tau u \partial_i \bar \nabla_j \bar \nu.
    \end{align*}
\end{lemma}

We introduce for each $\tau \in [0,1]$ an orthonormal frame $\{\tilde \nu_\alpha\}_{\alpha = 1}^{m}$ for the normal space of $\tilde N(\cdot, \tau)$ such that $\tilde \nu_1 = \tilde \nu = \tilde H/|\tilde H|$. We require that $\tilde \nu_\alpha$ is smooth in $\tau$ and 
    \[\langle\partial_\tau\tilde \nu_\alpha, \tilde \nu_\beta\rangle = 0 \; \forall \; \alpha,\beta \geq 2.\]
This system of ODEs can be solved using standard theory and will simplify our computations. For each $\alpha \in \{1, \dots, m\}$ we have $0 = \langle \partial_\tau\tilde\nu_\alpha, \partial_i \tilde N\rangle + \langle \tilde\nu_\alpha, \partial_\tau \partial_i \tilde N\rangle$, which together with Lemma~\ref{lem_derivatives of tilde N} gives
    \begin{align}\label{deformation of normal tangent}
    (\partial_\tau \tilde\nu_\alpha)^\top &=-u'\langle\tilde\nu_\alpha, \bar \nu\rangle\tilde g^{1j}\partial_j\tilde N + u\langle\tilde\nu_\alpha,\partial_k \bar N\rangle \bar h^{kj}\partial_j\tilde N - u\langle\tilde\nu_\alpha,\bar \nabla_i \bar \nu\rangle  \tilde g^{ij}\partial_j\tilde N.
    \end{align}
In addition, for $\alpha \geq 2$ we have $\partial_\tau \tilde \nu_\alpha = -\langle\tilde \nu_\alpha, \partial_\tau \tilde \nu_1\rangle \tilde \nu_1 + (\partial_\tau \tilde \nu_\alpha)^\top$ and hence
    \begin{align}\label{deformation of normal}
    \partial_\tau \tilde \nu_\alpha &=-|\tilde H|^{-1}\tilde g^{kl} \langle\tilde \nu_\alpha, \partial_\tau \tilde A_{kl}\rangle \tilde \nu_1 - |\tilde H|^{-1}\partial_\tau\tilde g^{kl} \langle\tilde\nu_\alpha, \tilde A_{kl}\rangle\tilde\nu_1 + (\partial_\tau \tilde \nu_\alpha)^\top.
    \end{align}

Let us define 
    \[Y^- := r_0^{-1}\sum_{i,j} |\bar A_{ij}^-| + \sum_{i} |\bar \nabla_i \bar \nu| + \sum_{i,j} |\partial_i \bar \nabla_j \bar \nu|.\]
Throughout our analysis we treat terms bounded by $Y^-$ as small errors, since in the worst case $Y^- = O(\varepsilon)$ due to Proposition~\ref{geom quants blended neck 2} and the following  consequence of \eqref{eq_principal torsion}. 

\begin{lemma}\label{estimate for Y}
If $\varepsilon$ is sufficiently small then
    \[|Y^-|^2 \leq C r_0^2|\bar A^-|^2_{\bar g} + Cr_0^4 |\bar \nabla \bar A^-|^2_{\bar g} + Cr_0^6|\bar \nabla^2 \bar A^-|^2_{\bar g}\]
where the constant $C$ depends only on $n$. 
\end{lemma}

\noindent \textbf{Estimates for $\tilde A^-$}. We now come to the most involved part of our analysis in this section. The claim is that, when we bend $\bar N$ to obtain $\tilde N$, the component $A^-$ can become larger but only be a very small amount. To be precise, we prove that
    \[|\tilde A^-|^2_{\tilde g} = |\bar A^-|^2_{\bar g} + O(\tau r^{-3}(|u| + |u'| + |u''|))|Y^-|^2.\]
The crucial point is that the error term on the right contains the factor $|Y^-|^2$ and hence can be bounded in terms of $A^-$ and its first two derivatives. The estimate holds trivially in codimension one and cannot be approached using the arguments developed in Section~3 of \cite{HuSi09}.

The structure of the proof is as follows. After some preliminary analysis, in Lemma~\ref{deformation A-} we estimate $\partial_\tau |\langle\tilde \nu_\alpha, \tilde A^-\rangle|_{\bar g}^2$ in terms of $\langle \tilde \nu_\alpha, \bar \nu_1\rangle$ and $Y^-$. We can of course bound $|\langle \tilde \nu_\alpha, \bar \nu_1\rangle| \leq 1$ and integrate, but this turns out to be too crude, so we proceed more carefully. We derive an estimate for $\partial_\tau \langle \tilde \nu_\alpha, \bar \nu_1\rangle^2$ and use ODE comparison to bound a suitable sum of the quantities $|\langle\tilde \nu_\alpha, \tilde A^-\rangle|_{\bar g}^2$ and $\langle \tilde \nu_\alpha, \bar \nu_1\rangle^2$ in Lemma~\ref{Gronwall}. This leads to the desired estimate for $|\tilde A^-|_{\tilde g}^2$ in Theorem~\ref{A- norm in deformed metric}.

We first estimate the component of $\tilde \nu_\alpha$ which is tangent to the initial neck. 

\begin{lemma}\label{bent normal with tangent}
If \eqref{u is small} holds and $\varepsilon$ is small we have $|\langle \tilde\nu_\alpha,\partial_i \bar N\rangle| \leq C\tau |u'| |\langle\tilde\nu_\alpha,\bar \nu\rangle| + C\tau |u| |Y^-|$.
\end{lemma}
\begin{proof}
Using $\langle \tilde \nu_\alpha, \partial_i \bar N\rangle = \langle \tilde \nu_\alpha, \partial_i \bar N - \partial_i \tilde N\rangle$ and Lemma~\ref{lem_derivatives of tilde N} we obtain 
    \[(\delta_j^i - \tau uh_j^i)\langle \tilde\nu_\alpha, \partial_i \bar N\rangle = -\tau\delta_j^1u'\langle\tilde\nu_\alpha,\bar \nu\rangle - \tau u \langle\tilde\nu_\alpha,\bar \nabla_j \bar \nu\rangle\]
and hence
    \[\langle \tilde\nu_\alpha,\partial_i \bar N\rangle = -\tau u'P_i^1\langle\tilde\nu_\alpha,\bar \nu\rangle -\tau u P_i^j \langle\tilde\nu_\alpha,\bar \nabla_j \bar \nu\rangle,\]
where $P_i^j$ denotes the matrix inverse to $\delta_j^i - \tau u\bar h_j^i$. From Proposition~\ref{geom quants blended neck} we know that $|\bar h_j^i| \leq (1+C\varepsilon)r_0^{-1}$, so by \eqref{u is small} if $\varepsilon$ is small then $|P_j^i| \leq C$. The claim follows. 
\end{proof}

The following estimates follow easily from Lemma~\ref{lem_derivatives of tilde N} and $\tilde g_{ij} = \langle \partial_i \tilde N, \partial_j \tilde N\rangle$. 

\begin{lemma}\label{deformation of metric}
Assuming \eqref{u is small} holds, we have $|\partial_\tau \tilde g_{ij}| \leq Cr_0|u| + C\tau|u'|^2$ and $|\partial_\tau \tilde g^{ij}| \leq Cr_0^{-3}|u| + C\tau r_0^{-4}|u'|^2$.
\end{lemma}

Each of the next two lemmas provides an estimate for one of the terms on the right-hand side of the identity $\partial_\tau \langle \tilde \nu_\alpha, \tilde A_{ij}\rangle = \langle \partial_\tau\tilde \nu_\alpha, \tilde A_{ij}\rangle + \langle \tilde \nu_\alpha, \partial_\tau\tilde A_{ij}\rangle$. 

\begin{lemma}\label{deformation normal with A}
Assuming \eqref{u is small} holds, for each $\alpha \geq 2$ we have
    \begin{align*}
    \sum_{i,j} |\langle \partial_\tau \tilde \nu_\alpha, \tilde A_{ij}\rangle| &\leq C\sum_{i,j}|\langle\tilde \nu_\alpha, \partial_\tau \tilde A_{ij}\rangle| + C(r_0^{-1}|u| + \tau r_0^{-2}|u'|^2)\sum_{i,j}|\langle\tilde\nu_\alpha, \tilde A_{ij}\rangle|.
    \end{align*}
\end{lemma}
\begin{proof}
Taking the inner product of \eqref{deformation of normal} with $\tilde A_{ij}$ gives
    \begin{align*}
    \langle \partial_\tau \tilde \nu_\alpha, \tilde A_{ij}\rangle &=-|\tilde H|^{-1}\tilde g^{kl} \langle\tilde \nu_\alpha, \partial_\tau \tilde A_{kl}\rangle \langle\tilde \nu_1,\tilde A_{ij}\rangle - |\tilde H|^{-1}\partial_\tau\tilde g^{kl} \langle\tilde\nu_\alpha, \tilde A_{kl}\rangle\langle\tilde\nu_1,\tilde A_{ij}\rangle.
    \end{align*}
Lemma~\ref{lem_derivatives of tilde N} and Lemma~\ref{deformation of metric} easily imply the crude bounds
    \[|\tilde g^{ij}| \leq Cr_0^{-2}, \qquad |\tilde A_{ij}|\leq Cr_0, \qquad |\tilde H|^{-1} \leq Cr_0,\]
so we obtain
    \begin{align*}
    \sum_{i,j} |\langle \partial_\tau \tilde \nu_\alpha, \tilde A_{ij}\rangle| &\leq C\sum_{i,j}|\langle\tilde \nu_\alpha, \partial_\tau \tilde A_{ij}\rangle| + Cr_0^2 \sum_{i,j}|\partial_\tau g^{ij}| \cdot \sum_{i,j} |\langle\tilde\nu_\alpha, \tilde A_{ij}\rangle |.
    \end{align*}
The claim now follows from Lemma~\ref{deformation of metric}. 
\end{proof}

\begin{lemma}\label{normal with deformation A}
Assuming \eqref{u is small} holds, if $\varepsilon$ is sufficiently small then we have
    \[\sum_{i,j}|\langle\tilde\nu_\alpha, \partial_\tau\tilde A_{ij}\rangle| \leq C(|u| + |u'| + |u''|)|\langle\tilde\nu_\alpha,\bar \nu\rangle|+C(|u| + |u'|)Y^-.\]
\end{lemma}
\begin{proof}
We first use $\tilde A_{ij} = \partial_i\partial_j\tilde N - \tilde \Gamma_{ij}^k \partial_k \tilde N$ to compute
    \[\langle\tilde\nu_\alpha, \partial_\tau\tilde A_{ij}\rangle = \langle \tilde \nu_\alpha, \partial_\tau\partial_i\partial_j\tilde N - \partial_\tau\tilde \Gamma_{ij}^k \partial_k \tilde N- \tilde \Gamma_{ij}^k \partial_\tau\partial_k \tilde N\rangle = \langle \tilde \nu_\alpha, \partial_\tau\partial_i\partial_j\tilde N - \tilde \Gamma_{ij}^k \partial_\tau\partial_k \tilde N\rangle\]
and then use $|\tilde \Gamma_{ij}^k| \leq C$ to deduce
    \[\sum_{i,j}|\langle\tilde\nu_\alpha, \partial_\tau\tilde A_{ij}\rangle| \leq C\sum_{i,j} |\langle\tilde\nu_\alpha,\partial_\tau\partial_i \partial_j\tilde N\rangle| + C\sum_i |\langle\tilde\nu_\alpha, \partial_\tau\partial_i\tilde N\rangle|.\]
Now we estimate each of the two terms on the right. Lemma~\ref{lem_derivatives of tilde N} yields
    \begin{align*}
    \langle\tilde\nu_\alpha, \partial_\tau\partial_i\tilde N\rangle &= \delta^1_i u' \langle\tilde\nu_\alpha, \bar \nu\rangle - u \bar h_i^k \langle\tilde\nu_\alpha,\partial_k \bar N\rangle + u \langle \tilde\nu_\alpha, \bar \nabla_i \bar \nu\rangle,
    \end{align*}
which together with Lemma~\ref{bent normal with tangent} implies
    \begin{align*}
    \sum_i |\langle\tilde\nu_\alpha, \partial_\tau\partial_i\tilde N\rangle| &\leq C|u'| |\langle\tilde\nu_\alpha, \bar \nu\rangle| +C|u| Y^-.
    \end{align*}
Next we use Lemma~\ref{lem_derivatives of tilde N} to bound
    \begin{align*}
    \sum_{i,j} |\langle\tilde\nu_\alpha,\partial_\tau\partial_i \partial_j\tilde N\rangle| &\leq C(|u| + |u''|) |\langle\tilde\nu_\alpha,\bar \nu\rangle| + Cr^{-1}(|u| + |u'|) \sum_i|\langle\tilde\nu_\alpha,\partial_i \bar N\rangle|\\
    &+C(|u| + |u'|)Y^-.
    \end{align*}
To claim follows once we estimate the final term on the first line using Lemma~\ref{bent normal with tangent}.
\end{proof}

By combining Lemma~\ref{deformation normal with A} with Lemma~\ref{normal with deformation A} we obtain:

\begin{lemma}\label{deformation A-}
Assuming \eqref{u is small} holds, if $\varepsilon$ is sufficiently small then for each $\alpha \geq 2$ we have 
    \begin{align*}
    \partial_\tau |\langle\tilde\nu_\alpha, \tilde A\rangle|^2_{\bar g} &= O(r_0^{-2}(|u| + |u'| + |u''|))|\langle\tilde\nu_\alpha,\bar \nu\rangle||\langle\tilde\nu_\alpha, \tilde A\rangle|_{\bar g}\\
    &+O(r_0^{-2}(|u| + |u'|))  Y^- |\langle\tilde\nu_\alpha, \tilde A\rangle|_{\bar g}\\
    &+ O(r_0^{-1}(|u|+|u'|))|\langle\tilde\nu_\alpha, \tilde A\rangle|^2_{\bar g}.
    \end{align*}
\end{lemma}
\begin{proof}
We use the Cauchy--Schwarz inequality to estimate 
    \begin{align*}
    |\partial_\tau |\langle\tilde\nu_\alpha, \tilde A\rangle|^2_{\bar g}| &\leq 2 |\langle\tilde\nu_\alpha, \tilde A\rangle|_{\bar g}|\langle\tilde\nu_\alpha, \partial_\tau \tilde A\rangle|_{\bar g} + 2 |\langle\tilde\nu_\alpha, \tilde A\rangle|_{\bar g}|\langle\partial_\tau\tilde\nu_\alpha, \tilde A\rangle|_{\bar g}.
    \end{align*}
Using Lemma~\ref{deformation normal with A} to bound the second term on the right gives
    \begin{align*}
    |\partial_\tau |\langle\tilde\nu_\alpha, \tilde A\rangle|^2_{\bar g}| &\leq C |\langle\tilde\nu_\alpha, \tilde A\rangle|_{\bar g}|\langle\tilde\nu_\alpha, \partial_\tau \tilde A\rangle|_{\bar g} + Cr_0^{-1}(|u| + |u'|) |\langle\tilde\nu_\alpha, \tilde A\rangle|^2_{\bar g}.
    \end{align*}
The claim now follows from Lemma~\ref{normal with deformation A}.
\end{proof}

The estimate stated in Lemma~\ref{deformation A-} contains the term $\langle\tilde \nu_\alpha, \bar \nu\rangle$. We need to study how this quantity changes with $\tau$.  

\begin{lemma}\label{deformation normal with principal}
Assuming \eqref{u is small} holds, if $\varepsilon$ is sufficiently small then for each $\alpha \geq 2$ we have 
    \begin{align*}
    \partial_\tau \langle\tilde\nu_\alpha,\bar \nu\rangle^2 &= O(r_0^{-1}(|u|+|u'|+|u''|))\langle\tilde\nu_\alpha,\bar \nu\rangle^2\\
    &+O(|u| +  r_0^{-1}|u'|^2)|\langle\tilde\nu_\alpha,\bar \nu\rangle||\langle\tilde\nu_\alpha,\bar A\rangle|_{\bar g}\\
    &+O(r_0^{-1}(|u| + |u'|))Y^-|\langle\tilde\nu_\alpha,\bar \nu\rangle|.
    \end{align*}
\end{lemma}
\begin{proof}
The evolution \eqref{deformation of normal} immediately implies
    \begin{align*}
    |\langle \partial_\tau \tilde \nu_\alpha, \bar \nu\rangle| &\leq Cr_0^{-1}\sum_{k,l}|\langle\tilde \nu_\alpha, \partial_\tau \tilde A_{kl}\rangle| + Cr_0 \sum_{k,l} |\partial_\tau\tilde g^{kl}| \cdot  \sum_{k,l}|\langle\tilde\nu_\alpha, \tilde A_{kl}\rangle| + C |\langle (\partial_\tau \tilde \nu_\alpha)^\top, \bar \nu\rangle|.
    \end{align*}
Inserting Lemma~\ref{deformation of metric}, we see that
    \begin{align*}
    |\langle \partial_\tau \tilde \nu_\alpha, \bar \nu\rangle| &\leq Cr_0^{-1}\sum_{k,l}|\langle\tilde \nu_\alpha, \partial_\tau \tilde A_{kl}\rangle| + C(r_0^{-2}|u| + r_0^{-3}|u'|^2) \sum_{k,l} |\langle\tilde\nu_\alpha, \tilde A_{kl}\rangle| + C|\langle (\partial_\tau \tilde \nu_\alpha)^\top, \bar \nu\rangle|.
    \end{align*}
To estimate the final term on the right we first use \eqref{deformation of normal tangent} to bound 
    \[|\langle (\partial_\tau \tilde \nu_\alpha)^\top, \bar \nu\rangle| \leq Cr_0^{-1}|u'||\langle\tilde\nu_\alpha, \bar \nu\rangle| + Cr_0^{-2}|u| \sum_k |\langle\tilde\nu_\alpha,\partial_k \bar N\rangle| + Cr^{-1}|u|Y^-\]
and then combine this inequality with Lemma~\ref{bent normal with tangent} to obtain
    \[|\langle (\partial_\tau \tilde \nu_\alpha)^\top, \bar \nu\rangle| \leq Cr_0^{-1}|u'||\langle\tilde\nu_\alpha, \bar \nu\rangle| + Cr_0^{-1}|u|Y^-.\]
Substituting this back in above gives 
    \begin{align*}
    |\langle \partial_\tau \tilde \nu_\alpha, \bar \nu\rangle| &\leq Cr_0^{-1} \sum_{k,l}|\langle\tilde \nu_\alpha, \partial_\tau \tilde A_{kl}\rangle| + C(r_0^{-2}|u| + r_0^{-3}|u'|^2) \sum_{k,l}|\langle\tilde\nu_\alpha, \tilde A_{kl}\rangle|\\
    &+ Cr_0^{-1}|u'||\langle\tilde\nu_\alpha, \bar \nu\rangle| + Cr_0^{-1}|u|Y^-.
    \end{align*}
Finally, we estimate the first term on the right using Lemma~\ref{normal with deformation A} to get
    \begin{align*}
    |\langle \partial_\tau \tilde \nu_\alpha, \bar \nu\rangle| &\leq Cr_0^{-1} (|u| + |u'| + |u''|)|\langle\tilde\nu_\alpha,\bar \nu\rangle| + C(r_0^{-2}|u| + r_0^{-3}|u'|^2) \sum_{k,l}|\langle\tilde\nu_\alpha, \tilde A_{kl}\rangle|\\
    &+ Cr_0^{-1}(|u| + |u'|)Y^-.
    \end{align*}
The claim follows from this estimate and $\partial_\tau \langle\tilde\nu_\alpha,\bar \nu\rangle^2 = 2 \langle\tilde\nu_\alpha,\bar \nu\rangle\langle\partial_\tau\tilde\nu_\alpha,\bar \nu\rangle$.
\end{proof}

We now use Lemma~\ref{deformation A-} and Lemma~\ref{deformation normal with principal} to compare $|\tilde A^-|_{\bar g}^2$ with $|\bar A^-|_{\bar g}^2$.

\begin{lemma}\label{Gronwall}
Assuming \eqref{u is small} holds, if $\varepsilon$ is sufficiently small then for each $\alpha \geq 2$ we have 
    \[|\langle\tilde\nu_\alpha,\tilde A\rangle|_{\bar g}^2 = |\langle \bar \nu_\alpha, \bar A\rangle|_{\bar g}^2 + O(\tau r_0^{-3}(|u| + |u'| + |u''|))|Y^-|^2\]
and 
    \[\langle\tilde\nu_\alpha,\bar \nu\rangle^2 = O(\tau r_0^{-1} (|u| + |u'| + |u''|))|Y^-|^2.\] 
\end{lemma}
\begin{proof}
We set $\Phi := r_0^2|\langle\tilde\nu_\alpha,\tilde A\rangle|_{\bar g}^2$ and $\Psi:= \langle\tilde\nu_\alpha,\bar \nu\rangle^2$. Lemma~\ref{deformation A-} provides the estimate
    \begin{align*}
    |\partial_\tau \Phi| &\leq Cr_0^{-1}(|u| + |u'| + |u''|)\sqrt{\Psi}\sqrt{\Phi}\\
    &+Cr_0^{-1}(|u| + |u'|) Y^- \sqrt{\Phi}\\
    &+ Cr_0^{-1}(|u|+ |u'|)\Phi.
    \end{align*}
Lemma~\ref{deformation normal with principal} provides the estimate
    \begin{align*}
    |\partial_\tau \Psi| &\leq Cr_0^{-1}(|u|+|u'|+|u''|)\Psi\\
    &+Cr_0^{-1}(|u| + |u'|)\sqrt{\Psi}\sqrt{\Phi}\\
    &+Cr_0^{-1}(|u| + |u'|)Y^-\sqrt{\Psi}.
    \end{align*}
Using Young's inequality these two estimates yield
    \[|\partial_\tau \Phi| \leq Cr_0^{-1}(|u|+|u'|+|u''|)(\Phi+\Psi + |Y^-|^2)\]
and 
    \[|\partial_\tau \Psi| \leq Cr_0^{-1}(|u|+|u'|+|u''|)(\Phi+\Psi + |Y^-|^2),\]
which together with $\Phi(0) \leq C|Y^-|^2$ and $r_0^{-1}(|u|+|u'|+|u''|) \leq C$ imply
    \[\partial_\tau (|\Phi - \Phi(0)| + \Psi) \leq C(|\Phi - \Phi(0)| + \Psi) + Cr_0^{-1}(|u|+|u'|+|u''|)Y^2.\]
Given that $|\Phi - \Phi(0)|$ might only be Lipschitz, this inequality may only hold for a.e.\ $\tau$, but this is sufficient. It follows that the quantity
    \[e^{-C\tau}(|\Phi(\tau) - \Phi(0)|+\Psi(\tau)) - C\tau r_0^{-1}(|u|+|u'|+|u''|)|Y^-|^2 \]
is nonincreasing for $\tau \leq 1$ and hence 
    \begin{align*}|\Phi(\tau) - \Phi(0)| +\Psi(\tau) &\leq  C\tau  r_0^{-1}(|u|+|u'|+|u''|) |Y^-|^2,
    \end{align*}
where we have used $\Psi(0) = 0$. This completes the proof. 
\end{proof}

In Lemma~\ref{Gronwall} the norm of $\langle \tilde \nu_\alpha, \tilde A^-\rangle$ is with respect to $\bar g$. One readily obtains the following estimate for the norm with respect to the deformed metric $\tilde g$. 

\begin{theorem}\label{A- norm in deformed metric}
Assuming \eqref{u is small} holds, if $\varepsilon$ is sufficiently small then 
    \[|\tilde A^-|^2_{\tilde g} = |\bar A^-|^2_{\bar g} + O(\tau r_0^{-3}(|u| + |u'| + |u''|))|Y^-|^2.\]
\end{theorem}
\begin{proof}
From Lemma~\ref{deformation of metric} it follows that $\tilde g^{ij} = \bar g^{ij} + O(\tau r_0^{-3} (|u| + |u'|))$. Combined with the inequality $|\langle\tilde \nu_\alpha, \tilde A\rangle|_{\bar g}^2 \leq Cr_0^{-2} |Y^-|^2$, which is a consequence of Lemma~\ref{Gronwall}, this implies 
    \begin{align*}
    |\tilde A^-|^2_{\tilde g} &= \sum_{\alpha \geq 2} \tilde g^{ik}\tilde g^{jl}\langle \tilde \nu_\alpha, \tilde A_{ij}\rangle\langle \tilde \nu_\alpha, \tilde A_{kl}\rangle = |\tilde A^-|^2_{\bar g} + O(\tau r_0^{-3}(|u| + |u'|))|Y^-|^2,
    \end{align*}
To conclude we estimate the first term on the right using Lemma~\ref{Gronwall}. 
\end{proof}

\noindent \textbf{Estimates for $\tilde h$.} With our estimates for $\tilde A^-$ in place we now proceed to study $\tilde h$. Unlike $\tilde A^-$ the component $\tilde h$ was already dealt with in the codimension-one case, so we follow the arguments in \cite{HuSi09}.

\begin{lemma}\label{deformed principal normal 1}
Assuming \eqref{u is small} holds, if $\varepsilon$ is sufficiently small then we have
    \[|\tilde \nu_1 - (\bar \nu - \tau u' \bar g^{1i}\partial_i \bar N)| \leq C\tau^2 r_0^{-2}(|u|^2 + |u'|^2) + C\tau r_0^{-1}(|u| + |u'| + |u''|)Y^-.\]
\end{lemma}
\begin{proof}
We have 
    \[\partial_\tau \tilde \nu_1 = (\partial_\tau \tilde \nu_1)^\perp + (\partial_\tau \tilde \nu_1)^\top = \sum_{\alpha \geq 2} \langle \partial_\tau \tilde \nu_1, \tilde\nu_\alpha\rangle \tilde \nu_\alpha - \langle \tilde \nu_1, \partial_\tau \partial_i \tilde N\rangle \tilde g^{ij}\partial_j \tilde N.\]
Using Lemma~\ref{deformation of metric} and Lemma~\ref{normal with deformation A} we get 
    \begin{align*}
    \bigg|\sum_{\alpha \geq 2} \langle \partial_\tau \tilde \nu_1, \tilde\nu_\alpha\rangle \tilde \nu_\alpha\bigg| & \leq Cr_0^{-2}(|u| + |u'|)\sum_{\alpha \geq 2} \sum_{i,j}|\langle\nu_\alpha, \tilde A_{ij}\rangle| + Cr_0^{-1}(|u| + |u'| + |u''|)\sum_{\alpha \geq 2} |\langle\tilde\nu_\alpha,\bar\nu\rangle|\\
    &+ Cr_0^{-1}(|u| + |u'|) Y^-.
    \end{align*}
Since $|\langle\nu_\alpha, \tilde A_{ij}\rangle| \leq Cr_0 Y^-$ and $|\langle \tilde \nu_\alpha, \bar \nu\rangle| \leq CY^-$ by Lemma~\ref{Gronwall},
    \[\partial_\tau \tilde \nu_1 = -\langle \tilde \nu_1, \partial_\tau \partial_i \tilde N\rangle \tilde g^{ij}\partial_j \tilde N + O(r_0^{-1}(|u| + |u'| + |u''|))Y^-.\]
Next we use Lemma~\ref{lem_derivatives of tilde N} to bound the final term on the right by $Cr_0^{-1}(|u| + |u'|)$ and integrate with respect to $\tau$ to obtain 
    \[|\tilde \nu_1 - \bar \nu| \leq C\tau r_0^{-1} (|u| + |u'|) + C\tau r_0^{-1}(|u| + |u'| + |u''|)Y^-.\]
Using this estimate we see that
    \begin{align*}
    - \langle \tilde \nu_1&, \partial_\tau \partial_i \tilde N\rangle \tilde g^{ij}\partial_j \tilde N = -\delta^1_i u' \langle \tilde \nu_1, \bar \nu\rangle \tilde g^{ij}\partial_j \tilde N + u \bar h_i^k \langle\tilde\nu_1, \partial_k \bar N\rangle \tilde g^{ij}\partial_j \tilde N + O(r_0^{-1} |u|)Y^-,
    \end{align*}
which can be fed back into our expression for $\partial_\tau \tilde \nu_1$ to give
    \begin{align*}
    \partial_\tau \tilde \nu_1 &= -u'\tilde g^{1j}\partial_j \tilde N + O(\tau r_0^{-2}(|u|^2+|u'|^2))+ O(r_0^{-1}(|u| + |u'| + |u''|))Y^-\\
    &= -u'\bar g^{1j}\partial_j  \bar N + O(\tau r_0^{-2}(|u|^2+|u'|^2))+ O(r_0^{-1}(|u| + |u'| + |u''|))Y^-,
    \end{align*}
where we have used Lemma~\ref{lem_derivatives of tilde N} and Lemma~\ref{deformation of metric}. Our claim follows after integrating.
\end{proof}

Using $\tilde g_{ij} = \langle \partial_i \tilde N, \partial_j \tilde N\rangle$ and $\tilde h_{ij} = \langle \partial_i \partial_j \tilde N,\tilde \nu_1\rangle$ we repeat the arguments in \cite{HuSi09} to obtain the following statements. Here Lemma~\ref{deformed principal normal 1} replaces \cite[(iv) of Corollary~3.16]{HuSi09}. 

\begin{lemma}\label{bent neck metric}
Assuming \eqref{u is small} holds, if $\varepsilon$ is sufficiently small then
\begin{align*}
\tilde g_{ij} &= \bar g_{ij} -2\tau u \bar h_{ij} + O(\tau^2 (|u|^2 + |u'|^2)),\\
\tilde g^{ij} &= \bar g^{ij} + 2\tau u \bar h^{ij} + O(\tau^2 r_0^{-4}(|u|^2 + |u'|^2)),\\
\tilde h_{ij} &= \bar h_{ij} + \tau \delta_i^1\delta_j^1 u'' - \tau u\bar h_j^k\bar h_{ik} + O(\e\tau|u'|)\\
&\qquad + O(\tau^2 r_0^{-1} (|u|^2 + |u'|^2 + |u''|^2))\\
&\qquad +O(\tau(|u| + |u'| + |u''|))Y^-.
\end{align*}
\end{lemma}

We now make use of the specific choice $u(z) = r_0 \exp(-\frac{B}{z-\Lambda})$ as in \cite{HuSi09, Hamilton1997}. Namely, we appeal to the fact that for any $\delta > 0$, if $B \geq 12\Lambda$ is large enough then for $z \in [\Lambda, 5\Lambda]$ we have
    \begin{equation}\label{u derivatives}
    u \geq 0, \;\;\;\;\; u' \geq 0, \;\;\;\;\; u'' \geq 0, \;\;\;\;\; u(z) + u'(z)\leq \delta u''(z), \;\;\;\;\; u''(z) \leq \delta r_0.
    \end{equation}
In addition, we will make use of the inequality
    \begin{equation}\label{u derivatives 2}
    |u'(z)|^2 \leq \frac{100}{B^2} r_0 u(z).
    \end{equation}

\begin{theorem}\label{deformed A and H squared 2}
Given any $\theta > 0$, if $B$ is large enough depending on $\Lambda$, and $\varepsilon$ and $\tau_0$ are both small depending on $n, \Lambda, \theta$ then for every $\tau \in [0,\tau_0]$ and $z \in [\Lambda, 5\Lambda]$ the submanifold $\tilde N$ satisfies:
\begin{enumerate}
    \item[(i)] $\big||\tilde A|^2_{\tilde g} - (|\bar A|^2_{\bar g} + 2\tau\bar h^{11}u'' + 2\tau u \bar h_j^l \bar h_l^i \bar h_i^j)\big| \leq \theta \tau r_0^{-3} u''$ and
    \item[(ii)] $\big| |\tilde H|^2 - (|\bar H|^2 +  2\tau\bar g^{11} u''|\bar H| + 2\tau u \bar h^l_i \bar h^i_l|\bar H|)\big| \leq \theta \tau r_0^{-3} u''$.
\end{enumerate}
\end{theorem}
\begin{proof}
We choose $B$ large enough so that \eqref{u derivatives} holds with $\delta = 1$. Recall that $Y^- \leq C\varepsilon$. From Lemma~\ref{bent neck metric} we obtain
    \begin{align*}
    |\tilde h|^2_{\tilde g} &= |\bar h|^2_{\bar g} + 2\tau\bar h^{11}u'' + 2\tau u \bar h_j^l \bar h_l^i \bar h_i^j +O(\e \tau r_0^{-3} u'') + O(\tau^2 r_0^{-4} |u''|^2).
    \end{align*}
From this the claim concerning $|\tilde A|^2_{\tilde g}$ follows since
    \[|\tilde A|^2_{\tilde g} = |\tilde h|_{\tilde g}^2 + |\tilde A^-|^2_{\tilde g} = |\tilde h|_{\tilde g}^2 + |\bar A^-|^2_{\bar g} + + O(\varepsilon \tau r_0^{-3}u'')\]
by Theorem~\ref{A- norm in deformed metric}. The claim concerning $|\tilde H|^2$ is immediate from Lemma~\ref{bent neck metric}.
\end{proof}
    
\begin{corollary}\label{deformed pinching quant}
For $B$ as in Theorem~\ref{deformed A and H squared 2}, if  $\varepsilon$ and $\tau_0$ are small enough depending on $n$ and $\Lambda$, then for $\tau \in [0,\tau_0]$ the submanifold $\tilde N$ has the following property: If $b$ is any constant in the range $[\frac{1}{n-1}, 1]$ then the inequalities
    \[b|\tilde H|^2 - |\tilde A|^2_{\tilde g} \geq b|\bar H|^2 - |\bar A|^2_{\bar g} + \tau r_0^{-3} u'' \qquad \text{and} \qquad |\tilde H|^2 \geq |\bar H|^2 + \tau r_0^{-3}u''\]
hold in $[\Lambda, 5\Lambda]$. 
\end{corollary}
\begin{proof}
Appealing to Theorem~\ref{deformed A and H squared 2} with $\theta = \frac{b(n-1)}{4}$, we may assume $\e$ and $\tau_0$ are such that for $\tau \in [0, \tau_0]$ we have
    \begin{align*}
    b|\tilde H|^2 - |\tilde A|^2_{\tilde g} &\geq  b|\bar H|^2 - |\bar A|^2_{\bar g} + 2\tau u''(b\bar g^{11} |\bar H| - \bar h^{11})\\
    &+2\tau u (b\bar h_i^j \bar h_j^i |\bar H| - \bar h_j^l\bar h^j_l\bar h_i^j) - \frac{b(n-1)}{4}\tau r_0^{-3} u''.
    \end{align*}
Using Lemma~\ref{lem_coord derivs blended neck} to compare $\bar h$ with $h$ we see that $b\bar g^{11} |\bar H| - \bar h^{11} \geq (b(n-1) - C\e)r_0^{-3}$ is positive when $\varepsilon$ is small, and $b\bar h_i^j\bar h_j^i |\bar H| - \bar h_i^j\bar h_j^l \bar h^j_l \geq (b(n-1)^2 - (n-1) - C\e) r_0^{-3}$. Now using $b \geq \frac{1}{n-1}$ and $0 \leq u \leq u''$ it follows that
    \begin{align*}
    b|\tilde H|^2 - |\tilde A|^2_{\tilde g} &\geq  b|\bar H|^2 - |\bar A|^2_{\bar g} +2\tau(b(n-1) - C\e - b(n-1)/4) r_0^{-3} u''.
    \end{align*}
If $\varepsilon$ is small enough the coefficient in brackets on the right-hand side is at least $\frac{b(n-1)}{2} \geq \frac{1}{2}$. This establishes the first claim. The second claim is analogous.
\end{proof}

\subsection{Blending with a convex axisymmetric cap}

We now consider Step~\ref{surgery blend 2} of the surgery procedure, in which we interpolate between $\tilde N$ and $\tilde C_{z_0}$ over the region $[3\Lambda, 4\Lambda]$ to obtain $\hat N$.

\begin{theorem}\label{geom quants N hat}
For $B$ as in Theorem~\ref{deformed A and H squared 2} we can choose $\tau_0$ depending on $n$ and $\Lambda$ and then restrict $\varepsilon \ll \tau_0$ depending on $n$ and $\Lambda$ so that the following holds. Taking $\tau = \tau_0$ in the definition of $\tilde N$, concerning $\hat N$ in $[0,5\Lambda]$ we have:
\begin{enumerate}[itemsep=0.1cm]
    \item[(i)] $b|\hat H|^2 - |\hat A|^2 \geq b|H|^2 - |A|^2$ for every constant $b \in [\frac{1}{n-1}, 1]$,
    \item[(ii)] $|\hat H|^2 \geq |H|^2$,
    \item[(iii)] $|\hat A^-|^2 = |A^-|^2 + O(|\Theta_2^-|^2 r_0^{-2})$,
    \item[(iv)] $\sqrt{\det \hat g} \leq \sqrt{\det g}$.
\end{enumerate}
Moreover, these inequalities are strict in the region $[\Lambda, 5\Lambda]$.
\end{theorem}
\begin{proof}
In the region $[0,\Lambda]$, $\hat N$ agrees with $N$ so there is nothing to prove. In the region $[\Lambda, 2\Lambda]$ the neck is not affected by either of the blending steps, i.e.\ $\hat N = \tilde N$ and $\bar N = N$, so in this region claims (i) and (ii) are immediate consequences of Corollary~\ref{deformed pinching quant}. In the region $[2\Lambda, 5\Lambda]$ we know that $b|\bar H|^2 - |\bar A|$ and $b|H|^2 - |A|^2$ differ by errors which are at worst of order $O(r_0^{-2}\varepsilon)$ because of Proposition~\ref{geom quants blended neck 2}. Moreover, when we blend $\tilde N$ with $\tilde C_{z_0}$ to obtain $\hat N$, this only introduces further errors of order $O(r_0^{-2} \varepsilon)$. Therefore, in the region $[2\Lambda, 5\Lambda]$, by Corollary~\ref{deformed pinching quant} we have 
    \[b|\hat H|^2 - |\hat A|^2 \geq b|H|^2 - |A|^2 + \tau_0 r_0^{-3} u'' - C\varepsilon r_0^{-2}.\]
Since $B$ has now been fixed depending only on $\Lambda$, because of \eqref{u derivatives} we have a positive lower bound for $r_0^{-1} u''$ in $[2\Lambda, 5\Lambda]$ which depends only on $\Lambda$. Therefore, in this region we may take $\varepsilon \ll \tau_0$ small enough to ensure 
    \[b|\hat H|^2 - |\hat A|^2 \geq b|H|^2 - |A|^2 + \frac{\tau_0 r_0^{-3} u'' }{2} > b|H|^2 - |A|^2.\]
This completes the proof of claims (i) and (ii).

We now turn to (iii). Since $\bar N$ is a hypersurface in $[3\Lambda, 5\Lambda]$, both $\tilde N$ and $\hat N$ are also hypersurfaces in this region. Therefore $|\hat A^-|^2 = 0$ identically in $[3\Lambda, 5\Lambda]$, meaning we can restrict attention to $[\Lambda, 3\Lambda]$, where we have $\hat N = \tilde N$. (This is why we included Step~\ref{surgery blend} in our procedure---if we had not done so, then further analysis would be required to estimate $|\hat A^-|^2$ in terms of $|\tilde A^-|^2$ in the region $[3\Lambda, 5\Lambda]$.) Using Theorem~\ref{A- norm in deformed metric} and Proposition~\ref{geom quants blended neck 2}, in $[\Lambda, 3\Lambda]$ we obtain
    \begin{align*}
    |\hat A^-|^2 &= |\tilde A^-|^2 = |\bar A^-|^2 + O(|Y^-|^2 r_0^{-2}) = |A^-|^2 + O(|\Theta_1^-|^2r_0^{-2}) + O(|Y^-|^2 r_0^{-2}).
    \end{align*}
One checks directly using Lemma~\ref{estimate for Y} and Proposition~\ref{geom quants blended neck 2} that $|Y^-|^2 \leq C|\Theta_2^-|^2$. Claim (iii) follows. 

To obtain the final claim we first combine Lemma~\ref{bent neck metric} with Proposition~\ref{geom quants blended neck} to see that 
    \begin{align*}
    \sqrt{\det(\tilde g)} &= \sqrt{\det(\bar g)}\Big[1 - \tau u |\bar H| + O(\tau^2 r_0^{-2}(|u|^2 + |u'|^2))\Big]\\
    &= \sqrt{\det(\bar g)}\Big[1 - \tau (n-1)r_0^{-1} u + O(\tau^2 r_0^{-2}(|u|^2 + |u'|^2)) + O(\varepsilon \tau r_0^{-1} u)\Big]\\
    &=\sqrt{\det(\bar g)}\Big[1 - \tau (n-1)r_0^{-1} u  + O(\tau^2 r_0^{-1}u) + O(\varepsilon \tau r_0^{-1} u)\Big],
    \end{align*}
where we have used \eqref{u derivatives 2} and $B \geq 12\Lambda$. Consequently, by choosing $\tau_0$ sufficiently small and then taking $\varepsilon \ll \tau_0$, we ensure that for $\tau = \tau_0$ we have
    \[\sqrt{\det(\tilde g)} = \sqrt{\det(\bar g)}\left(1 - \frac{\tau_0 (n-1)r_0^{-1} u}{2}\right)\]
in the region $[\Lambda, 5\Lambda]$. Given that $\hat g = \tilde g$ and $\bar g = g$ in the region $[\Lambda, 2\Lambda]$, we read off that (iv) holds there. In $[2\Lambda, 5\Lambda]$ we have a positive lower bound for $r_0^{-1}u$ depending only on $\Lambda$, so by making $\varepsilon$ a bit smaller, using Proposition~\ref{geom quants blended neck} we deduce 
    \begin{align*}
    \sqrt{\det(\hat g)} &= \sqrt{\det(\tilde g)}(1+O(\varepsilon))= \sqrt{\det(g)}\left(1 - \frac{\tau_0 (n-1)r_0^{-1} u}{2} + O(\varepsilon)\right) < \sqrt{\det(g)}.
    \end{align*}
This completes the proof. 
\end{proof}

\subsection{Preservation of the pinching class} We can now show that the class of submanifolds $ \mathcal C_{n,m}(R,\alpha)$ is preserved under standard surgery for appropriate choices of parameters. 

\begin{theorem}\label{thm_invariant}
Suppose $\mathcal{M} \subset \mathbb{R}^{n+m}$ is a closed submanifold in the class $\mathcal C_{n,m}(R,\alpha)$ and let $N : \mathbb{S}^{n-1} \times [-7\Lambda,7\Lambda] \to \mathcal{M}$ be a normal $(\varepsilon,k)$-cylindrical submanifold neck in $\mathcal{M}$. If $k \geq 4$ and the parameters $\Lambda \geq 10$, $(\tau, B)$ and $\varepsilon$ are in suitable ranges (depending only on $n$ and $\Lambda$), then the submanifold obtained from $\mathcal{M}$ by performing standard surgery on $N$ is still in the class $\mathcal C_{n,m}(R,\alpha)$. 
\end{theorem}
\begin{proof}
We need to show that the surgery preserves the inequalities
    \[|A|^2 + \alpha_2 R^{-2} - c_n |H|^2 \leq - \alpha_1 |H|^2, \qquad \mu(\mathcal{M}^n)\leq \alpha_3 R^n.\]
Recall that the regions $[-7\Lambda, -4\Lambda]$ and $[4\Lambda, 7\Lambda]$ are unchanged by surgery. In the regions $[-4\Lambda, -\Lambda]$ and $[\Lambda, 4\Lambda]$ the claim follows from Theorem~\ref{geom quants N hat} (note we may assume without loss of generality that $c_n - \alpha_1 > \frac{1}{n-1}$). For the regions $[-\Lambda,0]$ and $[0,\Lambda]$ we have to consider the final step of the surgery where $u$ is changed to $\hat u$. Assuming $\varepsilon$ is small enough (depending only on $n$ and $\Lambda$) this can be done so that $\frac{1}{n-1}|H|^2-|A|^2$ and $|H|^2$ both increase pointwise, hence the pinching estimate is preserved. Moreover, it is easy to see that the resulting convex caps will have less area than the cylindrical pieces they replace if $\Lambda$ is larger than some universal constant. 
\end{proof}

We henceforth consider $\Lambda \geq 10$ and $(\tau, B)$ to be fixed, and assume $\varepsilon_0 = \varepsilon_0(n)$ is small enough, so that Theorem~\ref{geom quants N hat} and Theorem~\ref{thm_invariant} are in effect whenever $k \geq k_0 := 4$ and $0 < \varepsilon \leq \varepsilon_0$. Whenever we use the standard surgery procedure, the parameters are assumed to be chosen in this way. In particular, every time we talk about a mean curvature flow with surgery, all surgeries are performed using these parameters. The estimates for $A$ and $g$ from Theorem~\ref{geom quants N hat} are essential for proving that our a priori estimates for quadratically pinched mean curvature flows are preserved under surgeries. We have the following immediate consequence of Theorem~\ref{thm_invariant}:

\begin{corollary}\label{pinching class under flow with surgeries} There exists a positive constant $\varepsilon_0 = \varepsilon_0(n)$ with the following property. Let $\{\mathcal{M}_t\}_{t\in[0,T]}$ be a mean curvature flow with surgery such that $\mathcal{M}_0$ is of class $\mathcal{C}_{n,m}(R,\alpha)$ and satisfies $|A|^2 \leq R^{-2}$. If all surgeries are performed on $(\varepsilon,k)$-cylindrical submanifold necks with $k \geq k_0$ and $\varepsilon \leq \varepsilon_0$ then $\mathcal{M}_t$ is in $\mathcal{C}_{n,m}(R,\alpha)$ for every $t \in [0,T]$.
\end{corollary}

To conclude this section let us note that at the topological level standard surgery amounts to a reverse connected sum, in the following manner. 

\begin{theorem}\label{thm_reverse}
Suppose the standard surgery procedure is performed on a normal $(\e,k)$-cylindrical submanifold neck in some closed, connected, immersed submanifold $\mathcal{M}$, resulting in a new submanifold $\widetilde{\mathcal{M}}$. If $\widetilde{\mathcal{M}}$ is connected then $\mathcal{M}$ is diffeomorphic either to $\widetilde{\mathcal{M}} \# (\mathbb{S}^{n-1}\times\mathbb{S}^1)$ or $\widetilde{\mathcal{M}} \# (\mathbb{S}^{n-1}\simtimes\mathbb{S}^1)$. If $\widetilde{\mathcal{M}}$ is disconnected with two components $\widetilde{\mathcal{M}}_1$ and $\widetilde{\mathcal{M}}_2$ then $\mathcal{M}$ is diffeomorphic to $\widetilde{\mathcal{M}}_1 \# \widetilde{\mathcal{M}}_2$. In particular, if $\widetilde{\mathcal{M}}$ is disconnected and $ \widetilde{\mathcal{M}}_2$ is diffeomorphic to $\mbb S^n$ then $ \widetilde{\mathcal{M}}_1$ is diffeomorphic to $\mathcal{M}$.
\end{theorem}
\begin{proof}
This follows from the surgery construction: The two open $n$-balls attached by the surgery are diffeomorphic to the standard ball, and the collar regions $\mathbb{S}^{n-1}\times(-5\Lambda, 0)$ and $\mathbb{S}^{n-1}\times(0,5\Lambda)$ in $\widetilde{\mathcal{M}}$ as well as the original neck $\mathbb{S}^{n-1}\times(-5\Lambda, 5\Lambda)$ in $\mathcal{M}$ are all diffeomorphic to the standard cylinder. When $\widetilde{\mathcal{M}}$ is connected, the two different cases for $\mathcal M$ arise as follows. Choose a region $U \subset \mathcal{M} \setminus \mathbb{S}^{n-1}\times [-4\Lambda, 4\Lambda]$ which contains both $\mathbb{S}^{n-1}\times[-5\Lambda, -4\Lambda)$ and $\mathbb{S}^{n-1}\times(4\Lambda,5\Lambda]$. Since $\widetilde{\mathcal{M}}$ is connected we can assume $U$ is diffeomorphic to a twice punctured $n$-ball and hence is orientable. If the orientations induced on the two collars $\mathbb{S}^{n-1}\times[-5\Lambda, -4\Lambda)$ and $\mathbb{S}^{n-1}\times(4\Lambda, 5\Lambda]$ by the neck parameterization are compatible with some orientation of $U$, then we have $\widetilde{\mathcal{M}}\# (\mathbb{S}^{n-1}\times\mathbb{S}^1)$. Otherwise, if one of these always disagrees with any orientation of $U$, then we have $\widetilde{\mathcal{M}} \# (\mathbb{S}^{n-1}\simtimes\mathbb{S}^1)$.  
\end{proof}


\section{Gradient estimates}\label{sec_estimates}

We prove a scale-invariant gradient estimate for the second fundamental form which is pointwise in nature. Standard parabolic estimates instead bound $\nabla A$ by the maximum of the curvature over a region of spacetime. The gradient estimate applies to the full range of quadratic pinching conditions shown to be preserved in \cite{Andrews2010}. An analogous gradient estimate played a central role in the hypersurface setting \cite{HuSi09}. Our proof in higher codimensions turns out to be remarkably similar given the far more complicated evolution equation for the second fundamental form. In one respect our proof is even simpler---we do not rely on a cylindrical estimate but rather deduce the gradient estimate directly from quadratic pinching. 

\begin{theorem}\label{thm_grad est}
Let $\{\mathcal{M}_t\}_{t \in [0,T)}$ be a closed mean curvature flow in $\mathbb{R}^{n+m}$ of dimension $n \geq 2$. Suppose $|A|^2 \leq R^{-2}$ and $|A|^2  + \alpha_2 R^{-2}\leq (\frac{4}{3n} - \alpha_1)|H|^2$ hold everywhere in $\mathcal{M}_0$ for some constants $\alpha_1,\alpha_2 > 0$. We then have
    \begin{align*}
    |\nabla A|^2 \leq \gamma |H|^4
    \end{align*}
at every point of $\mathcal{M}_t$ for all $ t\in [R^2/20, T)$, where $\gamma = \gamma(n, \alpha_1, \alpha_2)$.
\end{theorem}
\begin{proof}
We apply the maximum principle to $|\nabla A|^2/g^2$, where
    \[g:= \left(\frac{4}{3n}-\frac{\alpha_1}{2}\right)|H|^2 - |A|^2.\] 
Since the inequality $|A|^2  + \alpha_2R^{-2}\leq (\frac{4}{3n}-\alpha_1)|H|^2$ is preserved by the flow we have $g \geq \frac{\alpha_1}{2}|H|^2 + \alpha_2 R^{-2}$ for all $t \in [0,T)$.

First we bound $|\nabla A|^2$ at small times. Using $(\partial_t - \Delta)|A|^2 \leq 10|A|^4$, an ODE comparison argument shows that $|A|^2 \leq 2R^{-2}$ holds on $\mathcal{M}_t$ for $t \leq R^2/20$. Standard interior derivative estimates (see \cite[Proposition~4.8]{Andrews2010} and \cite{Ecker1991}) imply that, while the upper bound $|A|^2 \leq 2R^{-2}$ remains valid, we have $|\nabla A|^2 \leq K_0R^{-4}(1+R^2/t)$ for a constant $K_0 = K_0(n)$. Therefore, at time $t = R^2/20$ we have $|\nabla A|^2 \leq 21K_0 \alpha_2^{-2} g^2$.

Next we recall from Lemma~\ref{lem_pinching refinement} and Lemma~\ref{lem_Kato} the inequality 
    \[(\partial_t - \Delta) g \geq 2\left(|\nabla A|^2 - \left(\frac{4}{3n}-\frac{\alpha_1}{2}\right)|\nabla H|^2\right) \geq 2\left(1-\frac{4(n+2)}{9n}\right)|\nabla A|^2.\]
In \cite{Andrews2010} it was shown that 
    \[(\partial_t - \Delta)|\nabla A|^2 \leq-2 |\nabla^2 A|^2+K_1|A|^2 |\nabla A|^2\]
for a constant $K_1 = K_1(n)$. Combining these two inequalities gives 
    \begin{align*}
    (\partial_t - \Delta) \frac{|\nabla A|^2}{g} &= \frac{1}{g}(\partial_t - \Delta)|\nabla A|^2 - \frac{|\nabla A|^2}{g^2}(\partial_t - \Delta) g + \frac{2}{g}\left\langle\nabla g, \nabla \frac{|\nabla A|^2}{g}\right\rangle\\
    &\leq -2\frac{|\nabla^2 A|^2}{g} + K_1|A|^2\frac{|\nabla A|^2}{g} - 2\left(1-\frac{4(n+2)}{9n}\right)\frac{|\nabla A|^4}{g^2} +\frac{2}{g}\left\langle\nabla g, \nabla \frac{|\nabla A|^2}{g}\right\rangle.
    \end{align*}
Using $\langle \nabla g, \nabla |\nabla A|^2\rangle \leq 2|\nabla g||\nabla A||\nabla^2 A|$ and Young's inequality we obtain 
    \[-2\frac{|\nabla^2 A|^2}{g} +\frac{2}{g}\left\langle\nabla g, \nabla \frac{|\nabla A|^2}{g}\right\rangle = -2\frac{|\nabla^2 A|^2}{g} + \frac{2}{g^2} \langle \nabla g, \nabla |\nabla A|^2\rangle - \frac{2}{g^3} |\nabla A|^2 |\nabla g|^2 \leq 0\]
and hence 
    \begin{align*}
    (\partial_t - \Delta) \frac{|\nabla A|^2}{g} &\leq K_1|A|^2\frac{|\nabla A|^2}{g} - 2\left(1-\frac{4(n+2)}{9n}\right)\frac{|\nabla A|^4}{g^2}.
    \end{align*}
Dividing once more by $g$, we see that
    \begin{align*}
    (\partial_t - \Delta) \frac{|\nabla A|^2}{g^2} &\leq \frac{|\nabla A|^2}{g}\left[K_1\frac{|A|^2}{g} - 2\left(1-\frac{4(n+2)}{9n}\right)\frac{|\nabla A|^2}{g^2}\right] + \frac{2}{g}\left\langle\nabla g, \nabla \frac{|\nabla A|^2}{g^2}\right\rangle.
    \end{align*}
By appealing to the parabolic maximum principle, taking into account $|A|^2\leq 2\alpha_1^{-1} g$, we deduce the estimate
    \[\max_{\mathcal{M}_t}\frac{|\nabla A|^2}{g^2}\leq \max\left\{21K_0\alpha_2^{-2}, K_1\alpha_1^{-1}\left(1-\frac{4(n+2)}{9n}\right)^{-1} \right\}\]
for $t \geq R^2/20$. Since $g \leq |H|^2$, the claim follows. 
\end{proof}

Theorem~\ref{thm_grad est} applies to smooth flows. We now show that it also holds for flows with surgery. 

\begin{theorem}\label{thm_grad est with surgeries}
Let $\{\mathcal{M}_t\}_{t\in[0,T]}$ be a mean curvature flow with surgery of class $\mathcal C_{n,m}(R,\alpha)$ such that $|A|^2 \leq R^{-2}$ initially. If the surgery parameters of Section~\ref{sec_standard surgery} are in suitable ranges, then there exists a constant $ \gamma_1 = \gamma_1(n, \alpha)$ such that
	\begin{align*}
	|\nabla A|^2 \leq \gamma_1 |H|^4
	\end{align*}
at every point of $\mathcal{M}_t$ for $t \in [R^2/20, T]$. 
\end{theorem}
\begin{proof}
For $g := (\frac{4}{3n} - \frac{\alpha_1}{2})|H|^2 - |A|^2$ the same parabolic maximum principle argument used to prove Theorem~\ref{thm_grad est} shows that
    \[\max_{\mathcal{M}_t}\frac{|\nabla A|^2}{g^2}\leq \max\left\{21K_0\alpha_2^{-2}, K_1\alpha_1^{-1}\left(1-\frac{4(n+2)}{9n}\right)^{-1} \right\}\]
for all $t \geq R^2/20$ not exceeding the first surgery time. Now we discuss how the quantity $|\nabla A|^2/g^2$ behaves under surgery. On an $(\e, k)$-neck with $k \geq 1$, if $\varepsilon$ is small then $|\nabla A|^2/|H|^4 = O(\varepsilon)$. Therefore, for any choice of the transition functions $\psi$ and $\varphi$ and convex cap in the standard surgery procedure of Section~\ref{sec_standard surgery}, by assuming $\varepsilon$ is sufficiently small we can guarantee that $|\nabla A|^2 \leq \mu |H|^4$ in regions affected by surgery for a constant $\mu = \mu(n)$. Since $g \geq \frac{\alpha_1}{2}|H|^2$ this implies that $|\nabla A|^2 \leq 4\mu \alpha_1^{-2} g^2$ in regions affected by surgery. Applying the parabolic maximum principle as in Theorem~\ref{thm_grad est} we see that
    \[\max_{\mathcal{M}_t}\frac{|\nabla A|^2}{g^2}\leq \max\left\{21K_0\alpha_2^{-2}, K_1\alpha_1^{-1}\left(1-\frac{4(n+2)}{9n}\right)^{-1},4\mu \alpha_1^{-2}\right\}\]
holds up to the next surgery time. We may repeat this argument, iterating over the surgery times, to see that the same estimate holds for all $t \in [R^2/20, T]$. 
\end{proof}

Using the above theorem, we get estimates on higher-order derivatives up to order $k_0$ (the surgery regularity parameter). Strictly speaking, the estimates do not hold at surgery times as the flow is not smooth there, but by taking appropriate limits from below and above we recover the desired estimates

\begin{theorem}\label{thm_highergradest}
Let $\{\mathcal{M}_t\}_{t \in [0,T]}$, be a mean curvature flow with surgery of class $\mathcal C_{n,m}(R,\alpha)$ such that $|A|^2 \leq R^{-2}$ initially. If the surgery parameters of Section~\ref{sec_standard surgery} are in suitable ranges, then there exists a constant $\gamma = \gamma(n,\alpha)$ such that
	\begin{align*}
	|\nabla^\ell_t \nabla^k A|^2 \leq \gamma |H|^{4\ell +2 k + 2}
	\end{align*}
holds on $\mathcal M_t$ for all times, provided $2\ell + k \leq k_0$.
\end{theorem}
\begin{proof}
It suffices to prove the estimates for spatial derivatives, since the time derivatives can then be bounded by differentiating $\nabla_t A = \Delta A + A \ast A \ast A$. With Theorem~\ref{thm_grad est with surgeries} in place, the arguments needed to bound $|\nabla^k A|^2/|H|^{2+2k}$ are essentially contained in \cite{HuSi09}, so let us only provide a sketch in the case $k = 2$. We apply the maximum principle to
    \[Z := \frac{|\nabla^2 A|^2}{|H|^5} + \kappa_1 \frac{|\nabla A|^2}{|H|^3} - \kappa_2\sqrt{c_n|H|^2-|A|^2}\]
for appropriate constants $\kappa_1$ and $\kappa_2$. In \cite{Andrews2010} it was shown that
    \[(\partial_t - \Delta)|\nabla^2 A|^2 \leq - 2|\nabla^3 A|^2 + C(n)|A|^2|\nabla^2 A|^2 + C(n)|\nabla A|^4,\]
using \eqref{eqn:H2} one obtains
    \[(\partial_t - \Delta) |H|^p \geq -p(p-1)|H|^{p-2}|\nabla H|^2\]
for every $p > 2$, and because of Lemma~\ref{lem_pinching refinement} and Lemma~\ref{lem_Kato} we have
    \begin{align*}
    (\partial_t - \Delta)\sqrt{c_n|H|^2-|A|^2} &\geq  |h|^2\sqrt{c_n|H|^2-|A|^2} \geq \frac{\sqrt{\alpha_1}}{n} |H|^3.
    \end{align*}
Combining these three inequalities and using Young's inequality one finds that 
    \begin{align*}(\partial_t - \Delta) Z &\leq  \Big(C(1+\gamma_1) - \kappa_1\Big)\frac{|\nabla^2 A|^2}{|H|^3} + \bigg(C\gamma_1 + C\kappa_1\gamma_1(1+\gamma_1) - \frac{\kappa_2\sqrt{\alpha_1}}{n}\bigg)|H|^3
    \end{align*}
for some $C = C(n)$. We may choose $\kappa_1$ and $\kappa_2$ depending only on $n$ and $\alpha$ so that the right-hand side is nonpositive for $t \geq R^2/20$. Using the fact that $|A|$ is bounded for $t \leq R^2/20$ and the interior estimates from \cite{Andrews2010} we get $Z \leq C(n,\alpha)R^{-1}$ at time $t = R^2/20$. Moreover, on an $(\varepsilon,k_0)$-neck with $\varepsilon$ small we have $|\nabla^2 A|^2/|H|^6 + |\nabla A|^2/|H|^4 = O(\varepsilon)$, so if $\varepsilon$ is sufficiently small (depending only on $n$) then we can guarantee
    \[Z \leq C(n)(1+\kappa_1)|H| - \kappa_2\sqrt{\alpha_1} |H|\]
in regions affected by surgery, and by choosing $\kappa_2$ a bit larger we ensure that $Z \leq 0$ in regions affected by surgery. We can now iterate the maximum principle between surgery times as in Theorem~\ref{thm_grad est with surgeries} to get $Z \leq C(n,\alpha)R^{-1}$ for all $t \in [R^2/20, T]$. The claim follows since we can bound $R^{-1} \leq \alpha_2^{-1} |H|$.
\end{proof}

We have the following corollary, which will be used extensively in our analysis.

\begin{corollary}\label{cor_gradient}
In the setting of Theorem~\ref{thm_grad est with surgeries}, there exists a constant $c^\# = c^\#(n,\alpha)$ such that $|\nabla A| \leq n^{-1/2} c^\# |H|^2$ (implying $|\nabla H|\leq c^\# |H|^2$) and $|\nabla_t H|\leq c^\# |H|^3$. 
\end{corollary}

The following lemma allows us to compare the mean curvature at different points. It is proven by integrating $|\nabla H| \leq c^{\#}|H|^2$ along geodesics as in \cite{HuSi09}.

\begin{lemma}\label{lem_localHarnack}
Let $\mathcal{M}$ be an immersed $n$-submanifold of $\mathbb{R}^{n+m}$ and suppose the estimate $|\nabla H(p)| \leq c^{\#}|H(p)|^2$ holds for all $p\in \mathcal M$. For all $p, q \in \mathcal{M}$ we have
    \begin{align*}
    |H(q)|\geq \frac{|H(p)|}{1+c^{\#}d(p, q)|H(p)|}.
    \end{align*}
\end{lemma}


\section{The planarity estimate for flows with surgery}\label{sec_planarity}

In \cite{Naff_planarity} Naff showed that for $n \geq 5$ suitably pinched smooth mean curvature flows of arbitrary codimension become asymptotically codimension-one at singularities. This statement follows from a \emph{planarity estimate}, which shows that the components of $A$ orthogonal to $H$ are of lower order whenever the curvature becomes large. Naff's planarity estimate constituted a major breakthrough in the theory. Its proof uses only the parabolic maximum principle.\footnote{See also \cite{LNZ} where the planarity estimate was recently localised in time and used to classify ancient solutions.}

The planarity estimate applies in particular to smooth flows in the class $\mathcal{C}_{n,m}(R,\alpha)$. One of our major contributions in this work is to extend it to flows with surgery. This turns out to be extremely delicate, because the measure of planarity appearing in Naff's estimate deteriorates slightly under each surgery.

We first give a precise statement of the planarity estimate for a smooth mean curvature flow of class $\mathcal{C}_{n,m}(R,\alpha)$. Let us set
    \[Q=c|H|^2 - |A|^2\]
where $c = \frac{3(n+1)}{2n(n+2)}-\frac{1}{2}\min\{\alpha_1, \frac{3(n+1)}{2n(n+2)} -\frac{1}{n-1}\}$ if $n \in \{5,6,7\}$ and $c = \frac{4}{3n}$ if $n \geq 8$. In addition, let us set $\sigma = \min\{\frac{1}{32}, \frac{n(n+2)}{3(n-1)}\alpha_1\}$ if $n \in \{5,6,7\}$ and $\sigma = \frac{1}{5n-8}$ if $n \geq 8$.

\begin{theorem}[Planarity estimate for smooth flows  \cite{Naff_planarity}]\label{thm_planarity smooth} Let $\{\mathcal{M}_t\}_{t\in[0,T)}$ be a smooth mean curvature flow such that $\mathcal{M}_0$ is of class $\mathcal{C}_{n,m}(R,\alpha)$. We then have 
    \[\max_{\mathcal{M}_t} \frac{|A^-|^2}{Q^{1-\sigma}} \leq \max_{\mathcal{M}_0} \frac{|A^-|^2}{Q^{1-\sigma}}\]
for all $t \in [0,T)$. 
\end{theorem}

Since $Q \geq (c-c_n+\alpha_1)|H|^2$ for a flow in $\mathcal{C}_{n,m}(R,\alpha)$, assuming the normalisation $|A|^2 \leq R^{-2}$ at $t = 0$ we have 
    \begin{equation}\label{def C_pl}
    \max_{\mathcal{M}_0} \frac{|A^-|^2}{Q^{1-\sigma}} \leq \frac{n^\sigma c_n}{(c-c_n+\alpha_1)^{1-\sigma}} R^{-2\sigma}  =: C_{\pl} R^{-2\sigma},
    \end{equation}
and hence the planarity estimate implies 
    \begin{equation}\label{planarity estimate}
    |A^-|^2 \leq C_{\pl} R^{-2\sigma} Q^{1-\sigma}
    \end{equation}
on $\mathcal{M}_t$ for all $t \in [0,T)$. The constant $C_{\pl}$ can be bounded from above in terms of $n$ and $\alpha_1$ if $n \in \{5,6,7,8\}$, and purely in terms of $n$ if $n \geq 9$. 

We now prove the planarity estimate for flows with surgery. We need the estimate to hold with uniform constants independently of the number of surgeries; this is essential for the flow with surgery to terminate after finitely many steps. We might hope to proceed by iteratively applying Naff's maximum principle for the ratio $|A^-|^2/Q^{1-\sigma}$ in between surgery times, but this only works if $|A^-|^2/Q^{1-\sigma}$ admits a suitable uniform upper bound across every surgery. Here however there are fundamental issues:
\begin{itemize}
    \item The standard surgery procedure constructed in Section~\ref{sec_standard surgery} does not make every portion of the neck closer to a hypersurface: in some regions $A^-$ inevitably becomes larger. 
    \item Because the ratio appearing in the planarity estimate is not scale-invariant, any small errors introduced under standard surgery are amplified. Indeed, any small increase in $|A^-|^2$ compared with the natural scale $r^{-2}$ leads to an increase in $|A^-|^2/Q^{1-\sigma}$ which is of order $r^{-2\sigma}$ and hence explodes as the curvature of the neck becomes large. 
    \item Even with refined estimates for $|A^-|^2$ which overcome the scaling issue, iterated surgeries could force the ratio $|A^-|^2/Q^{1-\sigma}$ to blow up. This occurs even if each surgery only increases $|A^-|^2$ by a tiny amount relative to $r^{-2-2\sigma}$, since we cannot bound the number of surgeries.
\end{itemize}

To overcome this fundamental obstruction we need to ensure that on every neck where surgery is performed, the ratio $|A^-|^2/Q^{1-\sigma}$ is extremely small to begin with. Then, since the careful estimates of Section~\ref{sec_standard surgery} ensure that it does not increase by too much, it will remain suitably bounded so that we may iterate the maximum principle without our constants exploding. Thus we construct our surgery algorithm so that the following property holds:

\begin{enumerate}[label=(s\arabic*)]
\setcounter{enumi}{-1}
    \item\label{s0} Each time standard surgery is performed on a normal $(\varepsilon,k_0)$-cylindrical submanifold neck $\mathbb{S}^{n-1} \times [-7\Lambda, 7\Lambda]$, at every point of the neck we have
        \[r_0^{2-2\sigma}|A^-|^2 + r_0^{4-2\sigma}|\nabla A^-|^2 + r_0^{6-2\sigma} |\nabla^2 A^-|^2 \leq \delta_0 C_{\pl}  R^{-2\sigma}.\]
    Here, as in Section~\ref{sec_standard surgery}, $r_0$ is the mean radius of the cross-section $z = 0$. 
\end{enumerate}
Provided we fix the constant $\delta_0$ suitably (depending only on $n$), the planarity estimate holds for flows with surgery satisfying \ref{s0}, with the same constants as the smooth flow:

\begin{theorem}\label{thm_planarity surgeries}
There exists a constant $\varepsilon_0 = \varepsilon_0(n)$ with the following property. Let $\{\mathcal{M}_t\}_{t \in [0,T]}$ be a mean curvature flow with surgery such that $\mathcal{M}_0$ is of class $\mathcal{C}_{n,m}(R,\alpha)$ and satisfies $|A|^2 \leq R^{-2}$. Provided that all surgeries are performed on $(\varepsilon,k_0)$-cylindrical submanifold necks with $\varepsilon \in (0, \varepsilon_0)$, and assuming \ref{s0} holds, we have
    \[|A^-|^2 \leq C_{\pl} R^{-2\sigma} Q^{1-\sigma}\]
at every point of $\mathcal M_t$ for all $t \in [0,T]$.
\end{theorem}
\begin{proof}
Up to the first surgery time the estimate follows from Naff's work as in Theorem~\ref{thm_planarity smooth}. To complete the proof we show that, immediately after performing standard surgery with our parameters chosen as in Section~\ref{sec_standard surgery}, the maximum of $|A^-|^2/Q^{1-\sigma}$ does not exceed $C_{\pl}R^{-2\sigma}$. Indeed, if this is the case then we can appeal to Naff's result again to say that this bound will persist until the next surgery time. The claim then follows by iterating this argument across successive surgeries. 

Let us therefore consider our submanifold at a surgery time $t_0$. As usual we use $t_0-$ and $t_0+$ to distinguish between the submanifold immediately before and after surgery. Let $\mathbb{S}^{n-1}\times[-7\Lambda, 7\Lambda]$ be a portion of a normal $(\varepsilon,k_0)$-cylindrical submanifold neck in $\mathcal{M}_{t_0}$ where standard surgery is performed. The regions $[-7\Lambda, -4\Lambda]$ and $[4\Lambda,7\Lambda]$ are unchanged by surgery and can be ignored. By Theorem~\ref{geom quants N hat}, in the regions $[-4\Lambda, -\Lambda]$ and $[\Lambda,4\Lambda]$, we have 
    \[r_0^2 |A^-(\cdot,t_0+)|^2 \leq C |\Theta_2^-(\cdot, t_0-)|^2\]
for some $C = C(n)$ where, by definition,
    \[|\Theta_2^-|^2 \leq \sup_{z \in [-7\Lambda,7\Lambda]} \bigg(r_0^2|A^-|^2 + r_0^4|\nabla A^-|^2 + r_0^6|\nabla^2 A^-|^2\bigg).\]
Assuming \ref{s0} holds, we have 
	\[|\Theta_2^-(\cdot, t_0-)|^2 \leq \delta_0 C_{\pl}R^{-2\sigma}r_0^{2\sigma}\] 
in $[-7\Lambda, 7\Lambda]$ immediately prior to performing surgery. Consequently, in the regions $[-4\Lambda, -\Lambda]$ and $[\Lambda,4\Lambda]$, after surgery we have
    \[|A^-(\cdot, t_0+)|^2 \leq C \delta_0 C_{\pl}  R^{-2\sigma}r_0^{-2+2\sigma} \leq C\delta_0 C_{\pl} R^{-2\sigma} Q^{1-\sigma}\]
in the regions $[-4\Lambda, -\Lambda]$ and $[\Lambda,4\Lambda]$, where $C$ is now larger but still only depends on $n$. Choosing $\delta_0$ small enough to ensure $C\delta_0 \leq 1$, we get
    \[|A^-(\cdot, t_0+)|^2\leq C_{\pl} R^{-2\sigma}Q^{1-\sigma}\]
in the regions $[-4\Lambda, -\Lambda]$ and $[\Lambda,4\Lambda]$.

The final step is to consider the two strictly convex caps attached at the cross-sections $z = -\Lambda$ and $z = \Lambda$. But these are hypersurfaces, so $A^-$ vanishes in these regions. 

Repeating this argument for every surgery, we see that $|A^-|^2 \leq C_{\pl}R^{-2\sigma} Q^{1-\sigma}$ holds globally on the submanifold immediately after surgery, hence we can iterate the maximum principle as described above. 
\end{proof}

\subsection{The planarity improvement theorem} We have seen that the planarity estimate holds for the flow with surgery, provided that every surgery is performed on a neck which is close to being codimension-one in a strong quantitative sense, that is, if \ref{s0} holds. In order to construct our surgery algorithm so that \ref{s0} holds, we require a completely new analytic tool, which we refer to as the planarity improvement theorem. This theorem asserts that on any sufficiently long neck which is allowed to evolve for a long time under the mean curvature flow, the ratio $|A^-|^2/Q^{1-\sigma}$ becomes arbitrarily small at the centre of the neck. It is similar in spirit to the symmetry improvement theorem for codimension-one necks proven by Brendle--Choi \cite{Brendle--Choi_a, Brendle--Choi_b}, but its geometric content is completely different: here the point is not to show that the axial symmetry of the neck is improving, but rather that it is becoming increasingly close to a hypersurface. 

In what follows $Q = c|H|^2 - |A|^2$ where $c \in [\frac{1}{2}(\frac{3(n+1)}{2n(n+2)} + \frac{1}{n-1}), \frac{3(n+1)}{2n(n+2)}]$ if $n \in \{5,6,7\}$ and $c = \frac{4}{3n}$ if $n \geq 8$. 

\begin{theorem}[Planarity improvement]\label{thm_planarity improvement}
Given any $\delta > 0$ there are constants $L = L(n, \delta)$ and $\varepsilon = \varepsilon(n,\delta)$ with the following property. Let $\mathcal{M}_t$ be an $n$-dimensional mean curvature flow in $\mathbb{R}^{n+m}$, where $n \geq 5$, and suppose $(x_0, t_0)$ lies at the centre of an $(\varepsilon, 1, L, L^2)$-shrinking neck. In addition, let $\sigma \in (0, \frac{1}{32}]$ be fixed and set $K$ equal to the maximum value of $|A^-|^2/Q^{1-\sigma}$ over all points in $\hat {\mathcal P}(x_0,t_0, L, L^2)$. We then have
    \[|A^-(x_0,t_0)|^2 \leq \delta K Q(x_0,t_0)^{1-\sigma}.\]
\end{theorem}

To prove the planarity improvement theorem we study the evolution equation for $|A^-|^2$ on a neck and use the maximum principle to identify an upper barrier for $|A^-|^2/Q^{1-\sigma}$. First we recall the following evolution equation for $|A^-|^2$, which was derived by Naff in \cite{Naff_planarity}:
    \begin{align}\label{evolution |A-|^2}
    (\partial_t - \Delta) |A^-|^2 & = 2 \sum_{i,j,k,l}\langle A^-_{ij} , A_{kl}^-\rangle^2 + 2 \sum_{i,j} |A^-_{ik} \otimes A^-_{kj} - A^-_{jk} \otimes A^-_{ki}|^2 + 2 \sum_{i,j} |h_{ik} A^-_{kj} - h_{jk} A^-_{ki} |^2 \notag \\
    &\qquad - 2|\nabla A^-|^2 + 4 \sum_{i,j,k} (\nabla_k \circo h_{ij} - |H|^{-1} \circo h_{ij} \nabla_k |H| ) \langle A^-_{ij}, \nabla_k \nu\rangle.
    \end{align}
On a neck the right-hand side can be estimated as follows.

\begin{lemma}\label{evolution |A-|^2 neck}
Let $\mathcal{M}_t$ be an $n$-dimensional mean curvature flow in $\mathbb{R}^{n+m}$. At any point of $\mathcal{M}_t$ which lies in an $(\e, 1)$-neck we have
    \begin{align*}
    (\partial_t - \Delta) |A^-|^2 &\leq \bigg(\frac{2}{n-1} + C\e\bigg)|h|^2 |A^-|^2- (2-C\e)|\nabla A^-|^2
    \end{align*}
where $C = C(n)$. 
\end{lemma}
\begin{proof}
We begin with the reaction terms in \eqref{evolution |A-|^2}. On an $(\e,1)$-neck we have $|A^-|^2 \leq C\e |h|^2$, and hence
    \[2 \sum_{i,j,k,l}\langle A^-_{ij} , A_{kl}^-\rangle^2 + 2 \sum_{i,j} |A^-_{ik} \otimes A^-_{kj} - A^-_{jk} \otimes A^-_{ki}|^2 \leq C \e |h|^2 |A^-|^2.\]
Moreover, in an orthonormal frame such that $h$ is diagonal,
    \begin{align*} 2 \sum_{i,j} |h_{ik} A^-_{kj} - h_{jk} A^-_{ki} |^2= 2 \sum_{i , j}(h_{ii}-h_{jj})^2 |A^-_{ij}|^2.
    \end{align*}
Using $(h_{ii} - h_{jj})^2 \leq (\frac{1}{n-1} + C\e)|h|^2$ we conclude that
    \begin{align*}
    &2 \sum_{i,j,k,l}\langle A^-_{ij} , A_{kl}^-\rangle^2+ 2 \sum_{i,j} |A^-_{ik} \otimes A^-_{kj} - A^-_{jk} \otimes A^-_{ki}|^2 + 2 \sum_{i,j} |h_{ik} A^-_{kj} - h_{jk} A^-_{kj} |^2\\
    & \qquad \leq \bigg(\frac{2}{n-1} + C\e\bigg)|h|^2 |A^-|^2.
    \end{align*}
  
Concerning the gradient terms in \eqref{evolution |A-|^2}, on an $(\e,1)$-neck we have $|\nabla h| \leq C\e |h||H|$ and hence
    \[4 \sum_{i,j,k} (\nabla_k \circo h_{ij} - |H|^{-1} \circo h_{ij} \nabla_k |H|) \langle A^-_{ij}, \nabla_k \nu\rangle \leq C\e |h||H||A^-| |\nabla \nu|.\]
Using the Codazzi equations as in \eqref{eq_principal torsion} we get $|H||\nabla \nu| \leq C|\nabla A^-|$. Inserting this above yields
    \[4 \sum_{i,j,k} (\nabla_k \circo h_{ij} - |H|^{-1} \circo h_{ij} \nabla_k |H| ) \langle A^-_{ij}, \nabla_k \nu\rangle \leq C\e |h||A^-| |\nabla A^-|.\]
The claim follows after an application of Young's inequality. 
\end{proof}

Next we consider the ratio $|A^-|^2/ Q^{1-\sigma}$ on an $(\e,1)$-neck, where $Q:= c |H|^2 - |A|^2$ is the quantity appearing in Naff's planarity estimate, defined as above. We have $Q>0$ on any $(\varepsilon,1)$-neck with $\varepsilon$ sufficiently small, and $(\partial_t - \Delta) Q \geq 2|h|^2 Q$ by Lemma~\ref{lem_pinching refinement}.

\begin{proposition}\label{evolution planarity ratio neck}
Let $\mathcal{M}_t$ be an $n$-dimensional mean curvature flow in $\mathbb{R}^{n+m}$. At any point of $\mathcal{M}_t$ which lies in an $(\varepsilon,1)$-neck we have
    \begin{align*}
    (\partial_t - \Delta) \frac{|A^-|^2}{Q^{1-\sigma}} &\leq -\bigg(2(1-\sigma) -\frac{2}{n-1}- C\e\bigg)|h|^2 \frac{|A^-|^2}{Q^{1-\sigma}}-(2-C\e) \frac{|\nabla A^-|^2}{Q^{1-\sigma}}
    \end{align*}
where $C = C(n)$
\end{proposition}
\begin{proof}
We first observe that 
    \begin{align*}
    (\partial_t - \Delta) Q^{1-\sigma} &= (1-\sigma)Q^{-\sigma} (\partial_t - \Delta) Q + \sigma(1-\sigma) Q^{-1-\sigma}|\nabla Q|^2
    \end{align*}
and hence 
    \begin{align*}
    (\partial_t - \Delta) \frac{|A^-|^2}{Q^{1-\sigma}} &\leq \frac{1}{Q^{1-\sigma}} (\partial_t - \Delta) |A^-|^2 - (1-\sigma)\frac{|A^-|^2}{Q^{2-\sigma}}(\partial_t - \Delta) Q  +\frac{2}{Q^{1-\sigma}}\bigg\langle\nabla Q^{1-\sigma}, \nabla \frac{|A^-|^2}{Q^{1-\sigma}}\bigg\rangle.
    \end{align*}
Using the fact that $|\nabla Q| \leq C \e |h|^3$ on an $(\e,1)$-neck, we deduce
    \[\frac{2}{Q^{1-\sigma}}\bigg\langle\nabla Q^{1-\sigma}, \nabla \frac{|A^-|^2}{Q^{1-\sigma}}\bigg\rangle \leq C\e \frac{|h||A^-||\nabla A^-|}{Q^{1-\sigma}} + C\e |h|^2 \frac{|A^-|^2}{Q^{1-\sigma}},\]
and then use Young's inequality to obtain
    \begin{align*}
    (\partial_t - \Delta) \frac{|A^-|^2}{Q^{1-\sigma}} &\leq \frac{1}{Q^{1-\sigma}} (\partial_t - \Delta) |A^-|^2 - (1-\sigma)\frac{|A^-|^2}{Q^{2-\sigma}}(\partial_t - \Delta) Q  +C\e \frac{|\nabla A^-|^2}{Q^{1-\sigma}} + C\e |h|^2 \frac{|A^-|^2}{Q^{1-\sigma}}.
    \end{align*}
To conclude the proof we estimate the first two terms on the right using Lemma~\ref{evolution |A-|^2 neck} and the inequality $(\partial_t - \Delta) Q \geq 2|h|^2 Q$ respectively. 
\end{proof}

We can now establish the planarity improvement theorem. 

\begin{proof}[Proof of Theorem~\ref{thm_planarity improvement}]
Up to a parabolic rescaling and translation in time, we may assume $|H(x_0,t_0)|^2 = \frac{n-1}{2}$ and $t_0 = -1$. This ensures that the neck has radius approximately equal to $\sqrt{2(n-1)}$ at time $t_0=-1$. Let us assume that the axis of the neck is the $x_{n+1}$-axis and translate in space so that $x_0$ lies in the hyperplane $\{x_{n+1} = 0\}$. 

We restrict attention to the portion of the shrinking neck which lies in the slab $\{|x_{n+1}| \leq L/2\}$ for times $t \in [-L^2,-1]$. If $\varepsilon$ is sufficiently small (depending on $L$) then this portion of the neck is properly embedded, has length approximately $L$, and has two approximately round boundary components---one in each of the two hyperplanes $\{x_{n+1} = \pm L/2\}$. By Proposition~\ref{evolution planarity ratio neck} we have
    \begin{align*}
    (\partial_t - \Delta) \frac{|A^-|^2}{Q^{1-\sigma}} &\leq -\bigg(2(1-\sigma) -\frac{2}{n-1}- C\varepsilon\bigg)|h|^2 \frac{|A^-|^2}{Q^{1-\sigma}}-(2-C\varepsilon) \frac{|\nabla A^-|^2}{Q^{1-\sigma}}
    \end{align*}
in $\{|x_{n+1}| \leq L/2\}$ for $t \in [-L^2, -1]$, where $C = C(n)$. We may assume $\varepsilon$ is small enough so that the final term on the right-hand side is nonpositive. Then, setting $u = |A^-|^2/ Q^{1-\sigma}$ and using the fact that on a shrinking neck $2(-t)|h|^2 \geq 1 - C\varepsilon $, we obtain 
    \begin{align*}
    \left(\partial_t - \Delta + \frac{a}{2(-t)}\right) u &\leq 0, \qquad a = a(\varepsilon) := 2(1-\sigma) -\frac{2}{n-1}- C\varepsilon.
    \end{align*}
If $\varepsilon$ is small then $a$ is close to $2(1-\sigma) - \frac{2}{n-1}$. Since $n \geq 5$ and $\sigma \in (0,\frac{1}{32}]$ we may assume $a \in (0, 2)$.

We construct an upper barrier for $u$ of the form $v(x,t) = 1-(1-\eta(t))\varphi(x)$, where 
    \[\varphi(x) = \cos\left(\pi x_{n+1}/L\right).\]
Note that $\varphi$ is positive in the region $\{|x_{n+1}| < L /2\}$ and vanishes on $\{x_{n+1} = \pm L/2\}$. We require $\eta(-L^2) = 1$ so that $v \equiv 1$ on the parabolic boundary $\{x_{n+1} = \pm L/2\} \cup \{t = -L^2\}$. Using the evolution equation $(\partial_t - \Delta) x_{n+1} = 0$ we obtain 
    \[(\partial_t - \Delta)\varphi = \frac{\pi^2}{L^2} \varphi |\nabla x_{n+1}|^2  \]
and hence 
    \begin{align*}
    \bigg(\partial_t - \Delta + \frac{a}{2(-t)}\bigg) v &= \bigg(\dot \eta - \frac{\pi^2}{L^2}(1-\eta) + \frac{a}{2(-t)}\eta\bigg) \varphi  + \frac{\pi^2}{L^2}(1-\eta)\varphi(1-|\nabla x_{n+1}|^2)+ \frac{a}{2(-t)}(1-\varphi).
    \end{align*}
The final term on the right is nonnegative, since $a > 0$ and $\varphi \leq 1$. We force the first term to vanish by choosing $\eta$ to be the solution to the IVP
    \[\dot \eta - \frac{\pi^2}{L^2}(1-\eta) + \frac{a}{2(-t)}\eta = 0, \qquad \eta(-L^2) = 1,\]
for times $t \in [-L^2, -1]$. It is easy to see that $\eta$ then satisfies $\eta \leq 1$, so using $|\nabla x_{n+1}|^2 = |e_{n+1}^\top|^2 \leq 1$ we obtain
    \[\bigg(\partial_t - \Delta + \frac{a}{2(-t)}\bigg) v \geq \frac{\pi^2}{L^2}(1-\eta)\varphi(1-|\nabla x_{n+1}|^2) \geq 0.\]
From this inequality we deduce
    \[(\partial_t - \Delta) \frac{u}{v} = \frac{1}{v}(\partial_t - \Delta) u - \frac{u}{v^2}(\partial_t - \Delta)v +  \frac{2}{v} \bigg\langle \nabla v, \nabla \, \frac{u}{v} \bigg\rangle\leq \frac{2}{v} \bigg\langle \nabla v, \nabla \, \frac{u}{v} \bigg\rangle\]
and hence, by the parabolic maximum principle,
    \[\frac{u(x_0,-1)}{v(x_0,-1)} \leq \max_{\{x_{n+1} = \pm L/2\} \cup \{t = -L^2\}} \frac{u}{v} = \max_{\{x_{n+1} = \pm L/2\} \cup \{t = -L^2\}} u \leq K.\]

To conclude we estimate $v(x_0,-1)$. From the definition of $\eta$ we see that 
    \[\frac{d}{dt}\left((-t)^{-a/2}\eta\right) \leq \frac{\pi^2}{L^2}(-t)^{-a/2}.\]
Integrating in time and using $\eta(-L^2) = 1$, this gives 
    \[\eta(-1) \leq L^{-a} + \pi^2\left(1-\frac{a}{2}\right)^{-1}(L^{-a} - L^2) \leq CL^{-a}.\]
We therefore have $v(x_0,-1) = \eta(-1) \leq CL^{-a}$ and hence $u(x_0, -1) \leq CL^{-a} K$. Recalling that $\varepsilon$ has been chosen so that $a = a(\varepsilon) \in (0,2)$, we choose $L$ large enough to ensure $CL^{-a} \leq \delta$. This implies $u(x_0,-1) \leq \delta K$, completing the proof. 
\end{proof}

To conclude this section we argue that if a sufficiently large spacetime region of a flow with surgery is a smoothly evolving neck, then the final timeslice of the neck satisfies \ref{s0}. This uses the planarity improvement theorem and the interior derivative estimates for $A^-$ proved in Appendix~\ref{sec_interior estimates A-}.

\begin{lemma}\label{lem_separation implies s0}
Let $\{\mathcal M_t\}_{t \in [0,t_0]}$ be a mean curvature flow with surgery in $\mathcal C_{n,m}(R,\alpha)$ such that $\max |A(\cdot,0)|^2 \leq R^{-2}$. Suppose \ref{s0} holds for all surgeries performed before time $t_0$. There then exist $\hat L = \hat L(n)$ and $\hat \varepsilon = \hat \varepsilon(n)$ with the following property: If $x_0 \in \mathcal M_{t_0}$ is such that $\hat{\mathcal P}(x_0,t_0, \hat L, \hat L^2)$ is an $(\hat \varepsilon, 1)$-shrinking neck then we have
    \[r_0^{2-2\sigma}|A^-|^2 + r_0^{4-2\sigma}|\nabla A^-|^2 + r^{6-2\sigma} |\nabla^2 A^-|^2 \leq \delta_0 C_{\pl}  R^{-2\sigma}\]
at every point of any normal neck of length $14\Lambda$ centred at $x_0$. Here $r_0$ is the mean radius of the $z = 0$ cross-section and $\delta_0 = \delta_0(n)$ is the constant appearing in \ref{s0}. 
\end{lemma}
\begin{proof}
If $\hat \varepsilon$ is sufficiently small then we can find a normal neck $\mathcal N \subset \mathcal{M}_{t_0}$ of length $14\Lambda$  centred at $x_0$. Choosing $\hat L$ large with respect to $\Lambda$, we ensure that for every $x \in \mathcal N$ and $t \geq t_0 - 10r_0^2$ the region $\hat{\mathcal P}(x,t, \hat L/2, \hat L^2/4)$ is an $(\hat \varepsilon, 1)$-shrinking neck. In particular, given any $\delta$, we can choose $\hat L$ and $\hat \varepsilon$ so that by the planarity improvement theorem the estimate
    \[|A^-|^2 \leq \delta C_{\pl} R^{-2\sigma} Q^{1-\sigma} \leq C \delta C_{\pl} R^{-2\sigma} r_0^{-2+2\sigma}\]
holds in $\mathcal P(x, t_0, 10 r_0, 10 r_0^2)$ for every $x \in \mathcal N$ and a constant $C = C(n)$. Keeping in mind that we are working in a neck region, Corollary~\ref{interior estimates A- H bdd below} provides bounds for the derivatives of $A^-$ at $(x,t_0)$ in terms of the supremum of $|A^-|^2$ over $\mathcal P(x, t_0, 10 r_0, 10 r_0^2)$. In particular,
    \[r_0^{2-2\sigma}|A^-|^2 + r_0^{4-2\sigma}|\nabla A^-|^2 + r_0^{6-2\sigma} |\nabla^2 A^-|^2 \leq C \delta C_{\pl} R^{-2\sigma}\]
holds at $(x,t_0)$ for every $x \in \mathcal N$, where $C$ is now larger but still only depends on $n$. We now choose $\delta$ so that $C\delta \leq \delta_0$. Since this choice can be made with $\delta$ depending only on $n$, our choices of $\hat L$ and $\hat \varepsilon$ depend only on $n$. 
\end{proof}


\section{Cylindrical estimates}\label{sec_cylindrical estimates}

In this section we prove cylindrical estimates for flows in the class $\mathcal C_{n,m}(R,\alpha)$. These estimates show that the principal part $h$ of the second fundamental form either becomes strictly positive or cylindrical whenever the curvature is large. We first consider the smooth case, before explaining how to carry the analyis over to the flow with surgery. 

\begin{theorem}[Cylindrical estimates for smooth flows]\label{thm_cylindrical via Stampacchia}
Let $\{\mathcal M_t\}_{t \in [0,T)}$ be a smooth mean curvature flow such that $\mathcal{M}_0 \in \mathcal C_{n,m}(R,\alpha)$ and $\max|A(\cdot, 0)|^2 \leq R^{-2}$. Given any $\varepsilon > 0$ we have 
    \[|A|^2 - \frac{1}{n-1}|H|^2 \leq \varepsilon |H|^2 + C\]
on $\mathcal{M}_t$ for all $t \in [0,T)$, where $C$ is a constant such that $CR^2$ depends only on $n, \alpha, \varepsilon$. 
\end{theorem}

Given a mean curvature flow as in Theorem~\ref{thm_cylindrical via Stampacchia} we define $Q = c|H|^2 - |A|^2$ with $c$ chosen as in Section~\ref{sec_planarity}. Let $\varepsilon > 0$ be such that $\frac{1}{n-1} + \varepsilon \leq c_n$ and set
    \[u:=|A|^2 - \left(\frac{1}{n-1}+\varepsilon\right)|H|^2, \qquad u_\beta := u/Q^{1-\beta},\]
where $\beta$ is assumed to be a constant in $(0,1)$. 

\begin{lemma}\label{lem_cylindrical quantity evol}
There exist positive constants $\gamma_0 = \gamma_0(n,\e)$ and $\gamma_1 = \gamma_1(n,\e)$ such that 
    \begin{align*}
    (\partial_t - \Delta) u_\beta &\leq -\gamma_0\frac{|\nabla A|^2}{Q^{1-\beta}} +2\beta|h|^2 u_\beta + \gamma_1 |A^-|^2 u_\beta - \beta (1-\beta) u_\beta \frac{|\nabla Q|^2}{Q^2} + 2(1-\beta)\bigg\langle \frac{\nabla Q}{Q}, \nabla u_\beta\bigg\rangle.
    \end{align*}
\end{lemma}
\begin{proof}
Appealing to Lemma~\ref{lem_pinching refinement}, with $c = \frac{1}{n-1} + \e$ and $a = 0$, we obtain 
    \begin{align*}
    (\partial_t - \Delta)u &\leq -2\left(|\nabla A|^2 - \left(\frac{1}{n-1}+\e\right)|\nabla H|^2\right) + 2|h|^2 u + 2\left(n\varepsilon + \frac{1}{n-1}\right)^{-1} |A^-|^2 u \\
    &\leq -2\left(1-\frac{n+2}{3}\left(\frac{1}{n-1} + \e\right)\right)|\nabla A|^2 + 2|h|^2 u + 2\left(n\varepsilon + \frac{1}{n-1}\right)^{-1} |A^-|^2 u
    \end{align*}
where in the second line we have used $\frac{n+2}{3}|\nabla A|^2 \geq |\nabla H|^2$ (see Lemma~\ref{lem_Kato}). Lemma~\ref{lem_pinching refinement} also gives $(\partial_t - \Delta)Q \geq 2 |h|^2Q$. Combining these two inequalities with
    \begin{align*}
    (\partial_t - \Delta) u_\beta &= \frac{(\partial_t - \Delta)u}{Q^{1-\beta}} - (1-\beta)u_\beta\frac{(\partial_t - \Delta)Q}{Q} - \beta (1-\beta) u_\beta \frac{|\nabla Q|^2}{Q^2} + 2(1-\beta)\bigg\langle \frac{\nabla Q}{Q}, \nabla u_\beta\bigg\rangle,
    \end{align*}
we arrive at the claim.
\end{proof}

Our cylindrical estimates are obtained via a Stampacchia iteration applied to the functions $u_\beta$. This technique has played an important role in the hypersurface setting e.g.\ in \cite{Hu84, HuSi99a, HuSi99b, HuSi09}. Each of those works exploits a crucial Poincar\'{e} inequality derived from Simons' identity. Finding analogous inequalities in higher codimension becomes rather subtle, because the zeroth-order terms in Simons’ identity contain components that cannot be controlled without additional geometric input. Building on an idea which we introduced in our earlier work \cite{Lynch--Nguyen_convexity}, we are able to use the planarity estimate to absorb these components. This leads to the following crucial $L^p$-estimates. 

\begin{theorem}\label{thm_cylindrical Lp}
Given any $\varepsilon >0$ there exist constants $p_0$ and $\theta_0$ which depend only on $n,\alpha,\varepsilon$ and have the following property. Let $\{\mathcal{M}_t\}_{t\in[0,T)}$ be a smooth mean curvature flow in $\mathcal C_{n,m}(R,\alpha)$ such that $\max|A(\cdot,0)|^2 \leq R^{-2}$. If $p \geq p_0$ and $\beta p^{1/2} \leq \theta_0$, then the function $u_{\beta +} := \max\{u_\beta, 0\}$ satisfies
    \[\sup_{t \in [0,T)}\,\int_{\mathcal {M}_t} u_{\beta+}^p \leq C\]
for a constant $C = C(n, R, \alpha, \varepsilon, \beta, p)$. 
\end{theorem}
\begin{proof}
Throughout the proof, we write $C$ whenever we need to represent a positive constant that can be bounded in terms of $n$, $\alpha$ and $\e$. Let us abbreviate $\varphi := u_{\beta +}$. We assume throughout that $p \geq 2$. By Lemma~\ref{lem_cylindrical quantity evol},
    \begin{align*}
    (\partial_t - \Delta) \varphi^p &= p \varphi^{p-1} (\partial_t - \Delta) u - p(p-1)\varphi^{p-2} |\nabla \varphi|^2\\
    &\leq -\gamma_0 p \varphi^{p-1}\frac{|\nabla A|^2}{Q^{1-\beta}} + 2\beta p|h|^2 \varphi^p + \gamma_1 p |A^-|^2 \varphi^p\\
    &\qquad +2(1-\beta)p \varphi^{p-1}\bigg\langle\frac{\nabla Q}{Q}, \nabla \varphi\bigg\rangle - p(p-1)\varphi^{p-2}|\nabla \varphi|^2.
    \end{align*}
Using Young's inequality together with $u \leq |A|^2$ we obtain
    \begin{align*}
    -\gamma_0 p\varphi^{p-1}\frac{|\nabla A|^2}{Q^{1-\beta}} +2(1-\beta)p \varphi^{p-1}\bigg\langle\frac{\nabla Q}{Q}, \nabla \varphi\bigg\rangle &\leq -(\gamma_0 p - C) \varphi^p\frac{|\nabla A|^2}{|A|^2} + \frac{p(p-1)}{2} \varphi^{p-2} |\nabla \varphi|^2,
    \end{align*}
so we have 
    \begin{align*}
    (\partial_t - \Delta) \varphi^p &\leq -(\gamma_0 p - C) \varphi^p\frac{|\nabla A|^2}{|A|^2} + 2\beta p|h|^2 \varphi^p + \gamma_1 p |A^-|^2 \varphi^p - \frac{p(p-1)}{2}\varphi^{p-2}|\nabla \varphi|^2.
    \end{align*}
It follows that 
    \begin{align}\label{eq_cylindrical Lp evolution}
    \begin{split}
    \frac{d}{dt} \int_{\mathcal{M}_t} \varphi^p &\leq -(\gamma_0 p - C) \int_{\mathcal{M}_t} \varphi^p\frac{|\nabla A|^2}{|A|^2} + 2\beta p \int_{\mathcal{M}_t}|h|^2 \varphi^p + \gamma_1 p \int_{\mathcal{M}_t} |A^-|^2 \varphi^p\\
    &\qquad - \frac{p(p-1)}{2} \int_{\mathcal{M}_t} \varphi^{p-2}|\nabla \varphi|^2 - \int_{\mathcal{M}_t} |H|^2 \varphi^p.
    \end{split}
    \end{align}

The diffusion and gradient of curvature terms can be used to absorb the unfavourable reaction term $2\beta p\int_{\mathcal{M}_t}|h|^2 \varphi^p$ via the following Poincar\'{e} inequality.

\begin{claim}
If $f : \mathcal{M}_t \to [0,\infty)$ is such that $|A|^2 \geq \left(\frac{1}{n-1} + \e\right)|H|^2$ holds on $\{f > 0\}$ then we have
    \[\int_{\mathcal M_t} |h|^2 f^2 \leq C\int_{\mathcal M_t} f^2 \frac{|\nabla A|^2}{|A|^2} + C\int_{\mathcal M_t} f |\nabla f| \frac{|\nabla A|}{|A|} + C\int_{\mathcal M_t} |A||A^-| f^2\]
for a constant $C = C(n,\alpha,\varepsilon)$. 
\end{claim}
\begin{proof} Let us define
    \[E_{klij} := \nabla_k \nabla_l A_{ij} + \nabla_l\nabla_kA_{ij} - \nabla_i\nabla_j A_{kl} - \nabla_j \nabla_i A_{kl}.\]
According to Simons' identity the tensor $E$ is given by some cubic expression $A \ast A \ast A$. By splitting each factor of $A$ in this expression into components parallel and orthogonal to $H$, in \cite[Lemma~3.1]{Lynch--Nguyen_convexity} we showed that for an arbitrary submanifold with $|H| > 0$ we have the pointwise inequality 
    \begin{align*}
    |E|^2 &\geq 8|h|^2\tr(h^4) - 8\tr(h^3)^2 - K|A|^5|A^-|=4\sum_{i,j}\lambda_i^2\lambda_j^2(\lambda_i-\lambda_j)^2 - K|A|^5|A^-|,
    \end{align*}
where $K$ is a constant depending only on $n$ and we use $\lambda_i$ to denote the eigenvalues of $h$. The term 
    $\sum_{i,j}\lambda_i^2\lambda_j^2(\lambda_i-\lambda_j)^2$
vanishes precisely when $h$ coincides (up to scaling) with the second fundamental form of a cylinder $\mathbb{R}^k \times \mathbb{S}^{n-k}$, in which case $|h|^2 - \frac{1}{n-k} |H|^2 = 0$. Based on this observation a simple proof by contradiction shows that there exists a positive $\delta = \delta(n,\alpha_1,\varepsilon)$ with the following property: If the inequality 
    \begin{equation}\label{eq:pinching for poincare}
    \left(\frac{1}{n-2} - \alpha_1\right)|H|^2 \geq |A|^2 \geq \left(\frac{1}{n-1} + \varepsilon\right)|H|^2
    \end{equation}
holds at some point, then at that same point we have
    \[4\sum_{i,j}\lambda_i^2\lambda_j^2(\lambda_i-\lambda_j)^2  + |A|^5 |A^-| \geq \delta |A|^6\]
and hence 
    \begin{equation}\label{eq:E bound for poincare}
    |E|^2 \geq \delta |A|^6 - (K+1)|A|^5|A^-|.
    \end{equation}
    
Our discussion so far applies to a general submanifold. Returning to our solution $\mathcal{M}_t$, the first inequality in \eqref{eq:pinching for poincare} holds everywhere by hypothesis, so if the second inequality holds on $\{f > 0\}$ then from \eqref{eq:E bound for poincare} we obtain 
    \begin{align*}
    \delta \int_{\mathcal M_t} |A|^2 f^2 &\leq \int_{\mathcal M_t} \frac{|E|^2}{|A|^4} f^2 + (K+1)\int_{\mathcal M_t} |A||A^-|f^2\\
    & = \int_{\mathcal M_t} \bigg(\frac{A\ast A \ast A}{|A|^4} \ast \nabla^2 A\bigg)f^2 + (K+1)\int_{\mathcal M_t} |A||A^-|f^2.
    \end{align*}
The claim is deduced from this inequality using integration by parts and Young's inequality exactly as in \cite[Proposition~3.2]{Lynch--Nguyen_convexity}. 
\end{proof}

Setting $f = \varphi^{p/2}$ in the Poincar\'{e} inequality yields 
    \[\int_{\mathcal M_t} |h|^2 \varphi^p \leq C\int_{\mathcal M_t} \varphi^p \frac{|\nabla A|^2}{|A|^2} + C p \int_{\mathcal M_t} \varphi^{p-1} |\nabla \varphi| \frac{|\nabla A|}{|A|} + C\int_{\mathcal M_t} |A||A^-| \varphi^p.\]
We multiply this inequality by $3\beta p$ and then add $\gamma_1 p \int_{\mathcal M_t} |A^-|^2 \varphi^p$ to both sides to get
    \begin{align*}
    &3\beta p \int_{\mathcal M_t} |h|^2 \varphi^p + \gamma_1 p\int_{\mathcal M_t} |A^-|^2 \varphi^p \\
    &\qquad \leq C\beta p \int_{\mathcal M_t} \varphi^p \frac{|\nabla A|^2}{|A|^2} + C \beta p^2 \int_{\mathcal M_t} \varphi^{p-1} |\nabla \varphi| \frac{|\nabla A|}{|A|} + C p\int_{\mathcal M_t} |A||A^-| \varphi^p.
    \end{align*}
The last term on the right-hand side can be absorbed using the planarity estimate. Indeed, $|A^-|^2 \leq C_{\pl} R^{-2\sigma} Q^{1-\sigma}$ implies $|A^-| \leq C R^{-\sigma} |h|^{1-\sigma}$ and hence by Young's inequality
    \[C p\int_{\mathcal M_t} |A||A^-| \varphi^p \leq C R^{-\sigma} p\int_{\mathcal M_t} |h|^{2-\sigma} \varphi^p \leq C R^{-2} \beta^{1-2/\sigma}  p \int_{\mathcal M_t} \varphi^p  + \beta p \int_{\mathcal M_t} |h|^2 \varphi^p,\]
so we have
    \begin{align*}
    &2\beta p \int_{\mathcal M_t} |h|^2 \varphi^p + \gamma_1 p\int_{\mathcal M_t} |A^-|^2 \varphi^p\\
    &\qquad \leq C\beta p \int_{\mathcal M_t} \varphi^p \frac{|\nabla A|^2}{|A|^2} + C \beta p^2 \int_{\mathcal M_t} \varphi^{p-1} |\nabla \varphi| \frac{|\nabla A|}{|A|} + C R^{-2} \beta^{1-2/\sigma} p \int_{\mathcal M_t} \varphi^p.
    \end{align*}

We combine this inequality with \eqref{eq_cylindrical Lp evolution}, after estimating 
    \[\beta p^2 \int_{\mathcal M_t} \varphi^{p-1} |\nabla \varphi| \frac{|\nabla A|}{|A|} \leq \frac{\beta p^{3/2}}{2} \int_{\mathcal M_t} \varphi^p \frac{|\nabla A|^2}{|A|^2} + \frac{\beta p^{5/2}}{2} \int_{\mathcal M_t} \varphi^{p-2} |\nabla \varphi|^2,\]
in order to obtain
    \begin{align*}
    \frac{d}{dt} \int_{\mathcal M_t} \varphi^p &\leq -(\gamma_0 p - C - C\beta p^{3/2}) \int_{\mathcal M_t} \varphi^p\frac{|\nabla A|^2}{|A|^2} - (p(p-1)/2 - C\beta p^{5/2}) \int_{\mathcal M_t} \varphi^{p-2}|\nabla \varphi|^2\\
    &\qquad + C R^{-2} \beta^{1-2/\sigma} p \int_{\mathcal M_t} \varphi^p - \int_{\mathcal{M}_t} |H|^2 \varphi^p.
    \end{align*}
If $\beta p^{1/2}$ is sufficiently small and $p$ is sufficiently large then the first two terms on the right are nonpositive, in which case 
    \begin{align*}
    \begin{split}
    \frac{d}{dt} \int_{\mathcal M_t} \varphi^p &\leq C R^{-2} \beta^{1-2/\sigma} p \int_{\mathcal M_t} \varphi^p,
    \end{split}
    \end{align*}
and hence by ODE comparison
    \[\int_{\mathcal M_t} \varphi^p \leq \exp(C R^{-2} \beta^{1-2/\sigma} p T) \int_{\mathcal M_0} \varphi^p.\]
The maximal time $T$ can be bounded from above in terms of $n$, $R$ and $\alpha$ by a standard maximum principle argument using $(\partial_t - \Delta) Q \geq 2|h|^2 Q \geq \frac{2}{nc}Q^3$, whereas at time $t = 0$ we have $\varphi \leq CR^{-2\beta}$. The claim follows.
\end{proof}

We are now prepared to prove the cylindrical estimates for smooth flows.

\begin{proof}[Proof of Theorem~\ref{thm_cylindrical via Stampacchia}]
We establish an upper bound for $\max_{\mathcal M_t} u_\beta$ which is uniform in time using Stampacchia iteration. The claim follows from this estimate and Young's inequality. Let $u_{\beta, k} := \max\{u-k,0\}$ and set
    \[S(k, t) := \{x \in \mathcal M_t : u_{\beta, k}(x,t) > 0\}, \qquad |S(k)|:=\int_0^T \vol_{\mathcal M_t}(S(k,t))\,dt.\]
Let $p \geq 2$ be a constant which we will choose to be large in the course of the proof. The same arguments leading to \eqref{eq_cylindrical Lp evolution} show that  
    \begin{align*}
    \begin{split}
    \frac{d}{dt} \int_{\mathcal M_t} u_{\beta, k}^p &\leq -(\gamma_0 p - C) \int_{\mathcal M_t} u_{\beta, k}^{p-1}u_{\beta}\frac{|\nabla A|^2}{|A|^2} + 2\beta p \int_{\mathcal M_t}|h|^2 u_{\beta,k}^{p-1} u_\beta + \gamma_1 p \int_{\mathcal M_t} |A^-|^2 u_{\beta,k}^{p-1}u_\beta\\
    &\qquad - \frac{p(p-1)}{2} \int_{\mathcal M_t} u_{\beta, k}^{p-2}|\nabla u_{\beta}|^2 - \int_{\mathcal M_t} |H|^2 u_{\beta, k}^p.
    \end{split}
    \end{align*}
Assuming $p$ is large enough so that the first term on the right is nonpositive, we obtain
    \begin{align*}
    \begin{split}
    \frac{d}{dt} \int_{\mathcal M_t} u_{\beta, k}^p + \frac{p(p-1)}{2} \int_{\mathcal M_t} u_{\beta, k}^{p-2}|\nabla u_{\beta}|^2 &\leq C p \int_{S(k,t)}|h|^2 u_{\beta}^p.
    \end{split}
    \end{align*}
By integrating this inequality in time we find that for every $t_0 \in (0,T)$ we have
    \begin{align}\label{eq_cylindrical integral ineq}
    \int_{\mathcal M_{t_0}} u_{\beta, k}^p + \frac{p(p-1)}{2} \int_0^{t_0}\int_{\mathcal M_t} u_{\beta, k}^{p-2}|\nabla u_{\beta}|^2 &\leq C p \int_0^{t_0}\int_{S(k,t)}|h|^2 u_{\beta}^p.
    \end{align}
    
Using \eqref{eq_cylindrical integral ineq} together with the Michael--Simon--Sobolev inequality \cite{Michael--Simon} and the $L^p$-estimates from Theorem~\ref{thm_cylindrical Lp}, one deduces that 
    \begin{equation}\label{eq_cylindrical iteration}
    |S(\ell)| \leq \frac{C}{(\ell-k)^p} |S(k)|^{\gamma}
    \end{equation}
for every $\ell > k \geq k_0$, where $C$ and $k_0$ both depend only on $n, R, \alpha, \varepsilon, \beta$ and the constant $\gamma > 1$ depends only on $n$. (An important technical point is that, to obtain \eqref{eq_cylindrical iteration}, we need $p$ and $\beta$ to lie in more restrictive ranges than before, say $p \geq p_1(n,\alpha, \e)$ and $\beta p^{1/2} \leq \theta_1(n, \alpha, \e)$ where $p_1 > p_0$ and $\theta_1 < \theta_0$.) The desired global upper bound for $\max_{\mathcal{M}_t} u_\beta$ follows immediately from \eqref{eq_cylindrical iteration} and Stampacchia's lemma (see e.g. \cite[Lemma~B.1]{KS}). The arguments needed to obtain \eqref{eq_cylindrical iteration} from \eqref{eq_cylindrical integral ineq} are by now standard and have appeared many times in the literature, in e.g.\ \cite{Hu84,HuSi09,Lynch--Nguyen_convexity}. We refer to those texts for further details. 
\end{proof}

We now explain how the proof of Theorem~\ref{thm_cylindrical via Stampacchia} needs to be modified for a flow with surgery. A crucial point is that the parameters in the standard surgery procedure can be chosen so that the key quantity $u_\beta$ does not increase; this follows from Theorem~\ref{geom quants N hat}. In addition, we need to assume \ref{s0} so that the planarity estimate holds by Theorem~\ref{thm_planarity surgeries}. 

\begin{theorem}\label{thm_cylindrical via Stampacchia surgery}
Let $\{\mathcal M_t\}_{t \in [0,T]}$ be a mean curvature flow with surgery such that $\mathcal{M}_0 \in \mathcal C_{n,m}(R,\alpha)$ and $\max|A(\cdot, 0)|^2 \leq R^{-2}$. Suppose \ref{s0} holds for every surgery. Given any $\varepsilon > 0$ we have 
    \[|A|^2 - \frac{1}{n-1}|H|^2 \leq \varepsilon |H|^2 + C\]
on $\mathcal{M}_t$ for all $t \in [0,T]$, where $C$ is a constant such that $CR^2$ depends only on $n, \alpha, \varepsilon$. 
\end{theorem}
\begin{proof}
Let us define $u$ and $u_\beta$ exactly as before. We first demonstrate that Theorem~\ref{thm_cylindrical Lp} holds for the flow with surgery. As a consequence of \ref{s0} and Theorem~\ref{thm_planarity surgeries}, the planarity estimate $|A^-|^2 \leq C_{\pl}Q^{1-\sigma}$ holds on $\mathcal M_t$ at every time $t \in [0,T]$. Therefore, we may proceed exactly as in the proof of Theorem~\ref{thm_cylindrical Lp} to deduce that in between surgery times, i.e.\ when the flow is smooth, we have
    \begin{equation}\label{eq_cylindrical surgeries DE}
    \frac{d}{dt}\bigg(e^{-Ct} \int_{\mathcal M_t} u_{\beta+}^p\bigg) \leq 0, \qquad C := \tilde C(n,\alpha,\varepsilon) R^{-2} \beta^{1-2/\sigma} p,
    \end{equation}
provided $p \geq p_0$ and $\beta p^{1/2} \leq \delta_0$. We claim that $\int_{\mathcal M_t} u_{\beta+}^p$ does not increase when we perform surgery. If this is true then \eqref{eq_cylindrical surgeries DE} can be integrated to obtain the desired $L^p$-estimate for $u_{\beta+}$. Recall that the standard surgery procedure removes a portion of a neck and replaces it with two strictly convex axially symmetric hypersurface caps. The caps are constructed so that the function $u$ is negative on them, so $u_{\beta+}$ vanishes identically in these regions, and hence we can ignore them. On the two adjacent collar regions Theorem~\ref{geom quants N hat} ensures that $u$ is nonincreasing under surgery. Moreover, $Q$ is nondecreasing and $\sqrt{\det g}$ is nonincreasing. Therefore, the integral of $u_{\beta+}^p$ over these two collars does not increase and, consequently, performing any number of standard surgeries on necks in $\mathcal M_t$ does not increase $\int_{\mathcal M_t} u_{\beta+}^p$. After performing surgery some components are discarded, and this clearly only decreases $\int_{M_t} u_{\beta +}^p$. Therefore, $\int_{M_t} u_{\beta +}^p$ is indeed nonincreasing at each surgery time, and hence from \eqref{eq_cylindrical surgeries DE} we deduce an estimate of the form 
    \[\int_{\mathcal M_t} u_{\beta+}^p \leq C(n,R,\alpha, \varepsilon, \beta, p,T).\]
Since $Q$ does not decrease under surgery we can bound $T$ from above in terms of $n,\alpha,R$ exactly as in the smooth case. 

With the $L^p$-estimates at hand, for suitably small $\beta$ we establish a global upper bound for $u_\beta$ using Stampacchia iteration as in the smooth case. Indeed, proceeding as in Theorem~\ref{thm_cylindrical via Stampacchia} we find that whenever the flow is smooth we have
    \begin{align*}
    \frac{d}{dt} \int_{\mathcal M_t} u_{\beta, k}^p + \frac{p(p-1)}{2} \int_{\mathcal M_t} u_{\beta, k}^{p-2}|\nabla u_{\beta}|^2 &\leq C p \int_{S(k,t)}|h|^2 u_{\beta}^p.
    \end{align*}
As above, $\int_{\mathcal M_t} u_{\beta,k}^p$ is nonincreasing across surgeries, so we can inegrate this inequality in time to get \eqref{eq_cylindrical integral ineq} for every $t_0 \in [0,T]$. The rest of the proof is the same as in the smooth case. 
\end{proof}


\section{Neck detection}\label{sec_ND}

This section follows Section~7 of \cite{HuSi09}, except for a couple of results that are not needed in the hypersurface setting. These are the hypersurface detection lemma (Lemma~\ref{lem_HDL}) and a bound for the spacetime separation between surgery necks (Lemma~\ref{lem_ss}).

The following lemma is proven by integrating the gradient estimates of Corollary~\ref{cor_gradient} along curves as in \cite[Lemma~7.2]{HuSi09}. We refer to the constant $c^{\#}$ from Lemma~\ref{cor_gradient}.

\begin{lemma}\label{lem_Harnack2}
Set $d^{\#}= (8(n-1)^2 c^{\#})^{-1}$. Given $r, \theta \in (0, d^\#]$ we have $|H(q,s)|/|H(p,t)| \in [\frac{1}{2}, 2]$ for all $(q,s) \in \hat{\mathcal P }(p,t,r,\theta)$, provided this region does not contain surgeries. 
\end{lemma}

Next we state our neck detection lemma (cf.\ Lemma~7.4 in \cite{HuSi09}). 

\begin{lemma}[Neck detection]\label{lem_NDL}
Let $\{\mathcal M_t\}_{t \in [0,t_0]}$ be a mean curvature flow with surgery in $\mathcal C_{n,m}(R,\alpha)$ such that $\max|A(\cdot,0)|^2 \leq R^{-2}$. Suppose \ref{s0} holds at times prior to $t_0$. Given any $\varepsilon, \theta, L > 0$ and $k \geq k_0$ there exist constants $\eta_0 > 0$ and $H_0$ such that if $(p_0,t_0)$ satisfies
    \begin{equation}\label{eq_ND1}
    |H(p_0,t_0)|\geq H_0, \qquad |A(p_0,t_0)|^2 - \tfrac{1}{n-1}|H(p_0,t_0)|^2 \geq -\eta_0|H(p_0,t_0)|^2 \tag{ND1}
    \end{equation}
and
    \begin{equation}\label{eq_ND2}
    \hat{\mathcal{P}}(p_0,t_0,L,\theta) \text{ does not contain surgeries,} \tag{ND2}
    \end{equation}
then $\hat{\mathcal P}(p_0, t_0, L, \theta) $ is an $(\e,k_0-1)$-shrinking neck and $\hat{\mathcal P}(p_0,t_0,L-1,\frac{\theta}{2})$ is an $(\e,k)$-shrinking neck. The constants $\eta_0$ and $H_0 R$ depend only on $n, m, \alpha, \varepsilon, k, \theta, L$. 
\end{lemma}

\begin{proof}
We use a proof by contradiction similar to the proof of \cite[Lemma~7.4]{HuSi09}. We provide a sketch and refer to \cite{HuSi09} for further details. 

Suppose there exist $\varepsilon, \theta, L$ such that the first claim is not true. There then exists a sequence of flows $\mathcal M_t^j$ satisfying the hypotheses of the lemma and a sequence $p_j\in\mathcal{M}_{t_j}^j$ such that $|H(p_j, t_j)|\rightarrow\infty$, $|A(p_j ,t_j)|^2/|H(p_j,t_j)|^2 \rightarrow \frac{1}{n-1}$ and $\hat{\mathcal{P}}^j(p_j,t_j,L,\theta)$ does not contain surgeries, but $(p_j,t_j)$ does not lie at the centre of a $\left(\varepsilon,k_0-1,L,\theta\right)$-shrinking neck.

After shifting $(p_j,t_j) \to (0,0)$ and parabolically rescaling by $(n-1)/|H(p_j,t_j)|$ we can use the gradient estimates of Corollary~\ref{thm_highergradest} to extract a limit in $C^{k_0+1}$ of the regions $\mathcal P^j(0,0,d,d)$, for some $d > 0$. The cylindrical estimates imply that on the limiting mean curvature flow we have $|A|^2 - \frac{1}{n-1}|H|^2 \leq 0$ with equality at the spacetime point $(0,0)$. Then because of \eqref{eqn_pinchpres2} the strong maximum principle implies $|A|^2-\frac{1}{n-1}|H|^2 =0$, $|A^-|^2 = 0$ and $|\nabla A|^2 = 0$ everywhere on the limiting flow, which therefore coincides with part of a standard shrinking cylinder solution (see e.g.\ \cite{Huisken1993}). Since the curvature of the cylinder is constant at each time we can iterate this argument to extend the convergence to $\mathcal P^j(0,0,L,\theta)$, and thus conclude that $(0,0)$ lies at the centre of an $(\e,k_0-1,L,\theta)$-shrinking neck for all sufficiently large $j$. This is a contradiction, so the claim is proven. 

The proof of the second claim is analogous, except that in addition to the pointwise gradient estimates we use standard interior parabolic estimates (see e.g.\ Appendix~\ref{sec_interior estimates A-}) to get uniform bounds in $C^{k+2}$ in the appropriate spacetime regions
\end{proof}

The neck detection lemma establishes the existence of neck regions of large curvature. By definition, a neck is close to being a hypersurface in some $(n+1)$-dimensional affine subspace of $\mathbb{R}^{n+m}$. In the proof of our neck continuation theorem (Theorem~\ref{thm_neckcontinuation}) we will require this `almost-hypersurface' property to hold in \emph{all} regions where the curvature is large, not just on necks.

\begin{lemma}[Hypersurface detection]\label{lem_HDL}
Let $\{\mathcal M_t\}_{t \in [0, t_0]}$ be a mean curvature flow with surgery in $\mathcal C_{n,m}(R,\alpha)$ such that $\max|A(\cdot,0)|^2 \leq R^{-2}$. Suppose \ref{s0} holds at times prior to $t_0$. Let $\varepsilon, r, \theta> 0$ and $k \geq 0$ be given such that $r, \theta \in (0,d^\#]$, where $d^\#$ is the constant from Lemma~\ref{lem_Harnack2}. There exists a constant $H_0$ such that if $|H(p_0, t_0)|\geq H_0$ and the neighbourhood $\hat{\mathcal P}(p_0,t_0,r,\theta)$ does not contain surgeries, then the estimate $|\nabla^\ell A^-| + |\nabla^{\ell+1} \nu| \leq \e |H|^{\ell+1}$ holds in $\hat{\mathcal P} (p_0,t_0,\frac{r}{2},\frac{\theta}{2})$ for every $0 \leq \ell \leq k$. Moreover $H_0 R$ depends only on $n, m, \alpha, \varepsilon, k, r,\theta$. 
\end{lemma}
\begin{proof}
We argue by contradiction using the planarity estimate. Suppose that for some $\varepsilon, r, \theta$ and $k \geq 0$ we have a sequence of flows where the claim fails in $\hat{\mathcal P}^j(p_j,t_j,r,\theta)$ but $|H(p_j,t_j)| \to \infty$. After rescaling to arrange $|H(p_j,t_j)| = n-1$, using the gradient estimates and interior parabolic theory we can extract a smooth limit of the slightly smaller regions $\hat{\mathcal P} (p_0,t_0,\frac{r}{2},\frac{\theta}{2})$. By the planarity estimate, $|A^-|$ vanishes identically on the limiting flow, which is therefore a hypersurface by \cite[Proposition~2.5]{Naff_planarity}. At every point of the limit we have $|\nabla^\ell A^-| + |\nabla^{\ell+1} \nu| = 0$ for all $\ell \geq 0$, which yields a contradiction for large $j$. 
\end{proof}

The following lemma shows that \eqref{eq_ND2} is a consequence of the other assumptions of the neck detection lemma if the curvature at $(p_0,t_0)$ is larger by a fixed factor than the curvature in regions previously changed by surgeries. For the proof see \cite[Lemma 7.10]{HuSi09}. 

\begin{lemma}\label{lem_NDL2replace}
Consider a flow with surgery as in Lemma~\ref{lem_NDL}. Let $\varepsilon, k, \theta, L$ be given, where $\theta \leq d^{\#}$ and $k \geq k_0$. There exist $\eta_0 > 0$ and $H_0$ such that if $(p_0,t_0)$ satisfies
    \begin{align*}
    |H(p_0,t_0)|\geq \max\{H_0,5K\}, \qquad |A(p_0,t_0)|^2 - \tfrac{1}{n-1}|H(p_0,t_0)|^2 \geq -\eta_0|H(p_0,t_0)|^2,
    \end{align*}
where $K$ is the maximum of the curvature at points changed by surgeries at times before $t_0$, then $(p_0,t_0)$ satisfies \eqref{eq_ND2} and hence the conclusions of Lemma~\ref{lem_NDL} hold true. In addition $\mathcal P(p_0,t_0, \frac{n-1}{|H(p_0,t_0) |} L, \frac{(n-1)^2}{K^2} d^{\#})$ does not contain surgeries.
\end{lemma}

In the next result we assume our flow with surgery satisfies the following properties (as in \cite{HuSi09}). These will be consequences of the surgery algorithm defined in the next section.

\begin{enumerate}[label=(s\arabic*)]
\item\label{s1}
Each surgery is performed at a cross-section $\Sigma_{z_0}$ of a normal $(\varepsilon_*, k_0)$-cylindrical submanifold neck with $r(z_0) = r_* = (n-1)/K_*$, where $K_*$ is some fixed value (the same for every surgery). Moreover, on each side of $\Sigma_{z_0}$ there is a portion of the neck of length $L_*$, where $L_* \geq 5+10\Lambda$. 
\item\label{s2} After each surgery one of the two remaining portions of the neck belongs to a component which is discarded, while in the other remaining portion the part of the neck unchanged by surgery has the following properties: The first cross-section (which borders the region changed by surgery) has mean radius $r(z) \leq 11r^*/10$, the final cross-section is such that $r(z) \geq 2 r^*$, and for every cross-section between the two we have $r^* \leq r(z) \leq 2 r^*$.
\item\label{s3} If we consider any of the surgeries performed at any  surgery time $t$, there exists a component that is discarded afterwards which contains a point $p$ with curvature $|H(p,t)|\geq 10 K_*$ and which would not have been removed if not for that surgery.  
\end{enumerate}

For a proof of the following see \cite[Lemma~7.12]{HuSi09}. 

\begin{lemma}\label{prop_discborder}
Consider a mean curvature flow with with surgery as in Lemma~\ref{lem_NDL}. Suppose \emph{\ref{s0}-\ref{s3}} hold at times prior to $t_0$. Let $\theta, L >0$ be such that $\theta \leq d^\#$ and $L\geq 20$. Given $\varepsilon > 0$ there exist constants $\eta_0 > 0$ and $H_0$ with the following property. Suppose
    \begin{enumerate}[label=\emph{(\roman*)}]
    \item $(p_0, t_0)$ satisfies the hypotheses \eqref{eq_ND1} and \eqref{eq_ND2} of Lemma~\ref{lem_NDL}, and
    \item the parabolic neighbourhood $\hat{\mathcal P}(p_0,t_0,L,\theta)$ is adjacent to a surgery region, meaning it does not contain surgeries but there exists a point $p$ at distance $d_{g(t_0)}(p,p_0)=\frac{n-1}{|H(p_0,t_0)|}L$ which belongs to the boundary of a region changed by surgery at a time in the interval $[t_0-\frac{(n-1)^2}{|H(p_0,t_0)|^2} \theta, t_0]$.
    \end{enumerate}
Then $(p_0,t_0)$ lies at the centre of a normal $(\varepsilon, k_0-1)$-cylindrical submanifold neck $\mathcal N$ of length $L-3$ and one of the two components of $\partial \mathcal N$ is also the boundary of a closed domain $\mathcal D$ which is diffeomorphic to a standard $n$-ball and has no interior points in common with $\mathcal N$. In addition, $|H(q,t_0)| < 5K_*$ for all $q\in \mathcal N\cup \mathcal D$, where $K_*$ is the constant from \emph{\ref{s1}}.
\end{lemma}

We will construct our surgery algorithm so that the necks where we perform surgery are separated in spacetime. This property is needed to justify \ref{s0} and consequently ensure that the planarity estimate is preserved by surgery (as explained in Section~\ref{sec_planarity}).

\begin{lemma}\label{lem_ss}
Let $\{\mathcal M_t\}_{t\in [0,t_0]}$ be a mean curvature flow with surgery in $\mathcal C_{n,m}(R,\alpha)$ such that $\max|A(\cdot,0)|^2 \leq R^{-2}$ and \ref{s0} holds for all surgeries before time $t_0$. Given any $\hat \varepsilon > 0$, $\hat L \geq 5 +14\Lambda$ (where $\Lambda$ is the length parameter from the standard surgery procedure) there exist constants $\varepsilon_s, K_s$ with the following property. Suppose $(p_0,t_0)$ lies at the centre of a normal $(\varepsilon_s, 1)$-cylindrical submanifold neck of length $10\hat L$, on a cross-section with mean radius $r_*$. In addition, suppose all surgeries before time $t_0$ satisfy \ref{s1} with $\varepsilon_* \leq \varepsilon_s$, $K_* \geq K_s$. It then follows that $\hat{\mathcal P}(p_0, t_0, \hat L, \hat L^2)$ is an $(\hat \varepsilon,1)$-shrinking neck; in particular, this region does not contain surgeries. The constants can be chosen such that $\varepsilon_s$ and $K_sR$ both depend only on $n, m, \alpha, \hat L, \hat \varepsilon$.
\end{lemma}
\begin{proof}
Let $\mathcal N_{0}$ denote the neck in $\mathcal M_{t_0}$ containing $p_0$. We write $r_0 = (n-1)/|H(p_0,t_0)|$. We may assume $\varepsilon_s$ is small enough so that $r_0 \in [r_*/2, 2r_*]$. 

\begin{claim}\label{claim_small surgery free interval}
If $\varepsilon_s$ is sufficiently small, depending only on $n$ and $\hat L$, then we can find $\kappa = \kappa(n,\alpha)$ such that $\hat{\mathcal P}(p_0,t_0, \hat L, \kappa)$ does not contain surgeries. 
\end{claim}
\begin{proof}[Proof of Claim~\ref{claim_small surgery free interval}]
The claim is that $\mathcal B_{g(t_0)}(p_0, \hat L r_0)$ is unchanged by surgery in the interval $(t_0 - \kappa r_0^2, t_0]$. We establish this by proving that the larger ball $\mathcal B_{g(t_0)}(p_0, 2\hat Lr_*)$ is unchanged by surgery in the interval $(t_0 - \kappa r_*^2, t_0]$. We argue that if this region does contain a surgery then there must be a positively curved cap not far from $p_0$ at some time slightly before $t_0$, which then persists until time $t_0$, contradicting the fact that $\mathcal N_{0}$ is a neck.

Suppose there is a time $t_1 \in (t_0 - \kappa r_*^2, t_0)$ such that $\mathcal B_{g(t_0)}(p_0, 2\hat Lr_*)$ is unchanged by surgery in the interval $(t_1, t_0]$, but there exists a point $p_1 \in \mathcal B_{g(t_0)}(p_0, 2\hat Lr_*)$ which is affected by a surgery performed at time $t_1$. We can find a point $p_2$ at a controlled distance from $p_1$ such that $|A(p_2, t_1+)|^2 - \frac{1}{n}|H(p_2, t_1+)|^2 = 0$ (here $t_1+$ indicates the submanifold immediately after surgery). Indeed, $p_2$ is simply the `tip' of the axially symmetric convex hypersurface inserted in the final step of the standard surgery procedure. Since surgery is performed on necks with radius close to $r_*$, inspection of the standard surgery procedure shows that $d_{g(t_1+)}(p_2,p_1) \leq (2+10\Lambda) r_*$. Let us consider the region $\mathcal B := \mathcal B_{g(t_1+)}(p_2, (3+10\Lambda)r_*)$. We know that $|H(\cdot, t_1+)|$ is comparable to $K_*$ in this region. We have two cases to consider.  

\textbf{Case 1}. Suppose there is a time $\bar t \in (t_1, t_0)$ at which some point in $\mathcal B$ is affected by a surgery. We may assume $\bar t$ is the first such time and let $\bar p \in \mathcal B$ lie on a neck where surgery is performed at time $\bar t$. Using the gradient estimates and the evolution equation $\partial_t g = - H\cdot A$ we obtain
    \[\dist_{g(\bar t)}(\bar p, p_2) \leq \dist_{g(t_1+)}(\bar p, p_2) + C(n,\alpha)\kappa r_* \leq (4+10\Lambda)r_*\]
if $\kappa$ is sufficiently small. Because \ref{s1} holds at $\bar t$ with $\varepsilon_* \leq \varepsilon_s$ we know that $\bar p$ lies on a normal $(\varepsilon_s,1)$-cylindrical submanifold neck in $\mathcal {M}_{\bar t}$ which is long enough to contain $p_2$. In particular,
    \[|A(p_2, \bar t)|^2 - \tfrac{1}{n} |H(p_2, \bar t)|^2 \geq -C\varepsilon_s r_*^{-2} + \left(\tfrac{1}{n-1} - \tfrac{1}{n}\right) r_*^{-2}/4.\]
On the other hand, because $p_2$ is an umbilic point at time $t_1$, using the gradient estimates we get $|A(p_2, \bar t)|^2 - \tfrac{1}{n} |H(p_2, \bar t)|^2 \leq C(n,\alpha)\kappa r_*^{-2}$. These two estimates cannot hold simultaneously if $\kappa$ and $\varepsilon_s$ are both small, so we have reached a contradiction. 

\textbf{Case 2}. Now suppose $\mathcal B$ does not contain surgeries for $t \in (t_1, t_0)$. In this case, using the gradient estimates and evolution of the metric as before, we find that 
    \begin{align*}
    \dist_{g(t_0)}(p_2, p_0) &\leq \dist_{g(t_0)}(p_2, p_1) + \dist_{g(t_0)}(p_1, p_0)\\
    &\leq \dist_{g(t_1+)}(p_2, p_1) + C(n,\alpha)\kappa r_* + 2\hat Lr_*\\
    &\leq (2+10\Lambda)r_* + 3\hat Lr_*\\
    &\leq 4\hat L r_*.
    \end{align*}
In particular, at time $t_0$ the point $p_2$ is part of the neck $\mathcal N_{0}$, so we obtain a contradiction if $\kappa$ and $\varepsilon_s$ are both small exactly as in Case 1. This concludes the proof of Claim~\ref{claim_small surgery free interval}.
\end{proof}

Because of Claim~\ref{claim_small surgery free interval}, there are no surgeries in $\mathcal B_{g(t_0)}(p_0, \hat L r_0)$ for $t \in (t_0 - \kappa r_0^2, t_0]$. To complete the proof of Lemma~\ref{lem_ss} we show that there are also no surgeries in $\mathcal B_{g(t_0)}(p_0, \hat L r_0)$ for $t \in (t_0 - \hat L^2 r_0^2, t_0 - \kappa r_0^2]$. We use the neck detection lemma to argue that in this region of spacetime the mean curvature is smaller than $K_*$ by a definite amount and hence there can be no surgeries there. 

Let us denote by $\theta$ the maximal value such that $\hat {\mathcal P}(p_0, t_0, \hat L, \theta)$ does not contain surgeries and set $\hat \theta := \min\{\theta, \hat L^2\}$. Our assumptions concerning $(p_0,t_0)$ mean that
    \[|H(p_0,t_0)| \geq (1+O(\varepsilon_s))K_s, \qquad  |A(p_0,t_0)|^2 - \tfrac{1}{n-1}|H(p_0,t_0)|^2 = O(\varepsilon_s)|H(p_0,t_0)|^2.\]
Therefore, given any $\tilde \varepsilon$ and $\hat L$, we can choose $\varepsilon_s$ small and $K_s$ large so that by the neck detection lemma $\hat{\mathcal P}(p_0, t_0, \hat L, \hat \theta)$ is an $(\tilde \varepsilon, 1)$-shrinking neck. In particular, the mean curvature at every point of $\hat {\mathcal P}(p, t_0, \hat L, \hat\theta)$ is close to that of a standard shrinking cylinder whose radius equals $r_0$ at time $t_0$, that is:
    \begin{align*}
    |H(q,t)| &= (n-1 + O(\tilde \varepsilon))(r_0^2 + 2(n-1)(t_0-t))^{-1/2}.
    \end{align*}
for each $(q,t) \in \hat {\mathcal P}(p_0, t_0, \hat L, \hat\theta)$.

With the aim of deriving a contradiction we assume $\hat \theta < \hat L^2$. It follows that $\hat \theta = \theta$. We then have
    \begin{align*}
    |H(q, t_0 -\theta r_0^2)|&= (1+O(\tilde \varepsilon))(1+2(n-1)\theta)^{-1/2} |H(p_0, t_0)|\\
    &= (1+O(\tilde  \varepsilon) + O(\varepsilon_s))(1+2(n-1)\theta)^{-1/2} K_*
    \end{align*}
for each $q \in B_{g(t_0)}(p_0, t_0, \hat L r_0)$. By the definition of $\theta$ we know that $B_{g(t_0)}(p_0, \hat L r_0)$ contains a point modified by surgery at time $t_0 - \theta r_0^2$. All surgeries prior to $t_0$ are performed on necks where the mean curvature satisfies $|H| = (1+O(\varepsilon_s))K_*$, and standard surgery only increases the mean curvature, so there is a point $q_0 \in B_{g(t_0)}(p_0, \hat L r_0)$ where $|H(q_0,t_0 - \theta r_0^2)| \geq (1+O(\varepsilon_s))K_*$. In this way we deduce
    \[(1+2(n-1)\theta)^{1/2} \leq 1+O(\tilde \varepsilon) + O(\varepsilon_s).\]
On the other hand, by Claim~\ref{claim_small surgery free interval} we have $\kappa \leq \theta$, which contradicts the previous inequality if $\varepsilon_s$ and $\tilde \varepsilon$ are sufficiently small. Therefore, $\hat \theta = \hat L^2$ and $\hat{\mathcal P}(p_0, t_0, \hat L, \hat L^2)$ is an $(\tilde \varepsilon, 1)$-shrinking neck. Choosing $\tilde \varepsilon \leq \hat \varepsilon$ completes the proof. 
\end{proof}

We conclude this section with a couple of geometric lemmas which will be used in conjunction with the neck detection lemma in the sequel. The first of these should be compared with \cite[Theorem~7.14]{HuSi09}. 

\begin{lemma}\label{lem_pinching implies compact}
Let $\mathcal M$ be a connected, complete, smoothly immersed $n$-submanifold of $\mathbb{R}^{n+m}$. Suppose there are constants $c^{\#}$, $K$ and $\sigma > 0$ such that:
    \begin{itemize}[itemsep=0.2cm]
    \item $|\nabla A(p)|\leq n^{-1/2}c^{\#}|H(p)|^2$ holds at every $p\in \mathcal M$, and
    \item $|A^-(p)|^2 \leq K R^{-2\sigma}|H(p)|^{2-2\sigma}$ holds at every $p \in \mathcal{M}$ with $|H(p)| > 0$. 
    \end{itemize}
For every $\eta > 0$ there exists a constant $\Upsilon = \Upsilon(n, m, c^{\#}, K, \sigma, \eta)$, with the following property. Given any $p \in \mathcal M$ such that $|H(p)| \geq (1+c^\#\Upsilon)R^{-1}$, either the inequality
    \[|A|^2 < (\tfrac{1}{n-1}-\eta)|H|^2\]
holds everywhere in $\mathcal{M}$, or else there exists a point $q \in \mathcal M$ with $d_g(q,p) \leq \Upsilon/|H(p)|$ such that 
    \[|A(q)|^2 \geq (\tfrac{1}{n-1} - \eta)|H(q)|^2\]
and moreover $|H(q')| \geq |H(p)|/(1+c^\#\Upsilon)$ whenever $q'$  satisfies $d_g(q', p) \leq d_g(q, p)$.
\end{lemma}
\begin{proof}
Given $\Upsilon > 0$ and $p \in \mathcal {M}$, let us denote $\mathcal M_{p,\Upsilon} := \{q \in \mathcal{M} : d_g(q,p) \leq \Upsilon/|H(p)|\}$. We recall from Lemma~\ref{lem_localHarnack} that, because of the gradient estimate, $|H(q)| \geq |H(p)|/(1+c^\#\Upsilon)$
holds automatically for all points $q \in \mathcal M_{p,\Upsilon}$. Therefore, to prove the claim it suffices to demonstrate that if $\Upsilon$ is large enough, whenever $|A|^2 < (\tfrac{1}{n-1}-\eta)|H|^2$ holds everywhere in $\mathcal{M}_{p,\Upsilon}$ for some $p$ with $|H(p)|\geq (1+c^\#\Upsilon)R^{-1}$, it must be the case that $\mathcal M = \mathcal{M}_{p,\Upsilon}$. We prove this statement by contradiction. 

Supposing the claim is false for some $\eta > 0$, given any $\Upsilon_i \to \infty$ there exists a sequence of immersions $\mathcal M_i \to \mathbb{R}^{n+m}$, each satisfying the hypotheses of the lemma with the same constants $c^{\#}$, $K$ and $\sigma$, with the following property: For every $i$ there is a point $p_i \in \mathcal M_i$ such that 
    \[|H(p_i)| \geq (1+c^\#\Upsilon_i)R^{-1} \qquad \text {and} \qquad |A_i|^2 < \left(\tfrac{1}{n-1}-\eta\right)|H_i|^2 \;\; \text{in} \;\; \mathcal{M}_{p_i,\Upsilon_i},\]
but $\mathcal M_{p_i, \Upsilon_i}$ is a strict subset of $\mathcal M_i$. To derive a contradiction we will rescale and extract a complete noncompact limit which is a uniformly convex hypersurface, thus violating Hamilton's pinching theorem \cite{Hamilton_pinching}. The points around which we rescale need to be chosen carefully, however, so that the rescaled submanifolds have uniformly bounded curvature. 

Since $\mathcal{M}_{p_i,\Upsilon_i}$ is strictly contained in $\mathcal M_i$ it has nonempty boundary. For $q \in \mathcal{M}_{p_i,\Upsilon_i}$ we define $d_i(q) := d_{g_i}(q, \partial \mathcal{M}_{p_i,\Upsilon_i})$ and set 
    \[L_i := \max_{q \in \mathcal{M}_{p_i,\Upsilon_i}} d_i(q) |A_i(q)|.\]
Using $d_i(p_i) = \Upsilon_i/|H_i(p_i)|$ we see that $L_i \geq d_i(p_i) |A_i(p_i)| \geq n^{-1/2} \Upsilon_i$ and therefore $L_i \to \infty$. Let $q_i \in \mathcal{M}_{p_i,\Upsilon_i}$ be such that $L_i = d_i(q_i) |A_i(q_i)|$ and observe that since 
    \[d_i(q_i) \leq d_{g_i}(q_i,p_i) + d_i(p_i) \leq \frac{2\Upsilon_i}{|H_i(p_i)|} \leq \frac{2 R}{\Upsilon_i^{-1} + c^\#}\]
we have $|A_i(q_i)| \geq L_i/d_i(q_i) \to \infty$. Using the triangle inequality again, we see that
    \begin{align*} \sup_{q \in \mathcal{B}_{g_i}(q_i, d_i(q_i)/2)} |A_i|(q) &\leq L_i \sup_{q \in \mathcal{B}_{g_i}(q_i, d_i(q_i)/2)} \frac{1}{d_i(q)}  \leq 2\frac{L_i}{d_i(q_i)} \leq 2|A_i|(q_i).
    \end{align*}
Therefore, the rescaled immersions $\tilde {\mathcal M}_i := |A_i(q_i)|(\mathcal M_i - q_i)$ are such that $|\tilde A_i(q)| \leq 2$ for every $q \in \mathcal{B}_{\tilde g_i}(q_i, L_i/2)$. Since $L_i \to \infty$, using the uniform gradient estimate and standard compactness theorems we conclude that $\tilde{\mathcal M}_i$ subconverges locally uniformly in $C^2$ to a complete noncompact limiting immersion $\tilde{\mathcal M}_\infty : \mathcal{M}_\infty \to \mathbb{R}^{n+m}$ with $|\tilde A_\infty| \leq 2$. The gradient estimate also implies $|\tilde H_\infty| > 0$ everywhere in $\mathcal M_\infty$ (see Lemma~\ref{lem_localHarnack}), and we have
    \[\frac{|\tilde A_i^-|^2}{|\tilde H_i|^{2}} \leq K (R|A_i(q_i)|)^{-2\sigma} \leq K \left(\frac{1+c^\#\Upsilon_i}{n^{1/2}}\right)^{-2\sigma}\to 0,\]
so $\tilde {\mathcal M}_\infty$ is a hypersurface in some $\mathbb{R}^{n+1} \subset \mathbb{R}^{n+m}$ by \cite[Proposition~2.5]{Naff_planarity}. Moreover, from the inequality $|\tilde A_\infty|^2 \leq (\tfrac{1}{n-1} - \eta)|\tilde H_\infty|^2$, which holds everywhere in $\tilde {\mathcal M}_{\infty}$ by assumption, we deduce that $\tilde{\mathcal M}_\infty$ is uniformly convex via the following computation:
    \begin{align}\label{eq_pinching implies convex}
    \begin{split}
    \frac{1}{n-1}|H|^2 - |h|^2 & \leq \bigg(\frac{1}{n-1}\lambda_1 + \frac{2}{n-1}\sum_{i=2}^n \lambda_i - \lambda_1\bigg) \lambda_1 \leq \frac{2}{n-1}|H|\lambda_1.
    \end{split}
    \end{align}
Given that $\tilde{\mathcal M}_\infty$ is complete and noncompact, this contradicts Hamilton's pinching theorem for hypersurfaces \cite{Hamilton_pinching}. 
\end{proof}

The next lemma should be compared with \cite[(iii) of Proposition~7.18]{HuSi09}. 

\begin{lemma}\label{lem_axis approx const}
Let $\mathcal {N} : \mathbb{S}^{n-1} \times [a,b] \to \mbb{R}^{n+m}$ be a normal $(\varepsilon, 1)$-cylindrical submanifold neck and suppose $\langle \nu, \omega\rangle + |\omega^-| < \varepsilon$ at every point of $\mathcal N$ for some $\omega \in \mbb{S}^{n+m-1}$, where $\nu$ is the principal normal and $\omega^- = \omega^\perp - \langle\nu,\omega\rangle\nu$. In addition, suppose every point in $\mathbb{S}^{n-1}\times[a+10,b-10]$ lies at the centre of an $(\tilde\varepsilon,1,5)$-neck, i.e.\ it is a cylindrical graph in the sense of Definition~\ref{def_neck}. Given any $\tilde \varepsilon > 0$, we can choose $\varepsilon$ sufficiently small so that if $\tilde \omega \in \mbb{S}^{n+m-1}$ is parallel to the axis of any one of these cylindrical graphs then $\sqrt{1-\langle\omega,\tilde\omega\rangle^2} = O(\tilde \varepsilon)$.
\end{lemma}
\begin{proof}
We follow Proposition~7.8 of \cite{HuSi09}. We may assume $\tilde \omega \not = \omega$, so the vector $v := \omega - \langle\omega,\tilde\omega\rangle\tilde\omega$ is nonzero. Moreover, $v$ is orthogonal to $\tilde \omega$. It follows that there is a point $q$ on the cylindrical graph and a normal vector $\tilde \nu(q)$ such that $|\tilde \nu(q) - \frac{v}{|v|}| = O(\tilde \varepsilon)$. We may assume $\langle \nu(q), \omega\rangle > -\tilde \varepsilon$, for if not then by taking $q'$ to be the point in the cylindrical graph antipodal to $q$ we obtain 
    \[\varepsilon + |\nu(q') + \nu(q)|  > \langle \nu(q'),\omega\rangle -\langle \nu(q') + \nu(q),\omega\rangle  = -\langle \nu(q),\omega\rangle > \tilde \varepsilon,\]
which is a contradiction if $\varepsilon$ is sufficiently small. Combining $\langle \nu(q), \omega\rangle > -\tilde \varepsilon$ with $\langle \nu, \omega\rangle + |\omega^-| < \varepsilon$ we get $|\omega^\perp| < \tilde \varepsilon + 2\varepsilon$, and since $v \cdot \omega = |v|^2 = 1-\langle\omega,\tilde\omega\rangle^2$ this implies
    \[\sqrt{1-\langle\omega,\tilde\omega\rangle^2} = \frac{v}{|v|}\cdot\omega + \left(\tilde\nu(q) - \frac{v}{|v|}\right)\cdot\omega + O(\tilde \varepsilon) = \tilde \nu(q) \cdot \omega + O(\tilde \varepsilon) = O(\varepsilon) + O(\tilde \varepsilon).\]
We can choose $\varepsilon< \tilde \varepsilon$ to arrange $\sqrt{1-\langle\omega,\tilde\omega\rangle^2} = O(\tilde \varepsilon)$ as claimed. 
\end{proof}


\section{Neck continuation and the surgery algorithm}\label{sec_surgeries}

This final section is devoted to the proof of our main theorem, Theorem~\ref{thm_main}, a more precise restatement of which follows.

When $T$ is a surgery time we write $H(\cdot, T-)$ and $H(\cdot,T+)$ for the mean curvature before and after surgery, respectively. 

\begin{theorem}\label{thm_mainone}
Suppose $\mathcal M_0 \in \mathcal C_{n,m}(R,\alpha)$ is a closed immersed submanifold satisfying $|A|^2 \leq R^{-2}$. There then exist curvature thresholds $H_1< H_2<H_3$ and a mean curvature flow with surgery starting from $\mathcal M_0$ with the following properties.
\begin{itemize}
\item Each surgery occurs at a time $T_i$ such that $\max |H(\cdot, T_i-)|= H_3$.
\item After surgery, all components of the manifold satisfy $\max |H(\cdot, T_i+)|\leq H_2$, except for some components that are diffeomorphic to $\mathbb{S}^n$, $\mathbb{S}^{n-1} \times  \mathbb{S}^1$ or $\mathbb{S}^{n-1} \simtimes \mathbb{S}^1$ and are discarded.
\item Each surgery is performed on a normal cylindrical submanifold neck at a cross-section where the mean radius equals $(n-1)/H_1$.
\item The flow with surgery terminates after finitely many steps.
\end{itemize}
The constants $H_i$ can be arbitrarily large but satisfy $H_1 \geq \omega_1 R^{-1}, H_2=\omega_2 H_1$ and $H_3 = \omega_3 H_2$ with $\omega_i>1$ depending only on $n$, $m$ and $\alpha$.
\end{theorem}

In order to determine suitable ranges for the curvature thresholds $H_i$ appearing in Theorem~\ref{thm_mainone} we first specify a number of parameters. These choices mirror those in \cite[Section~8]{HuSi09}, except that we additionally make use of the planarity estimate via our hypersurface detection lemma (see \ref{P7} below), and we need to ensure that surgery necks are separated from each other in spacetime (see \ref{P0} below). 

\begin{enumerate}[itemsep=0.2cm, label=(P\arabic*)]
\setcounter{enumi}{-1}

\item\label{P0} \emph{Separation between surgery necks:} We define $\hat \varepsilon$ and $\hat L$, depending only on $n$, as in Lemma~\ref{lem_separation implies s0}. In addition, we require $\hat L \geq 5 + 14\Lambda$. We then choose $\varepsilon_s = \varepsilon_s$ and $K_s$ so that Lemma~\ref{lem_ss} can be applied with these choices of $\hat \varepsilon$ and $\hat L$. 

\item\label{P1}\emph{Neck parameters:} In Section~\ref{sec_standard surgery} we defined a standard surgery procedure for normal $(\varepsilon,k)$-cylindrical submanifold necks where $k \geq k_0 =4$ and $0 < \varepsilon \leq \varepsilon_0$. This requires the length of the neck to satisfy $L \geq 14\Lambda$. We continue to assume $\varepsilon_0$ is small enough so that, whenever standard surgery is performed, Theorem~\ref{geom quants N hat} and Theorem~\ref{thm_invariant} are in effect. We now further assume $L \geq 100 + \hat L$ and that $\e_0$ is small enough so that on any normal $(\e_0,1)$-cylindrical submanifold neck of length $10L$ the mean curvature at any two points can differ at most by a factor of $11/10$. Moreover, we choose $\varepsilon_0 \leq \varepsilon_s$.

\item\label{P2} We define $c^\#$ as in Corollary \ref{cor_gradient}, $d^\#$ as in Lemma \ref{lem_Harnack2} and $\Theta = 10+10^5(n-1)/c^\#$.

\item\label{P3} \emph{First application of neck detection:} We choose $\eta_0$ and $K_0$ such that if
    \begin{equation}\label{eqn_firstneck}
    |H(p,t_0)|\geq K_0, \quad |A(p,t_0)|^2 - \tfrac{1}{n-1}|H(p,t_0)|^2\geq -\eta_0 |H(p,t_0)|^2
    \end{equation}
and $\hat{\mathcal P}(p,t_0, L',\theta')$ does not contain surgeries for some $L' \in [L/4,L]$ and $\theta' \in [d^{\#}/1600, d^\#]$, then the region $\hat{\mathcal P}(p,t_0,L',\theta')$ is an $(\e_0,k_0-1)$-shrinking neck and $(p,t_0)$ lies at the centre of a normal $(\e_0,k_0)$-cylindrical submanifold neck of length $2L'-2$ (see Lemma~\ref{lem_NDL} and Theorem~\ref{neck overlapping properties}). In addition, we require $\eta_0$ and $K_0$ to be such that if $(p,t_0)$ satisfies \eqref{eqn_firstneck} and $|H(p,t_0)|\geq 5K$, where $K$ bounds from above the mean curvature in every region changed by a surgery, then the conclusions of Lemma~\ref{lem_NDL2replace} apply. Finally, we require $\eta_0$ and $K_0$ to be such that Lemma~\ref{prop_discborder} can be applied to the parabolic neighbourhood $\hat{\mathcal P}(p,t_0,L',\theta')$ for $L'$ and $\theta'$ as above.

\item\label{P4} \emph{Second application of neck detection:} We set $\e_1= \frac{(n-1)}{50}\eta_0$ and apply Lemma~\ref{lem_NDL} to find $\eta_1$ and $K_1$ such that if $(p,t_0)$ satisfies
\[|H(p,t_0)|\geq K_1, \qquad |A(p,t_0)|^2 - \tfrac{1}{n-1} |H(p,t_0)|^2 \geq - \eta_1 |H(p,t_0)|^2\]
and the parabolic neighbourhood $\hat{\mathcal P}(p,t_0, 10,d^{\#}/1600)$ does not contain surgeries, then $(p,t_0)$ lies at the centre of an $(\varepsilon_1, 1, 9)$-neck (in the sense of Definition~\ref{def_neck}). Moreover, we assume $\eta_1$ is small enough to guarantee that $|\omega^\perp| < \varepsilon_1$ at every point of this neck, where $\omega$ is a unit vector generating its axis. We assume that $K_1\geq K_0$ and $\eta_1\leq\eta_0$.

\item\label{P5} \emph{Application of Lemma~\ref{lem_pinching implies compact}:} We choose $\Upsilon$ and set $\gamma_0 = 1+c^\#\Upsilon$ such that if $|H(p,t_0)|\geq\gamma_0 R^{-1}$ then either $|A|^2- \frac{1}{n-1} |H|^2< -\eta_1|H|^2$ on the whole component of $\mathcal M_{t_0}$ containing $p_0$, or else there exists a $q\in\mathcal {M}_{t_0}$ with $d_{g(t_0)}(q,p) \leq \Upsilon/|H(p,t_0)|$ such that 
    \[|A(q,t_0)|^2- \tfrac{1}{n-1}|H(q,t_0)|^2 \geq -\eta_1|H(q,t_0)|^2\]
and moreover $|H(q',t_0)| \geq |H(p,t_0)|/\gamma_0$ whenever $q'$  satisfies $d_{g(t_0)}(q', p) \leq d_{g(t_0)}(q, p)$.

\item\label{P6} \emph{Third application of neck detection:} We set $\theta_2 = (10^5(n-1)\gamma_0^2\Theta^2)^{-1}$ and choose $\eta_2$ and $K_2$ such that if
    \[|H(p,t_0)|\geq K_2, \qquad |A(p,t_0)|^2 - \tfrac{1}{n-1} |H(p,t_0)|^2 \geq -\eta_2 |H(p,t_0)|^2\]
and $\hat{\mathcal P}(p,t_0,10,\theta_2)$ does not contain surgeries, then $(p,t_0)$ lies at the centre of an $(\varepsilon_1, 1, 9)$-neck. We assume $\eta_2$ is small enough to guarantee that this neck has the following property: If $q$ is antipodal to $p$, in the sense that the line through $p$ tangent to $\nu(p, t_0)$ meets the neck at $q$, then we have $|\nu(p,t_0) + \nu(q, t_0)| \leq \varepsilon_1$. We assume $K_2\geq K_1$ and $ \eta_2 < \eta_1$.

\item\label{P7} \emph{Hypersurface Detection:} We set $\e_2 = 10^{-3}\min\{(n-1)\eta_2, \eta_2\varepsilon_1^2\}$ and choose $K_3$ such that if $|H(p,t_0)|\geq K_3$ and $\hat{\mathcal P}(p,t_0,10,\theta_2)$ does not contain surgeries then 
    \[|A^-(p,t_0)| + |\nabla \nu(p,t_0)| \leq \varepsilon_2|H(p,t_0)|.\]
We assume $K_3 \geq K_2$.  

\item\label{P8} We finally restrict $H_1$ to be larger than $\max\{4\Theta K_3, K_s, 2R^{-1}\}$ and then define $H_2 = 10 \gamma_0H_1$ and $H_3 = 10 H_2$. It is possible to choose $H_1$, and hence $H_2$ and $H_3$, arbitrarily large. 
\end{enumerate}

\begin{remark} The curvature thresholds $H_i, K_i$ can be written as constants depending only on $n, m, \alpha$ multiplied by $R^{-1}$, while the remaining parameters depend only on $n, m, \alpha$. 
\end{remark}

\begin{remark}
It is natural to ask what happens as we send $H_i\rightarrow \infty$. In \cite{Lauer2013} and \cite{Head2013} it was shown that in this asymptotic regime the codimension-one mean curvature flow with surgery approaches the level-set flow.
\end{remark}

A priori, surgery times may accumulate. Therefore, we construct our surgery algorithm so that the following property holds:

\begin{enumerate}[label=(\Alph*)]
\setcounter{enumi}{18}
\item\label{property S} We perform surgery at times $T_i$ such that $\max |H(\cdot, T_i-)|= H_3$. After performing surgery at time $T_i$ and discarding suitable components whose topology is known, we have $\max|H(\cdot,T_i+)|\leq H_2$. In addition, all surgeries satisfy \ref{s0} and \ref{s1}-\ref{s3}, with the parameters in \ref{s1} given by $\varepsilon_* = \varepsilon_0$, $K_* = H_1$ and $L_* = L$. 
\end{enumerate}
Property \ref{property S} and the inequality 
    \[(\partial_t - \Delta)|H|^2 \leq 2|A|^2 |H|^2 \leq \frac{2}{n-2}|H|^4\]
imply a uniform lower bound for difference between any two consecutive times, namely
    \begin{align} \label{eq_waiting time}
    T_{i+1} - T_i \geq \frac{n-2}{2}\left(1-\frac{H_2^2}{H_3^2}\right)\frac{1}{H_2^2} \geq \frac{n-2}{10^3\gamma_0^2}\frac{1}{H_1^2}.
    \end{align}

The results in Section~\ref{sec_ND} are sufficient to conclude that whenever the maximum curvature becomes large, our submanifold is either a positively curved immersed sphere or else contains necks of large curvature. In order to perform controlled surgeries on these necks which remove all regions of large curvature and maintain topological control, we study the geometry of the submanifold as we travel outwards from a neck in each direction. The following \emph{neck continuation theorem} says that either we find a cross-section where the curvature is smaller by a fixed factor, or else the neck eventually closes up, forming a cap. It is inspired by the analogous Theorem~8.2 in \cite{HuSi09}. However in higher codimensions a fundamental new difficulty arises: a~priori there is no reason for the submanifold to remain close to a hypersurface after we leave the neck. Ruling out this possibility requires careful use of our planarity estimate for the flow with surgery. 

\begin{theorem}[{Neck continuation, cf.\ \cite[Theorem 8.2]{HuSi09}}] \label{thm_neckcontinuation}
Let $\{\mathcal M_t\}_{t\in [0,t_0]}$ be a mean curvature flow with surgery in $\mathcal C_{n,m}(R,\alpha)$ such that $\max|A(\cdot,0)|^2 \leq R^{-2}$ and \ref{property S} holds. Suppose $\max|H(\cdot, t_0)|\geq H_3$. If $\eta_1$ and $H_1$ are defined as in \ref{P0}-\ref{P8} and
    \begin{align}\label{eqn_cylinder}
    |H(p_0,t_0)|\geq 10 H_1, \quad |A(p_0,t_0)|^2-\frac{1}{n-1} |H(p_0,t_0)|^2\geq - \eta_1 |H(p_0,t_0)|^2,
    \end{align}
then $(p_0,t_0)$ lies on a normal $(\e_0, k_0)$-cylindrical submanifold neck $\mathcal N_0$ which either covers the whole component of $\mathcal M_{t_0}$ containing $p_0$, or else has a boundary consisting of two cross-sections $\Sigma_1$ and $\Sigma_2$, each of which satisfies one of the following:
\begin{enumerate}[label=(\roman*)]
\item\label{NC_mean radius} The mean radius of $\Sigma$ is $2(n-1)/H_1$.
\item\label{NC_convex cap} $\Sigma$ is the boundary of a region $\mathcal D$ diffeomorphic to a standard $n$-ball. Moreover, $\mathcal D$ lies ``after" the cross-section $\Sigma$, i.e.\ it is disjoint from $\mathcal N_0$.
\end{enumerate}
\end{theorem}
\begin{proof}
Suppose $p_0$ is such that \eqref{eqn_cylinder} is satisfied. From our definitions we have
    \begin{align*}
    &|H(p_0,t_0)|\geq 10 K_1\geq 10 K_0,\\
    &|A(p_0,t_0)|^2-\frac{1}{n-1} |H(p_0,t_0)|^2 \geq - \eta_1 |H(p_0,t_0)|^2\geq -\eta_0 |H(p_0,t_0)|^2.
    \end{align*}
Therefore, at $(p_0,t_0)$ we can apply neck detection at the ``finer" $\e_1$-level and the ``coarser" $\e_0$-level, by \ref{P4} and \ref{P3} respectively. We begin our analysis at the $\e_0$-level. Since surgeries are performed on neck regions where the curvature is close to $H_1$, $2 H_1$ is a bound from above for the mean curvature in every region changed by a surgery. Therefore, by \ref{P3} we can apply Lemma~\ref{lem_NDL2replace} with $K = 2H_1$ to ensure that $\hat{\mathcal P}(p_0,t_0,L,d^\#)$ does not contain surgeries, and then conclude that $(p_0,t_0)$ lies at the centre of a normal $(\e_0,k_0)$-cylindrical submanifold neck of length $2L-2$. Let us denote by $\mathcal N_0$ the maximal normal $(\e_0,k_0)$-cylindrical submanifold neck containing $p_0$. If $\mathcal N_0$ covers the entire component of $\mathcal M_{t_0}$ containing $p_0$ then we are done. If not, we follow the neck parameter $z$ in both directions from $p_0$ until we find cross-sections of $\mathcal N_0$ satisfying either \ref{NC_mean radius} or \ref{NC_convex cap}. 

Let us assume $p_0$ lies in the cross-section $z=0$, and follow the neck in the direction of increasing $z$. If there exists a cross section with mean radius $2(n-1)/H_1$, then \ref{NC_mean radius} is satisfied and so we are done. Therefore, let us assume $r(z) < 2(n-1)/ H_1$ for each $z \in [0,z_{\max}]$, where $z_{\max}$ corresponds to the last cross-section of the neck. This implies $|H|> H_1/4$ for $z \in [0,z_{\max}]$. We need to show that \ref{NC_convex cap} holds in this case. Let $\Omega$ denote the set of all $p \in \mathcal{N}_0$ with $z \in [0,z_{\max}]$ such that
    \begin{equation}\label{eq_Omega1}
    |A(p,t_0)|^2 - \frac{1}{n-1} |H(p,t_0)|^2\geq -\eta_0 |H(p,t_0)|^2 \tag{$\Omega$1}
    \end{equation}
and
    \begin{equation}\label{eq_Omega2}
    \mathcal P \left(p,t_0, \frac{n-1}{|H(p,t_0)|}L, \frac{(n-1)^2}{(10 H_1)^2} d^{\#} \right) \text{ does not contain surgeries.} \tag{$\Omega$2}
    \end{equation}
If $p$ satisfies \eqref{eq_Omega1} but not \eqref{eq_Omega2}, then \ref{P3} and Lemma~\ref{lem_NDL2replace} guarantee $|H(p,t_0)|< 5K = 10 H_1$. In particular, we see that $p_0\in \Omega$. 
    
\begin{claim}
Every $p \in \Omega$ lies at the centre of a normal $(\e_0,k_0)$-cylindrical submanifold neck of length $2L-2$.
\end{claim}
\begin{proof}
Given that points $p\in \mathcal N_0$ with $z \in [0,z_{\max}]$ satisfy $|H(p,t_0)|\geq H_1/4$, we have
    \begin{align*}
    \frac{(n-1)^2}{1600 |H(p,t_0)|^2} \leq \frac{(n-1)^2}{(10 H_1)^2},
    \end{align*}
which implies
    \begin{align*}
    \hat{\mathcal P}\left(p,t_0,L,\frac{d^{\#}}{40^2}\right) \subset \mathcal P \left(p,t_0, \frac{n-1}{|H(p,t_0)|}L, \frac{(n-1)^2}{(10 H_1)^2} d^{\#} \right).
    \end{align*}
Therefore, by \ref{P3}, every $p \in \Omega$ lies at the centre of a normal $(\e_0,k_0)$-cylindrical submanifold neck of length $2L-2$. 
\end{proof}

Since the neck ends when $z=z_{\max}$, there are no points of $\Omega$ in $(z_{\max}-(L-1), z_{\max}]$. We let $z^*$ be the maximal value of $z$ with the following property: the cross-section of $\mathcal N_0$ with coordinate $z^*$ contains a point $p_1\in \Omega$, but there are no points of $\Omega$ with $z \in (z^*, z^*+10)$. Clearly then we have $z^* \leq z_{\max} - (L-1)$. We consider two cases:
    \begin{enumerate}[label=(\alph*)]
    \item There exists at least one point $p_2$ with $z \in (z^*,z^*+10)$ which satisfies \eqref{eq_Omega1}.
    \item All points with $z\in (z^*,z^*+10)$ fail to satisfy \eqref{eq_Omega1}.
    \end{enumerate}

\textbf{Case (a)}. Here the argument is exactly the same as in \cite{HuSi09}. Since $p_2$ does not satisfy \eqref{eq_Omega2} we can find a neighbourhood $\hat{\mathcal P}(p_2,t_0,L',\theta')$ with $L' \in [L/4,L]$ and $d^{\#}/40^2 \leq \theta ' \leq d^{\#}$ which is adjacent to a surgery in the sense of Lemma~\ref{prop_discborder}. Applying that lemma as in \ref{P3} we find the desired cross-section of $\mathcal N_0$ such that \ref{NC_convex cap} holds.

\textbf{Case (b)}. We assume $|A|^2 - \frac{1}{n-1}|H|^2 < -\eta_0 |H|^2$ at every point of $\mathcal N_0$ with $z \in (z^*,z^*+10)$. Here again we want to establish \ref{NC_convex cap}. At this point in the proof the fact that we are in codimension $m \geq 2$ presents serious new obstacles not present in the hypersurface case, whereas until now we have essentially reproduced the argument from \cite{HuSi09}. As we leave the neck $\mathcal N_0$ our submanifold could in principle bend into the extra ambient dimensions, failing to lie close to a hypersurface in some plane. This possibility is encoded in additional terms involving $A^-$ and $\nabla \nu$ in the evolution equations \eqref{eq_omega nu_1 evol}, \eqref{eq_omega top evol} and \eqref{eq_omega - evol} which we analyse below. These additional terms need to be carefully controlled using the planarity estimate via \ref{P7}.

We let $\overline{z} \in [0,z^*]$ be the largest value of $z$ such that, for some $\bar q$ in the corresponding cross-section of $\mathcal{N}_0$, we have $|A(\bar q, t_0)|^2 - \frac{1}{n-1} |H(\bar q, t_0)|^2 \geq -\eta_1 |H(\bar q, t_0)|^2$.

\begin{claim}
The parabolic neighbourhood $\hat{\mathcal P}(\overline q, t_0, 10, d^\#/1600)$ does not contain surgeries.  
\end{claim}
\begin{proof}
By the definition of $z^*$, there exists a point $q\in\Omega$ whose $z$-coordinate satisfies $z \in [0, z^*]$ and $|\bar z - z|\leq 10$, so using $|H(q,t_0)|\leq 2|H(\overline q, t_0)|$, $L \geq 20$ and $|H(\bar q, t_0)| \geq H_1/4$ we obtain
    \begin{align*}
    \hat{\mathcal P}\left(\overline{q}, t_0, 10, \frac{d^\#}{1600} \right) & \subset \mathcal{P}\left(q,t_0, \frac{n-1}{|H(q,t_0)|}L, \frac{(n-1)^2}{(10H_1)^2} d^\# \right).
    \end{align*}
The right-hand side does not contain surgeries by the definition of $\Omega$, so the left-hand side does not contain surgeries either. 
\end{proof}

As a consequence of the above claim and \ref{P4}, there exists a region $\mathcal{G} \subset \mathcal N_0$ centred at $\overline{q}$ which is an $(\varepsilon_1,1,9)$-neck. Moreover, $|\omega^\perp|<\varepsilon_1$ on $\mathcal G$, where $\omega$ is the unit vector parallel to the axis of the neck and pointing in the direction of increasing $z$. Up to a rotation of the ambient Euclidean coordinates we may assume $\omega = e_{n+1}$. Let us write $y = x_{n+1}$. By translating $\mathcal{M}_{t_0}$ in the direction of $e_{n+1}$, we can arrange that the maximum value of $y$ over the cross-section $z = \bar z$ is zero. This ensures that $\bar q$, and in addition every other point of $\mathcal{N}_0$ where $|A|^2 - \tfrac{1}{n-1}|H|^2 \geq -\eta_1 |H|^2$, lies in the region of $\mathbb{R}^{n+m}$ where $y \leq 0$. We denote by $\Sigma_0$ the intersection of the subspace $\{y = 0\}$ and $\mathcal G$. For each
$p\in\Sigma_0$ we consider the curve $y\mapsto \gamma(y,p)$ in $\mathcal{M}_{t_0}$ which solves
    \begin{align}\label{eqn_trajectories}
    \frac{d}{dy} \gamma = \frac{\omega^\top}{|\omega^\top|^2}
    \end{align}
with initial value $\gamma(0,p) = p$. Each trajectory $\gamma$ is well defined for as long as $|\omega^\top|$ remains positive. We let $y_{\max}$ denote the largest value such that for every $p \in \Sigma_0$ the curve $\gamma(y,p)$ is well defined for all $0 \leq y < y_{\max}$. We set $\Sigma_y=\{\gamma(y,p): p \in \Sigma_0\}$ for $0\leq y < y_{\max}$. Given $0\leq y_1 < y_2 < y_{\max}$, we set
    \begin{align*}
    \Sigma(y_1,y_2)= \bigcup \{\Sigma_y: y_1\leq y \leq y_2\}.
    \end{align*}
   
Let us denote by $\mathcal N_0'$ the part of $\mathcal N_0$ corresponding to $z\in [\overline{z}, z^*+10]$. We briefly recall the important properties of $\mathcal N_0'$.
\begin{itemize}
    \item By the definition of $\overline{z}$, $|A|^2 - \frac{1}{n-1} |H|^2 \leq -\eta_1 |H|^2$ holds everywhere in $\mathcal N'_0$.
    \item By the definition of $z^*$, $|A|^2 - \frac{1}{n-1} |H|^2 \leq -\eta_0 |H|^2$ in the final portion of $\mathcal {N}_0'$ where $z \in [z^*, z^* + 10]$.
    \item Because of how we chose $\eta_0$ and $\eta_1$, it is clear that in practice $\eta_1 \ll \eta_0$.
    \item $|\omega^\perp| < \e_1$ holds everywhere in $\mathcal G$ and so in particular holds at every point of $\Sigma_0$.
\end{itemize}

Let $\omega^- := \omega^\perp - \langle\nu,\omega\rangle\nu$. To show that $\mathcal M_{t_0}$ is closing up as we follow the trajectories $\gamma$, we first compute how the components of $\omega$ evolve. For $0 \leq y < y_{\max}$ we have 
    \begin{align}\label{eq_omega nu_1 evol}
    \begin{split}
    \frac{d}{dy} \langle \nu, \omega \rangle &= - h(\dot \gamma, \omega^\top) + \langle D_{\dot \gamma} \nu, \omega^\perp\rangle=-h\left(\frac{\omega^\top}{|\omega^\top|},\frac{\omega^\top}{|\omega^\top|}\right) + \frac{1}{|\omega^\top|^2}\langle\nabla_{\omega^\top} \nu, \omega^-\rangle.
    \end{split}
    \end{align}
In addition,
    \begin{align}\label{eq_omega top evol}
    \begin{split}
    \frac{d}{dy}|\omega^\top| &=\frac{1}{|\omega^\top|}\langle D_{\dot \gamma}^\top \omega^\top, \omega^\top\rangle = -\frac{1}{|\omega^\top|}\langle D_{\dot \gamma}^\top \omega^\perp, \omega^\top\rangle = \frac{1}{|\omega^\top|}\left\langle A\left(\frac{\omega^\top}{|\omega^\top|},\frac{\omega^\top}{|\omega^\top|}\right), \omega^\perp\right\rangle.
    \end{split}
    \end{align}
These two formulae combine to give
    \begin{align}\label{eq_omega - evol}
    \begin{split}
    \frac{d}{dy} |\omega^-| &= -\left\langle A^-\left(\frac{\omega^\top}{|\omega^\top|},\frac{\omega^\top}{|\omega^\top|}\right),\frac{\omega^-}{|\omega^-|}\right\rangle - \frac{\langle\nu,\omega\rangle}{|\omega^\top|^2}\left\langle\nabla_{\omega^\top} \nu, \frac{\omega^-}{|\omega^-|}\right\rangle.
    \end{split}
    \end{align}
We recall from \eqref{eq_pinching implies convex} the estimate 
    \begin{equation}\label{eq_pinching implies convex 2}
    h(X,X) \geq \frac{n-1}{2|H|}\left(\frac{1}{n-1}|H|^2 - |A|^2\right)|X|^2,
    \end{equation}
which is valid for all tangent vectors $X$. In the following claim we use these formulae to derive a coarse estimate for the normal components of $\omega$ in $\mathcal {N}_0'$. 

\begin{claim}\label{claim_coarse est omega nu_1}
Given any $Q \in [0,100]$ and a trajectory $\gamma(\cdot,p)$, for as long as the trajectory remains in $\mathcal{N}_0'$, the quantity $\langle\nu,\omega\rangle+Q|\omega^-|$ is strictly decreasing with respect to $y$. In particular, $\langle\nu,\omega\rangle + Q|\omega^-| < (1+Q)\varepsilon_1$ on each trajectory for as long it remains in $\mathcal{N}_0'$.
\end{claim}
\begin{proof}
Let $z'' \in [\bar z, z^*+10]$ be the largest value such that $|\omega^\top| \geq 1/2$ on every cross-section $z \in [\bar z, z'']$, and denote by $\mathcal{N}_0''$ the portion of $\mathcal{N}_0'$ where $z \in [\bar z, z'']$. Using \eqref{eq_pinching implies convex 2} to estimate the first term on the right-hand side of \eqref{eq_omega nu_1 evol}, we obtain
    \begin{equation}\label{eq_omega nu_1 evol 2}
    \frac{d}{dy} \langle \nu, \omega \rangle \leq \frac{n-1}{2|H|} \left(|A|^2 - \frac{1}{n-1}|H|^2\right) + \frac{|\nabla \nu|}{|\omega^\top|}.
    \end{equation}
We have $|A|^2 - \tfrac{1}{n-1}|H|^2 \leq -\eta_1|H|^2$ at every point of $\mathcal{N}_0'$, and moreover $|\nabla \nu|\leq \e_2 |H|$ by \ref{P7}, so this gives
    \[\frac{d}{dy} \langle \nu, \omega \rangle \leq -\left(\frac{n-1}{2}\eta_1-2\e_2\right)|H|\]
on each trajectory for as long as it remains in $\mathcal{N}_0''$. Next we use \eqref{eq_omega - evol} and \ref{P7} to estimate
    \begin{equation}\label{eq_omega - evol 2}
    \frac{d}{dy}|\omega^-| \leq |A^-| +2|\nabla \nu_1|\leq 3\e_2|H|
    \end{equation}
at points in $\mathcal{N}_0'$. Given that $\e_2 \leq \tfrac{(n-1)}{1000}\eta_1$ this implies 
    \[\frac{d}{dy}(\langle\nu,\omega\rangle + Q|\omega^-|) \leq -\left(\frac{n-1}{2}\eta_1-302\e_2\right)|H| \leq -\eta_1|H|,\]
so $\langle\nu,\omega\rangle + Q|\omega^-|$ can only decrease from its inital value at $y=0$, and hence we have 
    \[\langle\nu,\omega\rangle + Q|\omega^-| < (1+Q)\e_1\]
on any given trajectory for as long as it remains in $\mathcal{N}_0''$.
This estimate holds for every $Q \in [0,100]$. In particular, $\langle\nu,\omega\rangle + |\omega^-| < 2\e_1$ holds everywhere in $\mathcal{N}_0''$.

To complete the proof we show that $\mathcal N_0'' = \mathcal N_0'$. By Lemma~\ref{lem_axis approx const}, we may assume that every portion of $\mathcal N_0$ with scale-free length 10 and which intersects $\mathcal N_0''$ has axis approximately equal to $\omega$. In particular, we may assume $\varepsilon_1$ is small enough so that $|\omega^\top| > 1/2$ everywhere in $\mathcal{N}_0''$. It follows that $\mathcal{N}_0'' = \mathcal{N}_0'$. 
\end{proof}

As we saw in the proof of Claim~\ref{claim_coarse est omega nu_1}, $|\omega^\top| \geq 1/2$ at every point of $\mathcal{N}_0'$. This ensures that the trajectories $\gamma$ are defined for at least as long as they remain within $\mathcal{N}_0'$, and hence there exists a smallest value $y'< y_{\max}$ such that $\gamma(y',p)\in \partial \mathcal N_0'$ for some $p\in \Sigma_0$. Next we get estimates at $\Sigma_{y'}$ using the fact that $|A|^2 - \tfrac{1}{n-1}|H|^2 \leq -\eta_0|H|^2$ on the final portion of $\mathcal N_0'$ where $z \in [z^*, z^* + 10]$. In particular, the following claim implies that on $\Sigma_{y'}$ we have $\langle \nu, \omega\rangle < -4\varepsilon_1$ and that $|\omega^-|$ is much smaller than $-\langle\nu,\omega\rangle$.

\begin{claim}\label{claim_fine est omega nu_1}
On $\Sigma_{y'}$ we have $\langle\nu,\omega\rangle + 100|\omega^-|\leq-4\e_1$. 
\end{claim}
\begin{proof}
We denote by $r^*=r(z^*)$ the mean radius of the cross-section $z = z^*$, and set $H^*=\frac{n-1}{r^*}$. We have $H^*\geq H_1/2$ and because the axis of the neck is close to $\omega$ (see Lemma~\ref{lem_axis approx const}), the coordinate $y$ is almost constant on each cross-section $z \in [z^*, z^* + 10]$. Moreover, the $y$-coordinates of the cross-sections $z=z^*$ and $z=z^*+10$ differ by approximately $10r^*$. It follows that at least the points of $\Sigma(y'-5r^*,y')$ lie in the region where $z\in[z^*,z^*+10]$. In addition, since $|H|$ varies slowly on the neck $\mathcal N_0$ we have $|H|\geq H^*/2 \geq H_1/4$ on $\Sigma(y'-5r^*,y')$. 

Using \eqref{eq_omega nu_1 evol 2}, \ref{P7} and the fact that $|\omega^\top| \geq 1/2$ in $\mathcal N_0'$, for $y\in[y'-5r^*,y']$ we obtain
    \[\frac{d}{dy} \langle \nu, \omega \rangle \leq -\left(\frac{n-1}{2}\eta_0-2\e_2\right)|H|.\]
Combining this inequality with \eqref{eq_omega - evol 2}, we see that the function $\psi(y):= \langle \nu(\gamma(y,p)),\omega\rangle + 100|\omega^-(\gamma(y,p))|$ satisfies 
    \[\frac{d}{dy}\psi \leq -\left(\frac{n-1}{2}\eta_0-302\e_2\right)|H| \leq -\eta_0|H| \leq -\frac{\eta_0}{2}H^*\]
for $y\in[y'-5r^*,y']$, where we have used $\varepsilon_2 \leq \frac{(n-1)}{1000}\eta_0$. Consequently, for every $p \in \Sigma_{0}$,
    \begin{align*}
    \psi(y')&= \psi(y'-5r^*) + \int_{y'-5r^*}^{y'}\frac{d}{dy} \psi \, dy \leq 101\e_1 - \frac{5}{2} \eta_0 r^* H^* \leq -4\e_1,
    \end{align*}
where we have used Claim~\ref{claim_coarse est omega nu_1} and $\varepsilon_1 = \frac{(n-1)}{50}\eta_0$. This proves the claim.
\end{proof}

We claim that as $y \to y_{\max}$ the trajectories $\gamma$ converge to a point in the ambient space. To prove this we need to understand how they behave after leaving the neck $\mathcal{N}_0$. Let $y^\#$ be the largest such value such that on all of the trajectories $\gamma$ the properties 
    \begin{enumerate}[itemsep=0.2cm, label=(\alph*)]
    \item\label{cap est a} $|\omega^\perp|< 1$,
    \item\label{cap est b} $|A|^2 - \frac{1}{n-1}|H|^2 < -\eta_2 |H|^2$,
    \item\label{cap est c} $\langle \nu, \omega\rangle + 100|\omega^-|< -2\varepsilon_1$ and
    \item\label{cap est d} $|H| \geq H_1/4\Theta$
    \end{enumerate}
hold for all $y' \leq y < y^\#$, with $\eta_2$ as in \ref{P6}. We will prove that $y^\# = y_{\max}$. Note that by the definition of $\theta_2$ in \ref{P6} if $|H(p,t_0)| \geq H_1/4\Theta$ then
    \begin{align*}
    \theta_2 \frac{(n-1)^2}{|H(p,t_0)|^2} \leq \theta_2 \frac{16 (n-1)^2 \Theta^2}{H_1^2} < \frac{n-2}{10^3\gamma_0^2H_1^2}.
    \end{align*}
Therefore, by the estimate \eqref{eq_waiting time} for the time between surgeries, \ref{cap est d} implies $ \hat{\mathcal P}(p,t_0,10,\theta_2)$ does not contain surgeries for every $p \in \Sigma(y', y^\#)$. We need the following statement. 

\begin{claim}\label{Sigma contracts}
The intrinsic diameter of $\Sigma_y$ decreases with respect to $y$ for $y \in [y', y^\#)$. In particular, the intrinsic diameter of $\Sigma_{y}$ never exceeds $4(n-1)/H_1$ for $y \in [y', y^\#)$. 
\end{claim}
\begin{proof}
Let us write $\hat g = \hat g(y)$ and $\hat A = \hat A(y)$ for the first and second fundamental forms of $\Sigma_y$. Pulling back by the diffeomorphism $\gamma(y, \cdot) : \Sigma_0 \to \Sigma_y$, we may view these tensors as being defined on $\Sigma_0$. Let $X$ be a vector tangent to $\Sigma_0$. The family of hypersurfaces $\Sigma_y$ moves with normal velocity $\omega^\top/|\omega^\top|^2$ in $\mathcal M_{t_0}$, so by the usual first variation formula for the induced metric we have 
    \begin{align*}
        \frac{\partial}{\partial y} \hat g(X,X) = - \frac{2}{|\omega^\top|^2} g(\omega^\top, \hat A(X,X)). 
    \end{align*}
Using $0 = D\omega = D\omega^\top + D\omega^\perp$ we may rewrite this as
    \begin{align*}
        \frac{\partial}{\partial y} \hat g(X,X)  = \frac{2}{|\omega^\top|^2} g(X, D_X \omega^\top) = -\frac{2}{|\omega^\top|^2} \langle X, D_X \omega^\perp\rangle = \frac{2}{|\omega^\top|^2} \langle \omega^\perp, A(X,X)\rangle.
    \end{align*}
We expand 
    \[\langle \omega^\perp, A(X,X)\rangle = \langle \nu,\omega\rangle h(X,X) + \langle \omega^-, A^-(X,X)\rangle.\] 
Because of \ref{cap est b} and \eqref{eq_pinching implies convex 2}, for $y \in [y', y^\#)$ we have 
    \[ h(X,X) \geq \frac{n-1}{2}\eta_2 |H| |X|^2,\]
and by \ref{cap est c} and \ref{P7} we have
    \[\langle \omega^-, A^-(X,X)\rangle \leq -\frac{\varepsilon_2}{100}\langle\nu,\omega\rangle |H||X|^2,\]
so we conclude that 
    \begin{align*}
       \frac{\partial}{\partial y} \hat g(X,X)  \leq \frac{2}{|\omega^\top|^2} \bigg(\frac{n-1}{2}\eta_2 - \frac{\varepsilon_2}{100}\bigg)\langle\nu,\omega\rangle|H||X|^2.
    \end{align*}
Our choice of $\varepsilon_2$ ensures that the right-hand side is non-positive i.e.\ $\hat g(X,X)$ is nonincreasing. We can apply this with $X = \dot c$ for any curve $c$ in $\Sigma_0$ to see that the length of the curve with respect to $\hat g$ is nonincreasing. Therefore, the intrinsic diameter of $\Sigma_y$ is nonincreasing. Recalling that $\Sigma_{y'}$ is close to the final cross-section of $\mathcal N_0'$, where the mean curvature is at least $H_1/2$, we see that the intrinsic diameter of $\Sigma_{y'}$ can be no larger than $4(n-1)/H_1$.
\end{proof}

Using Claim~\ref{Sigma contracts} we can bound the mean curvature from below in the region $\Sigma(y',y^\#)$. The bound we get is better than \ref{cap est d}. 

\begin{claim}\label{claim H lower bound cap}
In the region $\Sigma(y',y^\#)$ we have $|H| \geq H_1/2\Theta$.
\end{claim}
\begin{proof}
If the infimum of $|H|$ over the region $\Sigma(y',y^\#)$ is at least  $H_1/\Theta$ then there is nothing to prove, so suppose it is less than $H_1/\Theta$. Then since $|H| \geq H_1/2$ on $\Sigma_{y'}$, we can choose $\tilde y \in [y', y^\#)$ to be the smallest value such that $|H(\tilde p, t_0)| = H_1/\Theta$ for some $\tilde p \in \Sigma_{\tilde y}$. Our gradient estimate $|\nabla H| \leq c^\# |H|^2$ implies
    \[\frac{H_1}{2\Theta} \leq |H(p,t_0)| \leq \frac{2H_1}{\Theta}\]
for all $p \in \mathcal B_{g(t_0)}(\tilde p, \Theta/2c^\# H_1)$. In particular, 
    \[|A(p,t_0)| \leq \frac{2\sqrt{c_n} H_1}{\Theta}\]
for all $p \in \mathcal B_{g(t_0)}(\tilde p, \Theta/2c^\# H_1)$. Using this bound for the second fundamental form, elementary arguments show that the much smaller ball $\mathcal B_{g(t_0)}(\tilde p, \Theta/10^4 c^\# H_1)$ is a graph over the tangent space to $\mathcal M_{t_0}$ at $\tilde p$. In particular, $\mathcal B_{g(t_0)}(\tilde p, \Theta/10^4c^\# H_1)$ is diffeomorphic to a standard $n$-ball. 

By Claim~\ref{Sigma contracts}, the intrinsic diameter of $\Sigma_{\tilde y}$ is no larger than $4(n-1)/H_1$. Since we chose $\Theta \geq \frac{10^5(n-1)}{c^\#}$, this implies
    \[\Sigma_{\tilde y} \subset \mathcal B_{g(t_0)}(\tilde p, \Theta/10^4 c^\# H_1).\]
The region on the right is diffeomorphic to an $n$-ball, and $\Sigma_{\tilde y}$ is an embedded $(n-1)$-sphere, so it follows that $\Sigma_{\tilde y}$ bounds an $n$-ball in $\mathcal B_{g(t_0)}(\tilde p, \Theta/10^4 c^\# H_1)$, denoted $\tilde {\mathcal B}$. Note that in $\tilde{\mathcal B}$ we have $\frac{H_1}{2\Theta} \leq |H| \leq \frac{2 H_1}{\Theta}$. 

The region $\Sigma(\tilde y, y^\#)$ is either inside or outside of $\tilde{\mathcal B}$. If it is outside then $\Sigma(y', \tilde y)$ is inside, but since $|H| \geq H_1/2$ on $\Sigma_{y'}$ this is incompatible with the bound $|H| \leq 2H_1/\Theta$ in $\tilde{\mathcal B}$. Therefore, $\Sigma(\tilde y, y^\#)$ is inside $\tilde{\mathcal B}$ and hence $|H| \geq H_1/2\Theta$ holds there. Using the definition of $\tilde y$ we conclude that $|H| \geq H_1/2\Theta$ holds everywhere in $\Sigma(y',y^\#)$. 
\end{proof}

We are now ready to prove that $y^\# = y_{\max}$. Suppose, with the aim of deriving a contradiction, that $y^\# < y_{\max}$. Then since \ref{cap est a} holds for $y < y_{\max}$ by definition and \ref{cap est d} holds at $y^\#$ by Claim~\ref{claim H lower bound cap}, either \ref{cap est b} or \ref{cap est c} must fail at $y^\#$. Together the following two claims provide the desired contradiction. 

\begin{claim}\label{claim_c holds up to ymax}
Property \ref{cap est c} holds at $y^\#$. 
\end{claim}
\begin{proof}
From \eqref{eq_omega nu_1 evol}, \eqref{eq_omega - evol}, \eqref{eq_pinching implies convex 2} and \ref{P7} we see that for $y \in [y', y^\#)$ we have
    \[\frac{d}{dy}(\langle\nu,\omega\rangle + 100|\omega^-|) \leq 101\varepsilon_2 \frac{|H|}{|\omega^\top|}.\]
To bound $\langle\nu,\omega\rangle + 100|\omega^-|$ at $y^\#$ using the fundamental theorem of calculus, we need an estimate for the integral of the right-hand side. Since $y^\#$ may be close to $y_{\max}$, we cannot appeal to a positive lower bound for $|\omega^\top|$ as we did in our analysis of the region $y \leq y'$.

To bound the integral of $|H|/|\omega^\top|$, we first observe that by \eqref{eq_omega top evol}, \eqref{eq_pinching implies convex 2} and \ref{P7} we have
    \[\frac{d}{dy}|\omega^\top|^2 \leq \frac{(n-1)}{|H|}\left(\frac{1}{n-1}|H|^2 -|A|^2\right)\langle\nu,\omega\rangle + 2\e_2|\omega^-||H|.\]
Since \ref{cap est b} and \ref{cap est c} hold for $y < y^\#$, this shows that
    \begin{align*}
    2\frac{d}{dy}|\omega^\top| \leq ((n-1)\eta_2 - \e_2/50)\langle\nu,\omega\rangle\frac{|H|}{|\omega^\top|} \leq -\eta_2 \e_1 \frac{|H|}{|\omega^\top|}.
    \end{align*}
Given that $|\omega^\top| \leq 1$, by integrating with respect to $y$ we obtain
    \begin{equation}\label{eq_integral H / omega top}
    \int_{y'}^{y^\#}\frac{|H|}{|\omega^\top|} \leq \frac{2}{\eta_2 \e_1}.
    \end{equation}
    
Now, by integrating
    \[\frac{d}{dy}(\langle\nu,\omega\rangle + 100|\omega^-|) \leq 101\varepsilon_2 \frac{|H|}{|\omega^\top|},\]
and using Claim~\ref{claim_fine est omega nu_1} and \eqref{eq_integral H / omega top} we obtain the estimate
    \[\langle\nu ,\omega\rangle + 100|\omega^-|\leq -4\e_1 + 202\frac{\varepsilon_2}{\eta_2\varepsilon_1}\]
for every $y \in [y', y^\#)$. We chose $\varepsilon_2 \leq \eta_2 \varepsilon_1^2/1000$ so it follows that \ref{cap est c} cannot fail at $y^\#$.   
\end{proof}

\begin{claim}\label{claim_b holds up to ymax}
Property \ref{cap est b} holds at $y^\#$.   
\end{claim}
\begin{proof}
Suppose $|A(p^\#,t_0)|^2 -\frac{1}{n-1}|H(p^\#,t_0)|^2 = -\eta_2 |H(p^\#,t_0)|^2$ holds for some $p^\#\in \Sigma_{y^\#}$. By \ref{P6}, the point $p^\#$ lies on an $(\varepsilon_1,1,9)$-neck. Considering the point $q^\#$ of this neck which is antipodal to $p^\#$, by \ref{P6} we have 
    \[\langle \nu(q^\#),\omega\rangle = \langle \nu(q^\#) + \nu(p^\#),\omega\rangle - \langle \nu(p^\#),\omega\rangle \geq -\varepsilon_1 + 2\varepsilon_1 > 0,\]
which contradicts \ref{cap est c}. This completes the proof of Claim~\ref{claim_b holds up to ymax}.    
\end{proof}

It remains to show that the region traced out by the trajectories $\gamma$ is a standard $n$-ball. This follows from basic Morse theory. As $y \to y_{\max}$, at least one of the trajectories $\gamma$ reaches a point $\hat p$ in $\mathcal M_{t_0}$ where $\omega^\top(\hat p) = 0$. In particular, $\hat p$ is a critical point of the function $y$. Because of \ref{cap est b} and \ref{cap est c} this critical point is a nondegenerate maximum and hence the slices $\Sigma_y$ contract to $\hat p$ as $y \to y_{\max}$. Since there are no other critical points of $y$ in $(0, y_{\max})$, this region is diffeomorphic to a standard $n$-ball.
\end{proof}

With the neck continuation in place we now prove Theorem~\ref{thm_mainone}.

\begin{proof}[Proof of Theorem \ref{thm_mainone}]
We proceed using an iterative argument as in \cite[Theorem~8.1]{HuSi09}. Consider a mean curvature flow defined on $[0,t_0]$ which is either smooth or is a flow with surgery satisfying \ref{property S} at times prior to $t_0$. Let $t_0$ be the next time where the maximum of $|H|$ reaches $H_3$. The claim is that we can perform a finite number of surgeries on $\mathcal M_{t_0}$, which satisfy \ref{s0} and \ref{s1}-\ref{s3}, such that after these surgeries the maximum of $|H|$ is at most $H_2$ except on some components diffeomorphic to $\mathbb{S}^n$, $\mathbb{S}^{n-1} \times \mathbb{S}^1$ or $\mathbb{S}^{n-1} \simtimes \mathbb{S}^1$ which can be discarded. Let us therefore analyse the regions of our submanifold with large curvature. Let $p_0\in \mathcal{M}_{t_0}$ be any point at which $|H(p_0,t_0)|\geq H_2$. We have two cases to consider.

\textbf{Case (a)}. Suppose first that $|A(p_0,t_0)|^2-\frac{1}{n-1} |H(p_0,t_0)|^2\geq - \eta_1 |H(p_0,t_0)|^2$. We may then apply neck continuation theorem to obtain a neck $\mathcal N_0\ni p_0$ with the properties described there. We denote by $\mathcal A$ the region consisting of the neck $\mathcal N_0$, possibly together with the one or two regions diffeomorphic to standard $n$-balls which occur in case \ref{NC_convex cap} of the neck continuation theorem. The region $\mathcal A$ has one of the following possible structures:
\begin{itemize}
    \item $\mathcal A$ has two boundary components and is diffeomorphic to $\mathbb{S}^{n-1} \times [-1,1]$.
    \item $\mathcal A$ has one boundary component and is diffeomorphic to the standard $n$-ball.
    \item $\mathcal A$ has no boundary, is the entire connected component of $\mathcal M_{t_0}$ containing $p_0$, and is either diffeomorphic to $\mathbb{S}^n$ if it has two caps, or to one of the bundles  $\mathbb{S}^{n-1}\times\mathbb{S}^1$ or $\mathbb{S}^{n-1} \simtimes \mathbb{S}^1$ in case it is covered by the neck $\mathcal N_0$ (see Theorem~\ref{neck overlapping properties}). 
\end{itemize}
If $\partial \mathcal A\neq \emptyset$ then $\partial \mathcal A$ consists of either one or two cross-sections of $\mathcal N_0$ with mean radius equal to $2(n-1)/H_1$ and hence with mean curvature close to $H_1/2$.

\textbf{Case (b)}. Now suppose $|A(p_0,t_0)|^2-\frac{1}{n-1} |H(p_0,t_0)|^2< - \eta_1 |H(p_0,t_0)|^2$. If this inequality holds on the entire component of $\mathcal M_{t_0}$ containing $p_0$ then this component is diffeomorphic to $\mathbb{S}^{n}$ by \cite{Andrews2010}. Otherwise we can appeal to to Lemma~\ref{lem_pinching implies compact} as in \ref{P5} to see that there exists a $q_0$ with $d_{g(t_0)}(q_0, p_0) \leq \Upsilon/|H(p_0,t_0)|$ such that $ |A(q_0,t_0)|^2-\frac{1}{n-1} |H(q_0,t_0)|^2\geq - \eta_1 |H(q_0,t_0)|^2$. Moreover, for all $q$ satisfying $d_{g(t_0)}(q,p_0) \leq d_{g(t_0)}(q_0, p_0)$ we have
	\begin{align*}
	|H(q,t_0)| \geq \frac{|H(p_0,t_0)|}{\gamma_0} \geq \frac{H_2}{\gamma_0} = 10H_1.
	\end{align*}
In particular, we have $|H(q_0,t_0)| \geq 10 H_1$, so we can apply the neck continuation theorem at $q_0$ to obtain a region $\mathcal A \ni q_0$ as in Case~(a). Observe that $p_0$ is contained in $\mathcal A$;  if not then any minimizing geodesic from $p_0$ to $q_0$ would have to pass through $\partial\mathcal A$, and hence contain a point where the mean curvature is close to $H_1/2$, but by our choice of $q_0$ we know that $|H| \geq 10H_1$ at every point of such a geodesic. 

In both cases we have covered the point $p_0$ with a region $\mathcal A$ having the structure described above. Next, if there is a $p' \not \in \mathcal A$ with $|H(p',t_0)| \geq H_2$ then we repeat the procedure to obtain a region $\mathcal A' \ni p'$. The regions $\mathcal A$ and $\mathcal A'$ do not overlap---indeed, if $\partial \mathcal A$ is nonempty then it consists of one or two cross-sections of a neck with mean radius equal to $2(n-1)/H_1$, but in the neck continuation theorem if we meet such a cross section we stop there, so $\mathcal A$ and $\mathcal A'$ can only meet at boundary points. Repeating the procedure we obtain a collection of regions $\mathcal A, \mathcal A', \dots, \mathcal A^{(k)}$ covering all points of $\mathcal M_{t_0}$ where the mean curvature is at least $H_2$. There can only be finitely many of these because they do not overlap and each one has area bounded from below by a fixed multiple of $H_2^{-n}$. 

We now proceed with the surgery procedure. First, let us discard all of the regions $\mathcal A^{(i)}$ which have empty boundary. Next, consider any $\mathcal A^{(i)}$ with nonempty boundary. Each component of $\partial \mathcal A^{(i)}$ is a cross-section of a neck with mean radius $2(n-1)/H_1$. By construction, the neck contains a point where the curvature is at least $H_2 \geq 10 H_1$. Starting at this point and moving towards a boundary component, there exists a first cross-section $\Sigma^{(i)}$ such that the mean radius is $(n-1)/H_1$, namely half the mean radius of the boundary. We perform standard surgery at this cross-section $\Sigma^{(i)}$. If $\mathcal A^{(i)}$ has two boundary components then we apply this procedure on both sides. Whether we perform one or two surgeries on a region $\mathcal A^{(i)}$, we always create a connected component which is diffeomorphic to $\mathbb{S}^n$ and contains all points of $\mathcal A^{(i)}$ where the mean curvature is at least $H_2$. This component is discarded, so that after performing the surgeries the maximum of $|H|$ is less than $H_2$. 

To continue the iteration we need to show that the standard surgery at the cross-section $\Sigma^{(i)}$ can be performed such that \ref{s0} and \ref{s1}-\ref{s3} are satisfied. It is easy to check that, because of how we chose the cross-section $\Sigma^{(i)}$, properties \ref{s1}-\ref{s3} are satisfied with $\varepsilon_* = \varepsilon_0$, $K_* = H_1$ and $L_* = L$. In order to establish \ref{s0} we appeal to Lemma~\ref{lem_ss} and Lemma~\ref{lem_separation implies s0} as in \ref{P0}. Indeed, we know that $\Sigma^{(i)}$ lies at the centre of a normal $(\varepsilon_0, 1)$-cylindrical submanifold neck, in a cross-section of mean radius $r_* = (n-1)/H_1$, where the length of the neck is at least $5L \geq 5\hat L$. This is because the neck must be long enough for the mean radius to increase to $2(n-1)/H_1$ in one direction and for the mean curvature to increase to $H_2$ in the other. Appealing to Lemma~\ref{lem_ss} as described in \ref{P0} we see that $\hat {\mathcal P}(p,t_0, \hat L, \hat L^2)$ is an $(\hat \varepsilon,1)$-shrinking neck for every $p \in \Sigma^{(i)}$. By our choice of $\hat \varepsilon$ and $\hat L$, the planarity improvement theorem then ensures that \ref{s0} holds for the surgery performed at $\Sigma^{(i)}$, as we demonstrated in Lemma~\ref{lem_separation implies s0}. 

After performing surgeries at time $t_0$ as described above, we restart the flow. The next time the maximum of $|H|$ reaches $H_3$ we repeat the procedure. Since the mean curvature flow decreases area and each surgery removes a region of the submanifold whose area is a fixed multiple of $H_1^{-n}$, there can only be finitely many surgeries. That is, the flow with surgery terminates after finitely many steps. In particular, eventually all of the remaining components are recognised as being diffeomorphic to $\mathbb{S}^n$, $\mathbb{S}^{n-1}\times\mathbb{S}^1$ or $\mathbb{S}^{n-1}\simtimes\mathbb{S}^1$. 
\end{proof}

\begin{proof}[Proof of Corollary~\ref{cor_topological conclusion}]
When the flow with surgery terminates, we are left with a finite collection of submanifolds, each of which is diffeomorphic either to $\mathbb{S}^n$, $\mathbb{S}^{n-1}\times \mathbb{S}^1$ or $\mathbb{S}^{n-1} \simtimes \mathbb{S}^1$. This includes all components discarded at surgery times. Following the flow backwards from the final time, since each surgery is a reverse connected sum, we see that $\mathcal M$ is diffeomorphic either to $\mathbb{S}^n$ or a finite connected sum of copies of $\mathbb{S}^{n-1}\times\mathbb{S}^1$ and $\mathbb{S}^{n-1}\simtimes\mathbb{S}^1$. 
\end{proof}


\appendix

\section{Interior derivative estimates for $A^-$}\label{sec_interior estimates A-} Interior derivative estimates for mean curvature flows with bounded second fundamental form were proven in \cite{Ecker1991} and \cite{Andrews2010}. The estimates in \cite{Andrews2010} hold in every codimension but are only interior in time, not in space. We observe that they also localise in space. 

\begin{theorem}
Let $\mathcal M_t$ be a smooth $n$-dimensional mean curvature flow, properly immersed in $B(x_0,r) \subset \mathbb{R}^{n+m}$ and defined for all times $t \in [t_0 - r^2, t_0]$. We then have 
    \begin{equation}\label{eq interior estimates A}
    r^{2\ell+2}|\nabla^\ell A|^2(x_0,t_0) \leq C(n,\ell,K)
    \end{equation}
where $K := \sup_{B(x_0,r) \times [t_0 - r^2, t_0]} r^2 |A|^2$ 
\end{theorem}
\begin{proof}
The proof is easily adapted from \cite{Andrews2010} and \cite{Ecker1991}. 
\end{proof}

Whenever $|H|$ is nonzero we set $\nu = H / |H|$, $h = \langle A, \nu\rangle$ and $A^- = A - h\nu$ as usual. Our main result in this section says that if $A^-$ is small in a region of spacetime then all of its derivatives are small in the interior, provided we have $h \leq (1-\delta)|H|g$ for a constant $\delta > 0$.  

\begin{theorem}\label{interior estimates A-}
Let $\mathcal M_t$ be a smooth $n$-dimensional mean curvature flow, properly immersed in $B(x_0,r) \subset \mathbb{R}^{n+m}$ for all times $t \in [t_0 - r^2, t_0]$. Suppose $|H| > 0$ and $h \leq (1-\delta)|H|g$ hold everywhere, where $\delta > 0$ is a constant. We then have 
    \begin{equation*}
    r^{2\ell+2}|\nabla^\ell A^-|^2(x_0,t_0) \leq C(n,\ell,\delta, K,J_\ell) \sup_{B(x_0,r)\times[t_0 -r^2, t_0]} r^2|A^-|^2
    \end{equation*}
where $K := \sup_{B(x_0,r) \times [t_0 - r^2, t_0]} r^2 |A|^2$ and $J_\ell : =  \sup_{B(x_0,r) \times [t_0 - r^2, t_0]} \left(\sum_{i=0}^{\ell+1} |H|^{-2i-2}|\nabla^i A|^2\right)$. 
\end{theorem}

Using the evolution equation 
\begin{align*}
(\nabla_{t} - \Delta) h_{ij} &=  |h|^2 h_{ij}  + A^- \ast A^-\ast h -  h_{ij} |\nabla \nu|^2  - \frac{2}{|H|} \langle \nabla A_{ij}^-, \nabla H \rangle,
\end{align*}
which was derived in \cite{Lynch--Nguyen_convexity}, together with $(\nabla_{t}-\Delta) H =\langle H, A_{ij}\rangle A_{ij}$ and the definition $A^- = A - h \nu$, it is not difficult to see that $A^-$ solves an equation of the form 
\begin{align*}
(\nabla_t - \Delta) A^- = \mathbf B \ast \nabla A^- + \mathbf C \ast A^- + \mathbf D\ast \nabla \nu
\end{align*}
where 
    \[\mathbf B = \frac{h}{|H|} \ast \frac{\nabla h}{|H|}, \quad \mathbf C = A^- \ast h \ast \frac{h}{|H|} + h \ast h + A^- \ast A^-, \quad \mathbf D = \frac{h}{|H|}  \ast \nabla h + \nabla h.\]
We use the commutation identities\footnote{These can be deduced from the abstract spacetime-bundle formalism in \cite{Andrews2010}, or by direct computation. For the latter approach one needs the formula $\nabla_t X^i = \partial_t X^i - \langle H, A_j^i\rangle X^j$.}
    \[\nabla \nabla_t T = \nabla_t \nabla T + A \ast \nabla A \ast T + A \ast A \ast \nabla T\]
and 
    \[\nabla \Delta T = \Delta \nabla T + A \ast \nabla A \ast T + A \ast A \ast \nabla T\]
to obtain 
    \[(\nabla_t - \Delta)\nabla^\ell A^- = \nabla^\ell(\mathbf B \ast \nabla A^- + \mathbf C \ast A^- + \mathbf D \ast \nabla \nu) + \sum_{i+j+k=\ell} \nabla^i A \ast \nabla^j A \ast \nabla^k A^-.\]
For all $0\leq i \leq \ell$ we have
    \[r^{i+1}|\nabla^i \mathbf B|  + r^{i+2}|\nabla^i \mathbf C| + r^{i+1}|H|^{-1}|\nabla^i \mathbf D| \leq C(n,\ell, K, J_\ell),\]
so this implies 
    \begin{align}\label{A- derivs evol}
    (\nabla_t - \Delta)|\nabla^\ell A^-|^2 &\leq -2|\nabla^{\ell+1}A^-|^2 +C\sum_{i=0}^{\ell+1} r^{-\ell-2+i}|\nabla^i A^-||\nabla^\ell A^-| \notag\\
    & \qquad + C\sum_{i=1}^{\ell+1} r^{-\ell-2+i}|H||\nabla^i \nu||\nabla^\ell A^-|
    \end{align}
for a constant $C = C(n,\ell, K, J_\ell)$.

\begin{proof}[Proof of Theorem~\ref{interior estimates A-}] First we employ the pinching assumption $h \leq (1-\delta)|H|g$ to rewrite the derivatives of $\nu$ in \eqref{A- derivs evol} in terms of derivatives of $A^-$. Since $|H|g - h$ is positive-definite we may define $P$ to be the $(1,1)$-tensor such that $P^k_j(|H|\delta_k^i - h_k^i)=\delta^i_j$ and so use \eqref{eq_principal torsion} to express
    \[\nabla \nu = P \ast \nabla A^- + P \ast \nabla A^- \ast \frac{h}{|H|} \ast \frac{h}{|H|}.\]
Combined with the pinching assumption this gives $|H||\nabla \nu| \leq C(n,\delta)|\nabla A^-|$ and, moreover,
    \[|\nabla^i \nu| \leq C(n,\ell,\delta,K,J_\ell)\sum_{j=1}^i r^{j-i}\frac{|\nabla^j A^-|}{|H|}\]
for each $1 \leq i \leq \ell+1$. Substituting this into \eqref{A- derivs evol}, we obtain
    \begin{equation}\label{evolution higher derivatives A- 2}
    (\nabla_t - \Delta)|\nabla^\ell A^-|^2 \leq -2|\nabla^{\ell+1}A^-|^2 +C\sum_{i=0}^{\ell+1} r^{-\ell-2+i}|\nabla^i A^-||\nabla^\ell A^-|
    \end{equation}
where $C = C(n,\ell, \delta, K, J_\ell)$. 

From here on the proof is very similar to that of the Ecker--Huisken estimates. We use induction to establish that, for each $\ell \in \mathbb{N}$,
    \begin{equation}\label{induction higher interior estimate A-}
    \sup_{B(x_0, r/2^{\ell}) \times [t_0 - r^2/2^{2\ell}, t_0]} r^{2\ell+2}|\nabla^\ell A^-|^2\leq C(n,\ell,\delta, K, J_\ell) \sup_{B(x_0, r) \times [t_0 -r^2, t_0]} r^2|A^-|^2,
    \end{equation}
which implies the theorem. The case $\ell = 0$ is trivial, so let us assume \eqref{induction higher interior estimate A-} holds for all $0 \leq \ell \leq p - 1$ and prove that it then also holds for $\ell = p$. We introduce a cutoff function
    \[\varphi(x,t) = \bigg(1-\frac{t_0 - t}{r^2/2^{2p-2}}\bigg)\bigg(1-\frac{|x-x_0|^2}{r^2/2^{2p-2}} \bigg)_+,\]
such that
    \[0 \leq \varphi \leq 1, \qquad r |\nabla \varphi| + r^2|(\partial_t -\Delta)\varphi| \leq C(n,p),\]
and $\varphi$ vanishes outside of $B(x_0,r/2^{p-1})$ and at time $t = t_0 - r^2/2^{2p-2}$. Setting $G_i := r^{2i+2}|\nabla^i A|^2$ and using \eqref{evolution higher derivatives A- 2} together with Young's inequality we see that 
    \begin{align*}
    (\partial_t - \Delta) G_i &\leq -2r^{-2}G_{i+1} + Cr^{-2}\sqrt{G_{i+1}}\sqrt{G_i} + Cr^{-2} \sum_{j=0}^i G_jn\leq -r^{-2}G_{i+1} + Cr^{-2} \sum_{j=0}^i G_j.
    \end{align*}
for each $i \leq p$, where $C = C(n,p,\delta,K,J_p)$. 

We complete the induction by applying the maximum principle to $Z := \varphi^2 G_p + a G_{p-1}$ where $a > 0$ is a constant to be fixed in a moment. We have
    \begin{align*}
    (\partial_t - \Delta) Z &\leq (\partial_t - \Delta)\varphi^2 G_p  - \varphi^2 r^{-2} G_{p+1} + C\varphi^2 r^{-2} (G_p + \dots + G_0)\\
    &\qquad + Cr^{-1} \varphi |\nabla \varphi| \sqrt{G_p}\sqrt{G_{p+1}} -a r^{-2} G_p + C a r^{-2}(G_{p-1} + \dots + G_0)
    \end{align*}
for some $C = C(n,p,\delta, K, J_p)$. Using Young's inequality, and assuming $a$ is fixed large enough (depending only on $n$, $p$, $\delta$, $K$ and $J_p$) so that 
    \[r^2(\partial_t - \Delta)\varphi^2 + Cr^2 |\nabla \varphi|^2 + C\varphi^2 - a \leq - a/2,\]
this implies
    \begin{align*}
    (\partial_t - \Delta) Z &\leq -\frac{a}{2}r^{-2}G_p +C(1+a) r^{-2}(G_{p-1} + \dots + G_0).
    \end{align*}
Consider a point $(\bar x, \bar t)$ where $Z$ attains its maximum over $B(x_0,r/2^{p-1}) \times [t_0 - r^2/2^{2p-2}, t_0]$. This must be an interior maximum unless $Z$ vanishes identically, so $(\partial_t - \Delta)Z(\bar x, \bar t) \geq 0$ and hence 
    \[G_p(\bar x, \bar t) \leq C (G_{p-1}(\bar x, \bar t) + \dots + G_0(\bar x, \bar t)).\]
Using the definition of $Z$ and the inductive hypothesis we deduce 
    \[Z(\bar x, \bar t) \leq C \sup_{B(x_0,r)\times[t_0 - r^2,t_0]} r^2|A^-|^2\] 
for some $C = C(n,p,\delta, K, J_p)$. To complete the induction we observe that $\varphi \geq 1/C$ and hence $r^{2p+2}|\nabla^p A|^2 \leq C Z(\bar x, \bar t)$ in $B(x_0,r/2^p) \times [t_0 - r^2/2^{2p}, t_0]$.
\end{proof}

Assuming a uniform lower bound for the mean curvature, we obtain the following corollary of Theorem~\ref{interior estimates A-}. 

\begin{corollary}\label{interior estimates A- H bdd below}
Let $\mathcal M_t$ be a smooth $n$-dimensional mean curvature flow, properly immersed in $B(x_0,r) \subset \mathbb{R}^{n+m}$ for all times $t \in [t_0 - r^2, t_0]$. Suppose the inequalities $|H| \geq \alpha r^{-1}$ and $h \leq (1-\delta)|H|g$ hold everywhere for some constants $\alpha > 0$ and $\delta >0$. We then have 
    \begin{equation*}
    r^{2\ell+2}|\nabla^\ell A^-|^2(x_0,t_0) \leq C(n,\ell, \alpha, \delta, K) \sup_{B(x_0,r)\times[t_0 -r^2, t_0]} r^2|A^-|^2
    \end{equation*}
where $K := \sup_{B(x_0,r) \times [t_0 - r^2, t_0]} r^2 |A|^2$. 
\end{corollary}
\begin{proof}
Because of the Ecker--Huisken estimates, we can bound $J_\ell \leq C(n,\ell, \alpha, K)$ in the region $B(x_0, r/2)\times[t_0-r^2/4, t_0]$, and then apply Theorem~\ref{interior estimates A-} in this region. 
\end{proof}

\bibliography{surgery.bib}
\bibliographystyle{alpha}
\end{document}